\documentclass{amsart}
\usepackage[doi=false, isbn=false, url=false, style=alphabetic, backend=bibtex]{biblatex}
\AtEveryBibitem{\clearfield{note}} 
\AtEveryBibitem{\clearfield{MRCLASS}} 
\AtEveryBibitem{\clearfield{MRNUMBER}} 
\AtEveryBibitem{\clearfield{MRREVIEWER}}  
\usepackage{accents}
\usepackage{graphicx}
\usepackage{amsfonts}
\usepackage{amscd}
\usepackage{amssymb}
\usepackage{tikz}
\usepackage{calligra}
\usepackage{extpfeil}
\usetikzlibrary{shapes}
\usetikzlibrary{decorations.markings}
\usepackage{tikz-cd}
\usetikzlibrary{arrows}
\usepackage{multirow}
\usepackage[savepos]{zref}
\newtheorem{thm}{Theorem}[section]
\newtheorem{corollary}[thm]{Corollary}
\newtheorem{lemma}[thm]{Lemma}
\newtheorem{proposition}[thm]{Proposition}

\newtheorem{example}[thm]{Example}

\theoremstyle{definition}
\newtheorem{definition}[thm]{Definition}
\theoremstyle{remark}
\newtheorem{remark}[thm]{Remark}
\newtheorem{assumption}[thm]{Assumption}
\numberwithin{equation}{section}
\numberwithin{figure}{section}
\usepackage[mathscr]{euscript}
\usepackage[T2A,T1]{fontenc}

\newenvironment{notedescription}%
  {\begin{description}%
    \setlength{\itemsep}{2.5pt}%
    \setlength{\parskip}{5pt}}%
  {\end{description}}

\DeclareMathOperator{\GCD}{GCD}
\DeclareMathOperator{\Fr}{{Fr}}

\DeclareMathOperator{\prof}{prof}
\DeclareMathOperator{\mult}{m}
\DeclareMathOperator{\HW}{HW}

\DeclareMathOperator{\SC}{sc}
\DeclareMathOperator{\der}{der}

\DeclareMathOperator{\Comm}{Comm}

\DeclareMathOperator{\Br}{Br}

\DeclareMathOperator{\cent}{cent}

\DeclareMathOperator{\Hom}{Hom}

\DeclareMathOperator{\Par}{Par}

\DeclareMathOperator{\Aut}{Aut}
\DeclareMathOperator{\Int}{Int}

\DeclareMathOperator{\ab}{ab}

\DeclareMathOperator{\Gal}{Gal}
\DeclareMathOperator{\Span}{Span}

\DeclareMathOperator{\Ker}{Ker}

\DeclareMathOperator{\val}{val}

\DeclareMathOperator{\Hilb}{Hilb}

\DeclareMathOperator{\rec}{rec}

\DeclareMathOperator{\Spec}{Spec}

\DeclareMathOperator{\Id}{Id}

\DeclareMathOperator{\Tr}{Tr}
\DeclareMathOperator{\Norm}{N}
\DeclareMathOperator{\Kum}{\varkappa}
\DeclareMathOperator{\Obj}{Ob}

\DeclareMathOperator{\et}{\acute{e}t}

\newcommand{\From}{\colon}
\newcommand{\inar}{\ar@{^{(}->}}
\newcommand{\onar}{\ar@{->>}}

\newcommand{\defined}[1]{\underline{{#1}}}
\usepackage{accents}
\newlength{\dtildeheight}

\newcommand{\Baer}{\dotplus}

\newcommand{\Weil}{\mathcal{W}}
\newcommand{\weil}{\mathbf{w}}

\newcommand{\Inertia}{\mathcal{I}}
\newcommand{\Frob}{\Fr}

\newcommand{\limdir}{\varinjlim}

\makeatletter
\newcommand{\raisemath}[1]{\mathpalette{\raisem@th{#1}}}
\newcommand{\raisem@th}[3]{\raisebox{#1}{$#2#3$}}
\makeatother
\newcommand{\Vee}{ { \raisemath{-3pt}{\vee} } }

\newcommand{\EM}{{\mathsf M}}
\newcommand{\EL}{{\mathsf L}}

\newcommand{\Cat}[1]{ {\mathsf{#1}} }
\newcommand{\Fun}[1]{ {\mathsf{#1}} }

\newcommand{\Lie}[1]{ {\mathfrak{#1}} }

\newcommand{\sch}[1]{\underline{\boldsymbol{ \mathrm{#1}}}}
\newcommand{\alg}[1]{\boldsymbol{\mathrm{#1}}}

\newcommand{\sheaf}[1]{{\mathscr{#1}}}
\newcommand{\shom}{\mathscr{H}\mathit{om}}
\newcommand{\stors}{\boldsymbol{\Cat{Tors}}}
\newcommand{\ssym}{\mathscr{S}\mathit{ym}}

\newcommand{\mGal}{\widetilde{\Gal}}
\newcommand{\amu}{\sch{\mu}}
\newcommand{\sAut}{\mathscr{A}\mathit{ut}}

\newcommand{\sspl}{\mathscr{S}\mathit{pl}}
\newcommand{\Whit}{\mathscr{W}\mathit{hit}}
\newcommand{\whit}{\boldsymbol{w}}
\newcommand{\dgp}[1]{\boldsymbol{\mathscr{#1}}}
\newcommand{\gerb}[1]{\boldsymbol{\Cat{#1}}}
\newcommand{\gerbob}[1]{{\boldsymbol{#1}}}

\newcommand{\ZZ}{\mathbb Z}
\newcommand{\QQ}{\mathbb Q}
\newcommand{\RR}{\mathbb R}

\newcommand{\CC}{\mathbb C}

\renewcommand{\AA}{\mathbb A}
\newcommand{\VV}{\mathcal V}

\newcommand{\kk}{\mathfrak{f}}

\newcommand{\ident}{\equiv}

\newcommand{\FF}{\mathbb F}
\newcommand{\OO}{\mathcal{O}}

\newcommand{\Into}{\hookrightarrow}
\newcommand{\Onto}{\twoheadrightarrow}
\newcommand{\To}{\rightarrow}

\renewcommand{\H}{{\mathsf{H}}}

\newcommand{\isom}{\cong}
\newcommand{\half}{\tfrac{1}{2}}

\renewcommand{\th}{\text{th}}

\newcommand{\inarrow}{\arrow[hook]}
\newcommand{\onarrow}{\arrow[two heads]}

\makeatletter
\newcommand\@biprod[1]{%
  \vcenter{\hbox{\ooalign{$#1\prod$\cr$#1\coprod$\cr}}}}
\newcommand\biprod{\mathop{\mathpalette\@biprod\relax}\displaylimits}
\makeatother

\newcommand{\defeq}{:=}
\DeclareMathAlphabet{\mathcalligra}{T1}{calligra}{m}{n}
\DeclareMathOperator{\Zar}{Zar}

\begin{document}

\title{The L-group of a covering group}%
\author{Martin H. Weissman}%
\date{\today}

\address{Yale-NUS College, 6 College Ave East, \#B1-01, Singapore 138614}
\email{marty.weissman@yale-nus.edu.sg}%

\subjclass[2010]{11F70; 22E50; 22E55.}

% \subjclass{}%
% \keywords{}%
%\date{}%
%\dedicatory{}%
%\commby{}%
% ----------------------------------------------------------------
\begin{abstract}
We incorporate nonlinear covers of quasisplit reductive groups into the Langlands program, defining an L-group associated to such a cover.  This L-group is an extension of the absolute Galois group of a local or global field $F$ by a complex reductive group.  The L-group depends on an extension of a quasisplit reductive $F$-group by $\alg{K}_2$, a positive integer $n$ (the degree of the cover), an injective character $\epsilon \From \mu_n \Into \CC^\times$, and a separable closure of $F$.  Our L-group is consistent with previous work on covering groups, and its construction is contravariantly functorial for certain ``well-aligned'' homomorphisms.  An appendix surveys torsors and gerbes on the \'etale site, as they are used in a crucial step in the construction.
\end{abstract}

\maketitle

\tableofcontents

\section*{Introduction}

Let $\alg{G}$ be a quasisplit reductive group over a local or global field $F$.  Let $n$ be a positive integer and let $\mu_n$ denote the group of $n^{\th}$ roots of unity in $F$.  Assume that $\mu_n$ has order $n$.  Let $\alg{G}'$ be a central extension of $\alg{G}$ by $\alg{K}_2$, in the sense of Brylinski and Deligne \cite{B-D}.  We call the pair $\alg{\tilde G} \defeq (\alg{G}', n)$ a ``degree $n$ cover'' of $\alg{G}$.  Fix a separable closure $\bar F / F$ and write $\Gal_F = \Gal(\bar F / F)$.  Fix an injective character $\epsilon \From \mu_n \Into \CC^\times$.

The purpose of this article is to provide a robust definition of the L-group ${}^\EL \tilde G$ associated to this data.  It will be an extension,
$$\tilde G^\vee \Into {}^\EL \tilde G \Onto \Gal_F,$$
where $\tilde G^\vee$ is a complex reductive group.  The reader may replace $\CC$ by any $\ZZ[1/n]$-algebra $\Omega$ endowed with $\epsilon \From \mu_n \Into \Omega^\times$, if desired.

We have accumulated evidence that this L-group can be used to parameterize irreducible $\epsilon$-genuine representations (admissible in the local case, automorphic in the global case), in much the same way as the L-group is expected to parameterize representations in the linear case.  Our evidence, including cases of split tori, unramified representations (when $\alg{G}$ is unramified and $n$ is coprime to the residual characteristic), and double-covers of semisimple groups over $\RR$, will be forthcoming.  Further evidence is provided by the large literature on covering groups:  isomorphisms of Hecke algebras, theta correspondences, character lifting, and more, which may be viewed as examples of functoriality if one accepts our proposed L-group.

We devote this entire article to the definition of the L-group ${}^\EL \tilde G$, since it requires care to make sure that it is ``well-defined'' by the data given (independent of other choices like Borel subgroup and maximal torus).  We strive to give the cleanest definition possible and provide examples along the way.  A crucial step involves the construction of a gerbe $\gerb{E}_\epsilon(\alg{\tilde G})$ on $F_{\et}$; the appendix is meant to provide a suitable introduction to gerbes of this sort.  Another section defines the ``metaGalois group,''  a canonical extension $\mu_2 \Into \mGal_F \Onto \Gal_F$.  The metaGalois group may of intrinsic arithmetic interest, and provides hints of a motivic connection.  We do \textbf{not} assume that $F$ contains a primitive $(2n)^{\th}$ root of unity; such an assumption negates the need for the metaGalois group.

Our previous article \cite{MWCrelle} was limited to split reductive groups, and may be viewed as a proof-of-concept by Hopf algebra bludgening.  The ``dual group'' $\tilde G^\vee$ here is the same as that from \cite{MWCrelle}.  The ``first twist'' of \cite{MWCrelle} is unchanged in substance, but here it is encoded in a 2-torsion element of the center of $\tilde G^\vee$ and the metaGalois group.  The ``second twist'' of \cite{MWCrelle} provided the greatest challenge.  Among many reformulations that succeed or split groups, we found the gerbe $\gerb{E}_\epsilon(\alg{\tilde G})$, which extends to quasisplit groups.  

\subsection*{Philosophies}

A few principles are helpful when considering any putative Langlands program for covering groups.

\begin{enumerate}
\item
There is no $\epsilon$-genuine trivial (or Steinberg) representation for general covers, and so one should not expect a single distinguished splitting of the L-group.
\item
If some set of things is parameterized by cohomology in degree $2$, then that set of things should be viewed as the set of objects in a 2-category.
\item
Things which ``are trivial'' (e.g., extensions, gerbes) can be isomorphic to trivial things in interesting ways.
\end{enumerate}

\subsection*{Acknowledgments}

The ideas of this paper have evolved over the past few years, and I am very grateful for the numerous mathematicians who discussed the constructions in various stages of completeness and correctness.  The American Institute of Mathematics hosted a conference at which I spoke with Wee Teck Gan, Gordan Savin, Jeffrey Adams, Sergey Lysenko, Tamotsu Ikeda, Kaoru Hiraga, Tasho Kaletha and others.  Their previous work, and our discussions at the conference and elsewhere, have been very helpful.  I also appreciate the support of Harvard University during a short visit, where I gained insight from discussions with Dennis Gaitsgory, Dick Gross, and John Tate.  During a visit to the University of Michigan, I gained from feedback from Stephen DeBacker, and learned about gerbes from James Milne.

In \cite{GanGao} and \cite{GaoThesis}, Wee Teck Gan and Fan Gao have tested some of the conjectures of this paper, and they have gone further in developing the Langlands program for covering groups.  I have greatly appreciated our frequent conversations.  Their results provided constraints which kept the constructions of this paper on track. 

Pierre Deligne has kindly corresponded with me over the past few years, and his generosity has been incredibly helpful.  His ideas led me to a deeper understanding of the crucial questions, and his correspondence motivated me to pursue this project further.  

%\newpage
\section*{Notation}
\begin{notedescription}
\item[$F$]  A field, typically local or global.
\item[$\bar F$]  A separable closure of $F$.
\item[$\OO$] The ring of integers in $F$, in the nonarchimedean local case.
\item[$\AA$]  The ring of adeles of $F$, in the global case.
\item[${\Fr}$]  The geometric Frobenius automorphism.
\item[$S$]  A connected scheme, typically $\Spec(F)$ or $\Spec(\OO)$.
\item[$\bar s$]  The geometric point of $S$ corresponding to $\bar F$.
\item[$\Gal_S$]  The absolute Galois group $\pi_1^{\et}(S, \bar s)$.

\item[$\alg{X}$]  An algebraic variety over $S$, or sheaf on $S_{\Zar}$.
\item[$X$ or $X_F$]  The $F$-points $\alg{X}(F)$ for such a variety.

\item[$\sch{X}$]  A scheme over $\ZZ$.
\item[$\sch{G}_m$]  The multiplicative group.
\item[$\amu_n$] The group scheme over $\ZZ$ of $n^{\th}$ roots of unity.

\item[$\mu_n$]  The group $\amu_n(S)$, assumed to be cyclic of order $n$.

\item[$\sheaf{S}$]  A sheaf on $S_{\et}$.
\item[$\sheaf{S}{[U]}$] The sections of $\sheaf{S}$ over $U$ ($U \To S$ \'etale).
\item[$\dgp{G}$]  A local system on $S_{\et}$, of group schemes over $\ZZ$.
\item[$\Cat{C}$]  A category, with objects $\Obj(\Cat{C})$.
\item[$\gerb{E}$] A gerbe on $S_{\et}$.
\item[$\gerb{E}{[U]}$] The groupoid of sections of $\gerb{E}$ over $U$.

\item[$A$]  An abelian group.
\item[$A_{[n]}$]  Its $n$-torsion subgroup.
\item[$A_{/n}$]  The quotient $A / n A$.

\end{notedescription}

\newpage

\section{Covering groups}

Throughout this article, $S$ will be a scheme in one of the following two classes:  $S = \Spec(F)$ for a field $F$, or $S = \Spec(\OO)$ for a discrete valuation ring $\OO$ with fraction field $F$.  In the latter case, we assume that $\OO$ contains a field, or that $\OO$ has finite residue field.  We will often fix a positive integer $n$, and we will assume that $\mu_n = \amu_n(S)$ is a cyclic group of order $n$.  In Section \ref{MetaGaloisSection}, we will place further restrictions on $S$.

%Throughout this article, $S$ will be a scheme in one of the following three classes.
%\begin{description}
%\item[Global] $S = \Spec(F)$ for a global field $F$.
%\item[Local] $S = \Spec(F)$ for a local field $F$.
%\item[Local integral] $S = \Spec(\OO)$ for the ring of integers $\OO$ in a nonarchimedean local field $F$.
%\end{description}

\subsection{Reductive groups}

Let $\alg{G}$ be a reductive group over $S$.  We follow \cite{SGA3} in our conventions, so this means that $\alg{G}$ is a smooth group scheme over $S$ such that $\alg{G}_{\bar s}$ is a connected reductive group for all geometric points $\bar s$ of $S$.  Assume moreover that $\alg{G}$ is \textbf{quasisplit} over $S$.  

Let $\alg{A}$ be a maximal $S$-split torus in $\alg{G}$, and let $\alg{T}$ be the centralizer of $\alg{A}$ in $\alg{G}$.  Then $\alg{T}$ is a maximal torus in $\alg{G}$, and we say that $\alg{T}$ is a \defined{maximally split} maximal torus.  Let $\sheaf{X}$ and $\sheaf{Y}$ be the local systems (on $S_{\et}$) of characters and cocharacters of $\alg{T}$.  

Let $\alg{N}$ be the normalizer of $\alg{T}$ in $\alg{G}$.  Let $\sheaf{W}$ denote the Weyl group of the pair $(\alg{G}, \alg{T})$, viewed as a sheaf on $S_{\et}$ of finite groups.  Then $\sheaf{W}[S] = \alg{N}(S) / \alg{T}(S)$ (see \cite[Expos\'e XXVI, 7.1]{SGA3}).  Let $\alg{B}$ be a Borel subgroup of $\alg{G}$ containing $\alg{T}$, defined over $S$.  Let $\alg{U}$ be the unipotent radical of $\alg{B}$.
\begin{proposition}
\label{BTConj}
Assume as above that $\alg{G}$ is quasisplit, and $S$ is the spectrum of a field or of a DVR.  The group $\alg{G}(S)$ acts transitively, by conjugation, on the set of pairs $(\alg{B}, \alg{T})$ consisting of a Borel subgroup (defined over $S$) and a maximally split maximal torus therein.
\end{proposition}
\proof
As we work over a local base scheme $S$, \cite[Expos\'e XXVI, Proposition 6.16]{SGA3} states that the group $\alg{G}(S)$ acts transitively on the set of maximal split subtori of $\alg{G}$ (defined over $S$).

Every maximally split maximal torus of $\alg{G}$ is the centralizer of such a maximal split torus, and thus $\alg{G}(S)$ acts transitively on the set of maximally split maximal tori in $ \alg{G}$.  The stabilizer of such a maximally split maximal torus $\alg{T}$ is the group of $S$-points of its normalizer $\alg{N}(S)$.  The Weyl group $\sheaf{W}[S] = \alg{N}(S) / \alg{T}(S)$ acts simply-transitively on the minimal parabolic subgroups containing $\alg{T}$ by \cite[Expos\'e XXVI, Proposition 7.2]{SGA3}. This proves the proposition.  
\qed

The roots and coroots (for the adjoint action of $\alg{T}$ on the Lie algebra of $\alg{G}$) form local systems $\Phi$ and $\Phi^\vee$ on $S_{\et}$, contained in $\sheaf{X}$ and $\sheaf{Y}$, respectively.  The simple roots (with respect to the Borel subgroup $\alg{B}$) and their coroots form local systems of subsets $\Delta \subset \Phi$ and $\Delta^\vee \subset \Phi^\vee$, respectively.  In this way we find a local system on $S_{\et}$ of based root data (cf. \cite[\S 1.2]{BorelCorvallis}),
$$\Psi = \left( \sheaf{X}, \Phi, \Delta, \sheaf{Y}, \Phi^\vee, \Delta^\vee \right).$$
Write $\sheaf{Y}^{\SC}$ for the subgroup of $\sheaf{Y}$ spanned by the coroots.  

\subsection{Covers}

In \cite{B-D}, Brylinski and Deligne study central extensions of $\alg{G}$ by $\alg{K}_2$,  where $\alg{G}$ and $\alg{K}_2$ are viewed as sheaves of groups on the big Zariski site $S_{\Zar}$.  These extensions form a category we call $\Cat{CExt}_S(\alg{G}, \alg{K}_2)$,.  Such a central extension will be written $\alg{K}_2 \Into \alg{G}' \Onto \alg{G}$ in what follows.  We add one more piece of data in the definition below.
\begin{definition}
A degree $n$ \defined{cover} of $\alg{G}$ over $S$ is a pair $\alg{\tilde G} = (\alg{G}', n)$, where
\begin{enumerate}
\item
$\alg{K}_2 \Into \alg{G}' \Onto \alg{G}$ is a central extension of $\alg{G}$ by $\alg{K}_2$ on $S_{\Zar}$;
\item
$n$ is a positive integer;
\item
For all scheme-theoretic points $s \in S$, with residue field $\kk(s)$, $\# \amu_n( \kk(s) ) = n$.
\end{enumerate}
\end{definition}

Define $\Cat{Cov}_n(\alg{G})$ (or $\Cat{Cov}_{n/S}(\alg{G})$ to avoid confusion) to be the category of degree $n$ covers of $\alg{G}$ over $S$.  The objects are pairs $\alg{\tilde G} = (\alg{G}', n)$ as above, and morphisms are those from $\Cat{CExt}_S(\alg{G}, \alg{K}_2)$ (with $n$ fixed).

If $\gamma \From S_0 \To S$ is a morphism of schemes, then pulling back gives a functor $\gamma^\ast \From \Cat{Cov}_{n/S}(\alg{G}) \To \Cat{Cov}_{n/S_0}(\alg{G}_{S_0})$.  Indeed, a morphism of schemes gives inclusions of residue fields (in the opposite direction) and so Condition (3) is satisfied by the scheme $S_0$ when it is satisfied by the scheme $S$.

Central extensions $\alg{K}_2 \Into \alg{G}' \Onto \alg{G}$ are classified by a triple of invariants $(Q, \sheaf{D}, f)$.  For fields, this is carried out in \cite{B-D}, and the extension to DVRs (with finite residue field or containing a field) is found in \cite{MWIntegral}.  The first invariant $Q \From \sheaf{Y} \To \ZZ$ is a Galois-invariant Weyl-invariant quadratic form, i.e., $Q \in H_{\et}^0(S, \ssym^2(\sheaf{X})^{\sheaf{W}})$.  The second invariant $\sheaf{D}$ is a central extension of sheaves of groups on $S_{\et}$, $\sheaf{G}_m \Into \sheaf{D} \Onto \sheaf{Y}$.  The third invariant $f$ will be discussed later.  

A cover $\alg{\tilde G}$ yields a symmetric $\ZZ$-bilinear form $\beta_Q \From \sheaf{Y} \otimes_\ZZ \sheaf{Y} \To n^{-1} \ZZ$,
$$\beta_Q(y_1, y_2) \defeq n^{-1} \cdot \left( Q(y_1+y_2) - Q(y_1) - Q(y_2) \right).$$
This defines a local system $\sheaf{Y}_{Q,n} \subset \sheaf{Y}$,
$$\sheaf{Y}_{Q,n} = \{ y \in \sheaf{Y} : \beta_Q(y, y') \in \ZZ \text{ for all } y' \in \sheaf{Y} \}.$$

The category of covers is equipped with the structure of a Picard groupoid; one may ``add'' covers via the Baer sum.  If $\alg{\tilde G}_1, \alg{\tilde G}_2$ are two covers of $\alg{G}$ of degree $n$, one obtains a cover $\alg{\tilde G}_1 \Baer \alg{\tilde G}_2 = (\alg{G}_1' \Baer \alg{G}_2', n)$.

When $\alg{\tilde G} = (\alg{G}', n)$ is a degree $n$ cover of $\alg{G}$, and $\alg{H} \subset \alg{G}$ is a reductive subgroup defined over $S$, write $\alg{\tilde H}$ for the resulting cover of $\alg{H}$.  Following \cite{B-D}, the $\alg{G}'$ splits uniquely over any unipotent subgroup $\alg{U} \subset \alg{G}$, so we view $\alg{U}$ as a subgroup of $\alg{G}'$ in this case.

In three arithmetic contexts, a cover $\alg{\tilde G}$ yields a central extension of topological groups according to \cite[\S 10.3, 10.4]{B-D}.
\begin{description}
\item[Global] If $S = \Spec(F)$ for a global field $F$, then $\alg{\tilde G}$ yields a central extension $\mu_n \Into \tilde G_\AA \Onto G_\AA$, endowed with a splitting $\sigma_F \From G_F \Into \tilde G_\AA$.
\item[Local] If $S = \Spec(F)$ for a local field $F$, then $\alg{\tilde G}$ yields a central extension $\mu_n \Into \tilde G \Onto G$, where $G = \alg{G}(F)$.
\item[Local integral]  If $S = \Spec(\OO)$, with $\OO$ the ring of integers in a nonarchimedean local field $F$, then $\alg{\tilde G}$ yields a central extension $\mu_n \Into \tilde G \Onto G$, where $G = \alg{G}(F)$, endowed with a splitting $\sigma^\circ \From G^\circ \Into G$.
\end{description} 

Fix an injective character $\epsilon \From \mu_n \Into \CC^\times$.  This allows one to define $\epsilon$-genuine automorphic representations of $\tilde G_\AA$ in the global context, $\epsilon$-genuine admissible representations of $\tilde G$ in the local context, and $\epsilon$-genuine $G^\circ$-spherical representations of $\tilde G$ in the local integral context.

The purpose of this article is the construction an \defined{L-group} associated to such a $\alg{\tilde G}$ and $\epsilon$.  We believe that this L-group will provide a parameterization of irreducible $\epsilon$-genuine representations in the three contexts above.

\subsection{Well-aligned homomorphisms}

Let $\alg{G}_1 \supset \alg{B}_1 \supset \alg{T}_1$ and $\alg{G}_2 \supset \alg{B}_2 \supset \alg{T}_2$ be quasisplit groups over $S$, endowed with Borel subgroups and maximally split maximal tori.    Let $\alg{\tilde G}_1 = (\alg{G}_1', n)$ and $\alg{\tilde G}_2 = (\alg{G}_2', n)$ be covers (of the same degree) of $\alg{G}_1$ and $\alg{G}_2$, respectively.  Write $\sheaf{Y}_1$ and $\sheaf{Y}_2$ for the cocharacter lattices of $\alg{T}_1$ and $\alg{T}_2$, and $Q_1, Q_2$ for the quadratic forms arising from the covers.  These quadratic forms yield sublattices $\sheaf{Y}_{1,Q_1,n}$ and $\sheaf{Y}_{2,Q_2,n}$.  

\begin{definition}
A \defined{well-aligned homomorphism} $\tilde \iota$ from $(\alg{\tilde G}_1, \alg{B}_1, \alg{T}_1)$ to $(\alg{\tilde G}_2, \alg{B}_2, \alg{T}_2)$ is a pair $\tilde \iota = (\iota, \iota')$ of homomorphisms of sheaves of groups on $S_{\Zar}$, making the following diagram commute,
\begin{equation}
\label{CDiota}
\begin{tikzcd}
\alg{K}_2 \inarrow{r} \arrow{d}{=} & \alg{G}_1' \onarrow{r} \arrow{d}{\iota'} & \alg{G}_1 \arrow{d}{\iota} \\
\alg{K}_2 \inarrow{r} & \alg{G}_2' \onarrow{r} & \alg{G}_2
\end{tikzcd}
\end{equation}
and satisfying the following additional axioms:
\begin{enumerate}
\item
$\iota$ has normal image and smooth central kernel;
\item
$\iota(\alg{B}_1) \subset \alg{B}_2$ and $\iota(\alg{T}_1) \subset \alg{T}_2$.  Thus $\iota$ induces a map $\iota \From \sheaf{Y}_1 \To \sheaf{Y}_2$;
\item
$(\iota, \iota')$ realizes $\alg{G}_1'$ as the pullback of $\alg{G}_2'$ via $\iota$;
\item
The homomorphism $\iota$ satisfies $\iota(\sheaf{Y}_{1,Q_1,n}) \subset \sheaf{Y}_{2,Q_2,n}$.
\end{enumerate}
\end{definition}

\begin{remark}
Conditions (1) and (2) are inspired by \cite[\S 1.4, 2.1,2.5]{BorelCorvallis}, though more restrictive.  By ``normal image,'' we mean that for any geometric point $\bar s \To S$, the homomorphism $\iota \From \alg{G}_{1,\bar s} \To \alg{G}_{2, \bar s}$ has normal image.  Condition (3) implies that for all $y \in \sheaf{Y}_1$, $Q_1(y) = Q_2(\iota(y))$.  In other words, $Q_1$ is the image of $Q_2$ via the map
$$\iota^\ast \From H_{\et}^0(S, \ssym^2 (\sheaf{X}_2) ) \To H_{\et}^0(S, \ssym^2 (\sheaf{X}_1) ).$$
But Condition (3) does not imply Condition (4); one may cook up an example with $\alg{G}_1 = \alg{G}_{\mult}$ and $\alg{G}_2 = \alg{G}_{\mult}^2$ which satisfies (3) but not (4).
\end{remark}

\begin{proposition}
\label{ComposeWellaligned}
The composition of well-aligned homomorphisms is well-aligned.
\end{proposition}
\proof
Suppose that $(\iota_1, \iota_1')$ and $(\iota_2, \iota_2')$ are well-aligned homomorphisms, with $\iota_1 \From \alg{G}_1 \To \alg{G}_2$ and $\iota_2 \From \alg{G}_2 \To \alg{G}_3$.  Conditions (2), (3), and (4) are obviously satisfied by the composition $(\iota_2 \circ \iota_1, \iota_2' \circ \iota_1')$.  For condition (1), notice that the kernel of $\iota_2 \circ \iota_1$ is contained in the kernel of $\iota_1$, and hence is central.  The only thing left is to verify that $\iota_2 \circ \iota_1$ has normal image.  This may be checked by looking at geometric fibres; it seems well-known (cf. \cite[\S 1.8]{KotSTF}).
\qed

Inner automorphisms are well-aligned homomorphisms.
\begin{example}
Suppose that $\alg{\tilde G}$ is a degree $n$ cover of a quasisplit group $\alg{G}$.  Suppose that $\alg{B}_0 \supset \alg{T}_0$ and $\alg{B} \supset \alg{T}$ are two Borel subgroups containing maximally split maximal tori.  Suppose that $g \in \alg{G}(S)$, and write $\Int(g)$ for the resulting inner automorphism of $\alg{G}$.  As noted in \cite[0.N.4]{B-D}, $\Int(g)$ lifts canonically to an automorphism $\Int(g)' \in \Aut(\alg{G}')$.  If $\alg{B} = \Int(g) \alg{B}_0$ and $\alg{T} = \Int(g) \alg{T}_0$, then the pair $\left( \Int(g), \Int(g)' \right)$ is a well-aligned homomorphism from $(\alg{\tilde G}, \alg{B}_0, \alg{T}_0)$ to $(\alg{\tilde G}, \alg{B}, \alg{T})$.
\end{example}

While we focus on quasisplit groups in this article, the lifting of inner automorphisms allows one to consider ``pure inner forms'' of covers over a field.  
\begin{definition}
Let $\alg{\tilde G} = (\alg{G}', n)$ be a degree $n$ cover of a quasisplit group $\alg{G}$, over a field $F$.  Let $\xi \in Z_{\et}^1(F, \alg{G})$ be a 1-cocycle.  The image $\Int(\xi)$ in $Z_{\et}^1(F, \alg{Aut}(\alg{G}))$ defines an inner form $\alg{G}_\xi$ of $\alg{G}$.  These are called the \defined{pure inner forms} of $\alg{G}$.  On the other hand, we may consider the image $\Int(\xi)'$ in $Z_{\et}^1(F, \alg{Aut}(\alg{G}'))$, which by \cite[\S 7.1, 7.2]{B-D} defines a central extension $\alg{G}_\xi'$ of $\alg{G}_\xi$ by $\alg{K}_2$.  The cover $\alg{\tilde G}_\xi = (\alg{G}_\xi', n)$ of $\alg{G}_\xi$ will be called a \defined{pure inner form} of the cover $\alg{\tilde G}$.  
\end{definition}  

We have not attempted to go further in the study of inner forms for covers, but presumably one should study something like strong real forms as in \cite[Definition 1.12]{AdamsBarbaschVogan}, and more general rigid forms as in \cite{Kaletha}, if one wishes to assemble L-packets for covering groups.

The next example of a well-aligned homomorphism is relevant for the study of central characters of genuine representations.
\begin{example}
Let $\alg{\tilde G}$ be a degree $n$ cover of a quasisplit group $\alg{G} \supset \alg{B} \supset \alg{T}$.  Let $\alg{H}$ be the maximal torus in the center of $\alg{G}$, with cocharacter lattice $\sheaf{Y}_H \subset \sheaf{Y}$.  Let $\alg{C}$ be the algebraic torus with cocharacter lattice $\sheaf{Y}_H \cap \sheaf{Y}_{Q,n}$, and $\iota \From \alg{C} \To \alg{G}$ the resulting homomorphism (with central image).  Let $\alg{\tilde C}$ denote the pullback of the cover $\alg{\tilde G}$ via $\iota$.  Then $\iota$ lifts to a well-aligned homomorphism from $\alg{\tilde C}$ to $\alg{\tilde G}$.
\end{example}

The final example of a well-aligned homomorphism is relevant for the study of twisting genuine representations by one-dimensional representations of $\alg{G}$.
\begin{example}
Let $\alg{H}$ denote the maximal toral quotient of $\alg{G}$.  In other words, $\alg{H}$ is the torus whose character lattice equals $\Hom(\alg{G}, \alg{G}_m)$.  Let $p \From \alg{G} \To \alg{H}$ denote the canonical homomorphism, and write $\iota \From \alg{G} \To \alg{G} \times \alg{H}$ for the homomorphism $\Id \times p$.  Write $\alg{\tilde G} \times \alg{H}$ for the cover $(\alg{G}' \times \alg{H}, n)$.  

The homomorphism $\iota$ realizes $\alg{\tilde G}$ as the pullback via $\iota$ of the cover $\alg{\tilde G} \times \alg{H}$.  A Borel subgroup and torus in $\alg{G}$ determines a Borel subgroup and torus in $\alg{G} \times \alg{H}$.  In this way, $\iota$ lifts a well-aligned homomorphism of covers from $\alg{\tilde G}$ to $\alg{\tilde G} \times \alg{H}$.  
\end{example}

\section{The dual group}

In this section, fix a degree $n$ cover $\alg{\tilde G}$ of a quasisplit group $\alg{G}$ over $S$.  Associated to $\alg{\tilde G}$, we define the ``dual group,'' a local system on $S_{\et}$ of affine group schemes over $\ZZ$.  We refer to Appendix \S \ref{GroupsTorsors}, for background on such local systems.  We begin by reviewing the Langlands dual group of $\alg{G}$ in a framework suggested by Deligne (personal communication).

\subsection{The Langlands dual group}

Choose, for now, a Borel subgroup $\alg{B} \subset \alg{G}$ containing a maximally split maximal torus $\alg{T}$.  The based root datum of $(\alg{G}, \alg{B}, \alg{T})$ was denoted $\Psi$, and the dual root datum,
$$\Psi^\vee = \left( \sheaf{Y}, \Phi^\vee, \Delta^\vee, \sheaf{X}, \Phi, \Delta \right),$$
is a local system of root data on $S_{\et}$.  

This defines a unique (up to unique isomorphism) local system $\dgp{G^\vee}$ on $S_{\et}$ of pinned reductive groups over $\ZZ$, called the \defined{Langlands dual group} of $\alg{G}$.  The center of $\dgp{G^\vee}$ is a local system on $S_{\et}$ of groups of multiplicative type over $\ZZ$, given by
$$\dgp{Z^\vee} = \Spec \left( \ZZ [ \sheaf{Y} / \sheaf{Y}^{\SC} ] \right).$$
See Example \ref{LocSpec} for more on local systems and $\Spec$ in this context.

\subsection{The dual group of a cover}

Now we adapt the definition of the dual group to covers.  The ideas here are the same as those of \cite{MWCrelle}.  The ideas for modifying root data originate in \cite[\S 2.2]{LusztigQuantumGroups} in the simply-connected case, in \cite[Theorem 2.9]{FinkelbergLysenko} for the almost simple case, in \cite[\S 11]{McNamara} and \cite{Reich} in the reductive case.  This dual group is also compatible with \cite{ABPTV} and the Hecke algebra isomorphisms of Savin \cite{SavinUnramified}, and the most recent work of Lysenko \cite{LysenkoSatake}.

Associated to the cover $\alg{\tilde G}$ of degree $n$, recall that $Q \From \sheaf{Y} \To \ZZ$ is the first Brylinski-Deligne invariant, and $\beta_Q \From \sheaf{Y} \otimes \sheaf{Y} \To n^{-1} \ZZ$ a symmetric bilinear form, and
$$\sheaf{Y}_{Q,n} = \{ y \in \sheaf{Y} : \beta_Q(y, y') \in \ZZ \text{ for all } y' \in \sheaf{Y} \} \subset \sheaf{Y}.$$
Define $\sheaf{X}_{Q,n} = \{ x \in n^{-1} \sheaf{X} : \langle x, y \rangle \in \ZZ \text{ for all } y \in \sheaf{Y}_{Q,n} \} \subset n^{-1} \sheaf{X}$.  For each root $\phi \in \Phi$, define constants $n_\phi$ and $m_\phi$,
\begin{equation}
\label{nm}
n_\phi = \frac{n}{\GCD(n, Q(\phi^\vee))} , \quad m_\phi = \frac{Q(\phi^\vee)}{\GCD(n, Q(\phi^\vee))}.
\end{equation}
Define \defined{modified roots} and \defined{modified coroots} by
$$\tilde \phi = n_\phi^{-1} \phi, \quad \tilde \phi^{\vee} = n_\phi \phi^\vee.$$
These define subsets $\tilde \Phi = \{ \tilde \phi : \phi \in \Phi \} \subset \sheaf{X}_{Q,n}$ and $\tilde \Phi^\vee = \{ \tilde \phi^\vee : \phi^\vee \in \Phi^\vee \} \subset \sheaf{Y}_{Q,n}$, as in \cite{MWCrelle}.  Modifying the simple roots and their coroots, we have subsets $\tilde \Delta \subset \tilde \Phi$ and $\tilde \Delta^\vee \subset \tilde \Phi$.

By \cite[Construction 1.3]{MWCrelle}, this defines a local system of based root data $\tilde \Psi$ on $S_{\et}$.  Write $\tilde \Psi^\vee$ for its dual,
$$\tilde \Psi^\vee = (\sheaf{Y}_{Q,n}, \tilde \Phi^\vee, \tilde \Delta^\vee, \sheaf{X}_{Q,n}, \tilde \Phi, \tilde \Delta).$$
Write $\sheaf{Y}_{Q,n}^{\SC}$ for the subgroup of $\sheaf{Y}_{Q,n}$ spanned by the modified coroots $\tilde \Phi^\vee$.  

Define $\dgp{\tilde G^\vee}$ to be the (unique up to unique isomorphism) local system on $S_{\et}$ of pinned reductive groups over $\ZZ$, associated to the local system of based root data $\tilde \Psi^\vee$.  Its maximal torus is a local system on $S_{\et}$ of split tori over $\ZZ$,
$$\dgp{\tilde T^\vee} = \Spec \left( \ZZ[ \sheaf{Y}_{Q,n}] \right).$$
The center of $\dgp{\tilde G^\vee}$ is a local system on $S_{\et}$ of groups of multiplicative type over $\ZZ$,
$$\dgp{\tilde Z^\vee} = \Spec \left( \ZZ \left[ \sheaf{Y}_{Q,n} /  \sheaf{Y}_{Q,n}^{\SC} \right] \right).$$
We call $\dgp{\tilde G^\vee}$ (endowed with its pinning) the \defined{dual group} of the cover $\alg{\tilde G}$.

\begin{proposition}
\label{DualGroupModn}
Suppose that $\alg{\tilde G}_0$ is another cover of $\alg{G}$ of degree $n$, with first Brylinski-Deligne invariant $Q_0$.  If $Q \equiv Q_0$ modulo $n$, i.e. $Q(y) - Q_0(y) \in n \ZZ$ for all $y \in \sheaf{Y}$, then the resulting modified root data are equal:  $\tilde \Psi^\vee = \tilde \Psi_0^\vee$.  Thus the dual groups are equal, $\dgp{\tilde G}^\Vee = \dgp{\tilde G}_0^\Vee$.
\end{proposition} 
\proof
One checks directly that $\beta_Q \equiv \beta_{Q_0}$ modulo $\ZZ$, from which it follows that
$$\sheaf{Y}_{Q,n} = \sheaf{Y}_{Q_0,n}, \quad \sheaf{X}_{Q,n} = \sheaf{X}_{Q_0, n}.$$
Similarly, one checks that the constants $n_\phi$ are equal,
$$\frac{n}{\GCD(n, Q_0(\phi^\vee))} = \frac{n}{\GCD(n, Q(\phi^\vee))}.$$
The result follows.
\qed

The Weyl group of $\dgp{\tilde G^\vee}$ with respect to $\dgp{\tilde T^\vee}$ forms a local system $\sheaf{\tilde W}$ on $S_{\et}$ of finite groups, generated (locally on $S_{\et}$) by reflections $s_{\tilde \phi}$ for every $\tilde \phi \in \tilde \Phi$.  The action of $\sheaf{\tilde W}$ on $\sheaf{Y}_{Q,n}$ is given by the formula
$$s_{\tilde \phi}(y) = y - \langle \tilde \phi, y \rangle \tilde \phi^\vee = y - \langle \phi, y \rangle \phi^\vee.$$
This identifies the root reflections $s_{\tilde \phi}$ with the root reflections $s_{\phi}$, and hence identifies the Weyl group $\sheaf{\tilde W}$ with the Weyl group $\sheaf{W}$ of $\alg{G}$ with respect to $\alg{T}$ (where both are viewed as local systems on $S_{\et}$ of finite groups).

The dual group $\dgp{\tilde G^\vee}$ comes with a distinguished 2-torsion element in its center, described here.  From the quadratic form $Q \From \sheaf{Y} \To \ZZ$, observe that $2 Q(y) = n \beta_Q(y,y) \in n \ZZ$ for all $y \in \sheaf{Y}_{Q,n}$.  Moreover, we have 
$$Q(\tilde \phi^\vee) = n_\phi^2 Q(\phi^\vee) = n_\phi m_\phi n \in n \ZZ,$$
for all $\phi \in \Phi$.  Of course, $Q( n y) \in n \ZZ$ as well, for all $y \in \sheaf{Y}_{Q,n}$.  We find a homomorphism of local systems of abelian groups on $S_{\et}$,
$$\overline{n^{-1} Q} \From \frac{\sheaf{Y}_{Q,n}}{\sheaf{Y}_{Q,n}^{\SC} + n \sheaf{Y}_{Q,n}} \To \half \ZZ / \ZZ, \quad y \mapsto n^{-1} Q(y) \text{ mod } \ZZ.$$

Applying $\Spec$ yields a homomorphism of local systems on $S_{\et}$ of diagonalizable group schemes over $\ZZ$,
$$\tau_Q \in \Hom ( \amu_2, \dgp{\tilde Z}_{[n]}^\Vee ).$$
Thus $\tau_Q(-1)$ is a distinguished Galois-invariant 2-torsion element in the center of $\dgp{\tilde G^\vee}$.  If $n$ is odd, then $\tau_Q(-1) = 1$.

\subsection{Well-aligned functoriality}
\label{WAFDualGroup}

Consider a well-aligned homomorphism $\tilde \iota \From \alg{\tilde G}_1 \To \alg{\tilde G}_2$ of covers, each endowed with Borel subgroup and maximally split maximal torus.  Here we construct a corresponding homomorphism of dual groups,
$$\iota^\vee \From \dgp{\tilde G}_2^\Vee \To \dgp{\tilde G}_1^\Vee.$$
These dual groups are constructed, locally on $S_{\et}$, from root data:
\begin{align*}
\tilde \Psi_1^\vee &= \left( \sheaf{Y}_{1,Q_1,n}, \tilde \Phi_1^\vee, \tilde \Delta_1^\vee, \sheaf{X}_{1,Q_1,n}, \tilde \Phi_1, \tilde \Delta_1 \right); \\
\tilde \Psi_2^\vee &= \left( \sheaf{Y}_{2,Q_2,n}, \tilde \Phi_2^\vee, \tilde \Delta_2^\vee, \sheaf{X}_{2,Q_2,n}, \tilde \Phi_2, \tilde \Delta_2 \right).
\end{align*}

For the construction of $\iota^\vee$, it suffices to work locally on $S_{\et}$, on a finite \'etale cover over which $\alg{G}_1$ and $\alg{G}_2$ split.  The well-aligned Condition (4) gives a homomorphism $\iota \From \sheaf{Y}_{1,Q_1,n} \To \sheaf{Y}_{2,Q_2,n}$, and its dual homomorphism $\iota^\ast \From \sheaf{X}_{2,Q_2,n} \To \sheaf{X}_{1,Q_1,n}$.  As $\iota$ has normal image, the coroots from $\Phi_1^\vee$ map to coroots from $\Phi_2^\vee$.  Condition (3) implies that the scaled coroots in $\tilde \Phi_1^\vee \subset \sheaf{Y}_{1,Q_1,n}$ map to scaled coroots in $\tilde \Phi_2^\vee \subset \sheaf{Y}_{2,Q_2,n}$.  Since the Borel subgroups are aligned, the simple scaled coroots in $\tilde \Delta_1^\vee$ map to simple scaled coroots in $\tilde \Delta_2^\vee$.  Dually, the map $\iota^\ast \From \sheaf{X}_{2,Q_2,n} \To \sheaf{X}_{1,Q_1,n}$ sends $\tilde \Phi_2$ to $\tilde \Phi_1$ and $\tilde \Delta_2$ to $\tilde \Delta_1$.  

This allows us to assemble a homomorphism $\iota^\vee \From \dgp{\tilde G}_2^\Vee \To \dgp{\tilde G}_1^\Vee,$ (cf.~\cite[\S 2.1, 2.5]{BorelCorvallis}).  On tori, let $\iota^\vee \From \dgp{\tilde T}_2^\Vee \To \dgp{\tilde T}_1^\Vee$ be the homomorphism dual to the map of character lattices $\iota \From \sheaf{Y}_{1,Q_1,n} \To \sheaf{Y}_{2,Q_2,n}$.  Using the pinnings on $\dgp{\tilde G}_2^\Vee$ and $\dgp{\tilde G}_1^\Vee$, and the map of roots to roots, we obtain a homomorphism from the simply-connected cover $\dgp{\tilde G}_{2,\SC}^\Vee$ of the derived subgroup, $\dgp{\tilde G}_{2, \der}^\Vee$,
$$\iota_{\SC}^\vee \From \dgp{\tilde G}_{2,\SC}^\Vee \To \dgp{\tilde G}_1^\Vee.$$
Let $\dgp{\tilde T}_{2,\SC}^\Vee$ be the pullback of $\dgp{\tilde T}_2^\Vee$.  The following diagram commutes.
$$\begin{tikzcd}
\dgp{\tilde T}_{2, \SC}^\Vee \arrow{r} \inarrow{d} & \dgp{\tilde T}_2^\Vee \arrow{d}{\iota^\vee} \\
\dgp{\tilde G}_{2, \SC}^\Vee \arrow{r}{\iota_{\SC}^\Vee} & \dgp{\tilde G}_1^\Vee
\end{tikzcd}$$

The homomorphism $\iota_{\SC}^\Vee$ descends to the derived subgroup $\dgp{\tilde G}_{2,\der}^\Vee$, since it is trivial on the kernel of $\dgp{\tilde T}_{2,\SC}^\Vee \To \dgp{\tilde T}_2^\Vee$.  In this way, we have a pair of homomorphisms of groups over $\ZZ$,
$$\iota_{\der}^\Vee \From \dgp{\tilde G}_{2,\der}^\Vee \To \dgp{\tilde G}_1^\Vee, \quad \iota^\vee \From \dgp{\tilde T}_{2}^\Vee \To \dgp{\tilde T}_1^\Vee \subset \dgp{\tilde G}_1^\Vee.$$
These homomorphisms agree on their intersection, giving the desired homomorphism $\iota^\vee \From \dgp{\tilde G}_2^\Vee \To \dgp{\tilde G}_1^\Vee$.

Since modified coroots are sent to modified coroots, we find that $\iota^\vee$ sends $\dgp{\tilde Z}_2^\Vee$ to $\dgp{\tilde Z}_1^\Vee$.  The quadratic forms $Q_1$ and $Q_2$ induce two group homomorphisms:
$$\overline{n^{-1} Q_1} \From \frac{\sheaf{Y}_{1,Q_1,n}}{\sheaf{Y}_{1,Q_1,n}^{\SC} + n \sheaf{Y}_{1,Q_1,n}} \To \half \ZZ / \ZZ, \quad \overline{n^{-1} Q_2} \From \frac{\sheaf{Y}_{2,Q_2,n}}{\sheaf{Y}_{2,Q_2,n}^{\SC} + n \sheaf{Y}_{2,Q_2,n}} \To \half \ZZ / \ZZ.$$
These define two homomorphisms of group schemes over $\ZZ$,
$$\tau_{Q_1} \From \amu_2 \To \dgp{\tilde Z}_{1, [n]}^\Vee, \quad \tau_{Q_2} \From \amu_2 \To \dgp{\tilde Z}_{2, [n]}^\Vee.$$
As $Q_2(\iota(y)) = Q_1(y)$ for all $y \in \sheaf{Y}_1$, $\iota^\vee \circ \tau_{Q_1} = \tau_{Q_2}$.  In other words, the homomorphism $\iota^\vee \From \dgp{\tilde G}_2^\Vee \To \dgp{\tilde G}_1^\Vee$ sends the center to the center, and respects the 2-torsion elements therein, $\iota^\vee(\tau_{Q_2}(-1)) = \tau_{Q_1}(-1)$.

Given a pair of well-aligned homomorphisms,
$$\alg{\tilde G}_1 \xrightarrow{(\iota_1, \iota_1')} \alg{\tilde G}_2 \xrightarrow{(\iota_2, \iota_2')} \alg{\tilde G}_3,$$
their composition is a well-aligned homomrphism $(\iota_3, \iota_3') = (\iota_2 \iota_1, \iota_2' \iota_1')$ from $\alg{\tilde G}_1$ to $\alg{\tilde G}_3$ by Proposition \ref{ComposeWellaligned}.  

This gives a commutative diagram of sheaves of abelian groups.
$$\begin{tikzcd}
\sheaf{Y}_{1,Q_1,n} \arrow{r}[swap]{\iota_1} \arrow[bend left=20]{rr}{\iota_3} & \sheaf{Y}_{2,Q_2,n} \arrow{r}[swap]{\iota_2} & \sheaf{Y}_{3,Q_3,n}.
\end{tikzcd}$$
We find such a commutative diagram for dual groups, in the opposite direction.
$$\begin{tikzcd}
\dgp{\tilde G}_3^\Vee \arrow{r}[swap]{\iota_2^\vee} \arrow[bend left=20]{rr}{\iota_3^\vee} & \dgp{\tilde G}_2^\Vee \arrow{r}[swap]{\iota_1^\vee} & \dgp{\tilde G}_1^\Vee.
\end{tikzcd}$$

Let $\Cat{DGp}_S^\ast$ denote the category whose objects are local systems $\dgp{G}^\vee$ on $S_{\et}$ of group schemes over $\ZZ$, endowed with central morphisms $\amu_2 \To \dgp{G}^\vee$ (where $\amu_2$ is viewed as the constant local system of group schemes).  Morphisms in $\Cat{DGp}_S^\ast$ are morphisms of local systems of group schemes over $\ZZ$, compatible with the central morphisms from $\amu_2$.  

Let $\Cat{WAC}_S$ (Well-Aligned-Covers) denote the category whose objects are triples $(\alg{\tilde G}, \alg{B}, \alg{T})$ where $\alg{\tilde G}$ is a cover of a reductive group $\alg{G}$ over $S$, $\alg{B}$ is a Borel subgroup of $\alg{G}$, and $\alg{T}$ is a maximally split maximal torus of $\alg{G}$ contained in $\alg{B}$.  Morphisms in $\Cat{WAC}_S$ are well-aligned homomorphisms of covers.

We have proven the following result.
\begin{proposition}
\label{FuncDualGroup}
The construction of the dual group defines a contravariant functor from $\Cat{WAC}_S$ to $\Cat{DGp}_S^\ast$.
\end{proposition}

\begin{thm}
\label{WellDefinedDualGroup}
Let $\alg{\tilde G}$ be a degree $n$ cover of a quasisplit group $\alg{G}$ over $S$.  The dual group $\dgp{\tilde G^\vee}$ is well-defined, up to unique isomorphism, by $\alg{\tilde G}$ alone.
\end{thm}
\proof
The construction of the dual group depends on the choice of Borel and torus $\alg{B} \supset \alg{T}$.  So it suffices to construct a canonical isomorphism of dual groups $\dgp{\tilde G}^\Vee \To \dgp{\tilde G}_0^\vee$ for any pair of choices $\alg{B}_0 \supset \alg{T}_0$ and $\alg{B} \supset \alg{T}$.  Such ``well-definedness'' is discussed in more detail in \cite[\S 1.1]{DeligneLusztig}.

By Proposition \ref{BTConj}, there exists $g \in \alg{G}(S)$ such that $\Int(g) \alg{T}_0 = \alg{T}$ and $\Int(g) \alg{B}_0 = \alg{B}$.  This automorphism $\Int(g)$ lifts to an automorphism of $\alg{G}'$ and defines a well-aligned isomorphism of covers, 
$$\begin{tikzcd}
\alg{K}_{2} \inarrow{r} \arrow{d}{=} & \alg{G}' \onarrow{r} \arrow{d}{\Int(g)'} & \alg{G} \arrow{d}{\Int(g)} & & \alg{B}_{0} \arrow{d} & \alg{T}_{0} \arrow{d}  \\
\alg{K}_{2} \inarrow{r} & \alg{G}' \onarrow{r} & \alg{G} & & \alg{B} & \alg{T}
\end{tikzcd}$$
As a well-aligned isomorphism of covers, this yields an isomorphism (of local systems on $S_{\et}$ of reductive groups over $\ZZ$)
$$\Int(g)^\vee \From \dgp{\tilde G^\vee} \xrightarrow{\sim} \dgp{\tilde G}_0^\Vee.$$

If $g' \in \alg{G}(S)$ also satisfies $\Int(g') \alg{T}_0 = \alg{T}$ and $\Int(g') \alg{B}_0 = \alg{B}$, then $g' g^{-1} \in \alg{N}(S) \cap \alg{B}(S) = \alg{T}(S)$.  Thus $g' = t g$ for some $t \in \alg{T}(S)$.  Hence $\Int(g') = \Int(t) \Int(g)$, and so by Proposition \ref{FuncDualGroup},
$$\Int(g')^\vee = \Int(g)^\vee \Int(t)^\vee \From \dgp{\tilde G^\vee} \To \dgp{\tilde G}_0^\Vee.$$
But $\Int(t)^\vee = \Id$, since $\Int(t)$ leaves all relevant data unchanged.  Thus $\Int(g')^\vee = \Int(g)^\vee$.  Hence we find a canonical isomorphism $\dgp{\tilde G}^\Vee \xrightarrow{\sim} \dgp{\tilde G}_0^\Vee$.
\qed

\subsection{Change of base scheme}
\label{BaseChangeDualGroup}

Let $\alg{\tilde G}$ be a degree $n$ cover of a quasisplit group $\alg{G} \supset \alg{B} \supset \alg{T}$ over $S$, as before.  Let $\gamma \From S_0 \To S$ be a morphism of schemes, with $S_0 = \Spec(F_0)$ for some field $F_0$ or $S_0 = \Spec(\OO_0)$ for some DVR $\OO_0$ (with finite residue field or containing a field, as usual).  Then $\gamma$ gives rise to a pullback functor $\gamma^\ast$ from sheaves on $S_{\et}$ to sheaves on $S_{0,\et}$.  Pullback via $\gamma$ defines a degree $n$ cover $\alg{\tilde G}_0$ of a quasisplit group $\alg{G}_0 \supset \alg{B}_0 \supset \alg{T}_0$ over $S_0$. 

The cocharacter lattice $\sheaf{Y}_0$ of $\alg{T}_0$ is a sheaf on $S_{0,\et}$.  It is naturally isomorphic to the pullback $\gamma^\ast \sheaf{Y}$, with $\sheaf{Y}$ the cocharacter lattice of $\alg{T}$.  Write $N \From \gamma^\ast \sheaf{Y} \To \sheaf{Y}_0$ for the natural isomorphism.  The quadratic form $Q \From \sheaf{Y} \To \ZZ$ pulls back to a quadratic form $\gamma^\ast Q \From \gamma^\ast \sheaf{Y} \To \ZZ$.  The compatibility of Brylinski-Deligne invariants with pullbacks implies that $\gamma^\ast Q = N^\ast Q_0$, with $Q_0$ the first Brylinski-Deligne invariant of $\alg{\tilde G}_0$.  

\begin{remark}
This compatibility follows straightforwardly from \cite[\S 3.10]{B-D}; if $\alg{\tilde T}$ arises (after a finite \'etale $U \To S$) from the cocycle attached to $C \in \sheaf{X} \otimes \sheaf{X}$, then $\alg{\tilde T}_0$ arises from the pullback of this cocycle, i.e., from an element $C_0 \in \sheaf{X}_0 \otimes \sheaf{X}_0$ with $N^\ast C_0 = \gamma^\ast C$.  The quadratic form $Q$ is given by $Q(y) = C(y,y)$ and similarly $Q_0(y) = C_0(y,y)$.  Since $N^\ast C_0 = \gamma^\ast C$, we find that $N^\ast Q_0 = \gamma^\ast Q$.
\end{remark}

In this way, we find that $N$ restricts to an isomorphism from $\gamma^\ast \sheaf{Y}_{Q,n}$ to $\sheaf{Y}_{0,Q_0,n}$, sending roots and coroots (the sheaves of sets $\gamma^\ast \tilde \Phi$, $\gamma^\ast \tilde \Phi^\vee$ on $S_{0,\et}$) to the corresponding roots and coroots in $\sheaf{Y}_{0,Q_0,n}$.  As $\alg{B}_0 = \gamma^\ast \alg{B}$, simple roots and coroots are identified as well.  

By our construction of the dual group, we find that $N$ gives an isomorphism of local systems on $S_{0,\et}$ of pinned reductive groups over $\ZZ$,
$$N^\vee \From \gamma^\ast \dgp{\tilde G^\vee} \xrightarrow{\sim} \dgp{\tilde G}_0^\Vee.$$

\subsection{Parabolic subgroups}
\label{ParabolicDualGroup}

We keep the degree $n$ cover $\alg{\tilde G}$ of the quasisplit group $\alg{G} \supset \alg{B} \supset \alg{T}$ over $S$.  Let $\alg{P} \subset \alg{G}$ be a parabolic subgroup defined over $S$, containing $\alg{B}$.  Suppose that $\alg{P} = \alg{M} \alg{N}$ is a Levi decomposition defined over $S$, with $\alg{N}$ the unipotent radical of $\alg{P}$ and $\alg{M}$ a Levi factor containing $\alg{T}$.  Let $\alg{B}_{\alg{M}}$ denote the Borel subgroup $\alg{B} \cap \alg{M}$ of $\alg{M}$.  Write $\alg{\tilde M} = (\alg{M}', n)$ for the cover of $\alg{M}$ arising from pulling back $\alg{\tilde G}$.

The first Brylinski-Deligne invariant is the same for $\alg{\tilde G}$ as for $\alg{\tilde M}$, as it depends only upon the cover $\alg{\tilde T}$ of their common maximal torus.  Write $\Phi_{\alg{M}}$ and $\Phi_{\alg{M}}^\vee$ for the roots and coroots of $\alg{M}$; these are subsets of $\Phi$ and $\Phi^\vee$, respectively, and by agreement of the first Brylinski-Deligne invariant,
$$\tilde \Phi_{\alg{M}} \subset \tilde \Phi, \quad \tilde \Phi_{\alg{M}}^\vee \subset \tilde \Phi^\vee.$$
As $\alg{B}_{\alg{M}} = \alg{B} \cap \alg{M}$, we have $\tilde \Delta_{\alg{M}} \subset \tilde \Delta$ and $\tilde \Delta_{\alg{M}}^\vee \subset \tilde \Delta^\vee$.  We find a pair of local systems on $S_{\et}$ of based root data,
$$\left( \sheaf{Y}_{Q,n}, \tilde \Phi^\vee, \tilde \Delta^\vee, \sheaf{X}_{Q,n}, \tilde \Phi, \tilde \Delta \right), \quad \left( \sheaf{Y}_{Q,n}, \tilde \Phi_{\alg{M}}^\vee, \tilde \Delta_{\alg{M}}^\vee, \sheaf{X}_{Q,n}, \tilde \Phi_{\alg{M}}, \tilde \Delta_{\alg{M}} \right).$$
The first root datum defines a local system $\dgp{\tilde G^\vee}$ on $S_{\et}$ of pinned reductive groups over $\ZZ$.  The second root datum defines a local system $\dgp{\tilde M^\vee}$ on $S_{\et}$ of pinned reductive Levi subgroups of $\dgp{\tilde G^\vee}$.

Thus the dual group of $\alg{\tilde M}$ is naturally a Levi subgroup of the dual group $\dgp{\tilde G^\vee}$.  Moreover, by agreement of the first Brylinski-Deligne invariants, the central 2-torsion element of $\dgp{\tilde M^\vee}$ coincides with the central 2-torsion element $\tau_Q(-1) \in \dgp{\tilde G^\vee}$.

\subsection{Specific cases}

Here we give many specific cases of covers, associated dual groups, and 2-torsion elements in their centers corresponding to $\tau_Q(-1)$.  

\subsubsection{Method for simply-connected groups}

Let $\alg{G}$ be a simply-connected semisimple group over $S$, with Borel subgroup $\alg{B}$ containing a maximally split maximal torus $\alg{T}$ over $S$.  Let $S' / S$ be a Galois cover over which $\alg{T}$ splits.  If $t$ is an integer, there exists a unique Weyl-invariant quadratic form $Q_t$ with value $t$ on all short coroots.  This $Q_t$ is a multiple of the Killing form.

Corresponding to $Q_t$, there is a unique, up to unique isomorphism, object $\alg{G}^{(t)} \in \Cat{CExt}_S(\alg{G}, \alg{K}_2)$ with first Brylinski-Deligne invariant $Q_t$ (by \cite[\S 7.3(i)]{B-D} when working over a field or \cite[\S 3.3.3]{MWIntegral} over a DVR).  Write $\beta_t$ for the resulting $n^{-1} \ZZ$-valued bilinear form.  Weyl-invariance of $Q_t$ implies that
\begin{equation}
\beta_t(\phi^\vee, y) = n^{-1} Q_t(\phi^\vee) \langle \phi, y \rangle, \text{ for all } \phi \in \Phi, \text{ and all } y \in \sheaf{Y}.
\end{equation}
(See the proof of Lemma \ref{DeltaM} for a derivation.)  It follows that, for any $y \in \sheaf{Y}$,
\begin{equation}
\label{BetaObs}
\beta_t(\phi^\vee, y) \in \ZZ \text{ if and only if } \langle \phi,y \rangle \in n_\phi \ZZ.
\end{equation}

Let $Y = \sheaf{Y}[S']$, $X = \sheaf{X}[S']$, and $\Delta = \Delta[S'] = \{ \alpha_1, \ldots, \alpha_\ell \}$ a basis of simple roots corresponding to $\alg{B}$.  Since we assume $\alg{G}$ is simply connected, $Y = \bigoplus_{i = 1}^\ell \alpha_i^\vee \ZZ$.  Write $n_i = n_{\alpha_i}$.  From \eqref{BetaObs}, we find a characterization of $Y_{Q,n} = \sheaf{Y}_{Q,n}[S']$:
\begin{equation}
\label{YQnCharacterization}
y \in Y_{Q,n} \text{ if and only if } \langle \alpha_i, y \rangle \in n_i \ZZ \text{ for all } 1 \leq i \leq \ell.
\end{equation}
This, in turn, can be used to tabulate dual groups.  We provide tables here reference, noting that such information can also be found in the examples of \cite[\S 2.4]{FinkelbergLysenko}.    But we also include data on the central 2-torsion elements that we have not found in the literature.  Our tabulation was greatly assisted by using SAGE \cite{sage}, especially the recently updated package which deftly handles root data.

In what follows, we write $\alg{\tilde G}^\Vee = \dgp{\tilde G}^\Vee[S']$ for the dual group over $S'$ (a pinned reductive group scheme over $\ZZ$); it is convenient to view the dual group as a group scheme over $\ZZ$ endowed with a $\Gal(S'/S)$-action.

\subsubsection{$\alg{SL}_{\ell + 1}$}
For $\alg{G} = \alg{SL}_{\ell + 1}$, the standard Borel subgroup and maximal torus, and system of roots $\alpha_1, \ldots, \alpha_\ell$, the Dynkin diagram is
\begin{center}
\begin{tikzpicture}
\draw (1,0) -- (3.7,0);
\draw (5.3,0) -- (6,0);
\draw[dotted] (3.7,0) -- (5.3,0);
\foreach \i/\l in {1/1,2/2,3/3,6/\ell}
{ 
\filldraw[fill=gray, draw=black] (\i,0) circle (0.1) node[below=2pt] {$\alpha_{\l}$};
}
\end{tikzpicture}
\end{center}
Consider the cover $\alg{\tilde G}$ of degree $n$ arising from the quadratic form satisfying $Q(\phi^\vee) = 1$ for all $\phi \in \Phi$.  Then $n_i  = n$ and $\tilde \alpha_i^\vee = n \alpha_i^\vee$ for all $1 \leq i \leq \ell$.  Hence
$$Y_{Q,n}^{\SC} = n Y = n \alpha_1^\vee \ZZ + \cdots + n \alpha_\ell^\vee \ZZ.$$

The Cartan matrix of $\sch{\tilde G}^\Vee$ has entries $\tilde C_{ij} = \langle \tilde \alpha_i, \tilde \alpha_j^\vee \rangle = \langle \alpha_i, \alpha_j^\vee \rangle$, and so $\sch{\tilde G}^\Vee$ is isogenous to $\sch{SL}_{\ell + 1}$.  To determine the dual group up to isomorphism, it suffices to compute the order of the center, since the center of $\sch{SL}_{\ell + 1}$ is $\amu_{\ell+1}$.  The order of the center is equal to the index $\# Y_{Q,n} / Y_{Q,n}^{\SC}$, and this is computable from \eqref{YQnCharacterization}.  The results are given in Table \ref{DualATable}.

\begin{table}[!htbp]
\begin{tabular}{l|ccccc}
 & \multicolumn{5}{c}{Group $\alg{G}$} \\
$n$ \phantom{M} & $\alg{SL}_2$ & $\alg{SL}_3$ & $\alg{SL}_4$ & $\alg{SL}_5$ & $\alg{SL}_6$ \\ \hline \hline
1 & $\sch{PGL}_2$ & $\sch{PGL}_3$ & $\sch{PGL}_4$ & $\sch{PGL}_5$ & $\sch{PGL}_6$ \\
2 & ${}^{\ast} \sch{SL}_2$  & $\sch{PGL}_3$ & $\sch{SL}_4 / \amu_2$ & $\sch{PGL}_5$ & ${}^{\ast}\sch{SL}_6 / \amu_3$ \\  
3 & $\sch{PGL}_2$ & $\sch{SL}_3$ & $\sch{PGL}_4$ & $\sch{PGL}_5$ & $\sch{SL}_6 / \amu_2$ \\
4 & $\sch{SL}_2$ & $\sch{PGL}_3$ & ${}^\ast \sch{SL}_4$ & $\sch{PGL}_5$ & $\sch{SL}_6 / \amu_3$ \\
5 & $\sch{PGL}_2$ & $\sch{PGL}_3$ & $\sch{PGL}_4$ & $\sch{SL}_5$ & $\sch{PGL}_6$ \\
6 & ${}^\ast \sch{SL}_2$ & $\sch{SL}_3$ & $\sch{SL}_4 / \amu_2$ & $\sch{PGL}_5$ & ${}^\ast \sch{SL}_6$ \\

\end{tabular}
\caption{Table of dual groups for degree $n$ covers of $\alg{SL}_{\ell+1}$.  Groups marked with $\ast$ have $\tau_Q(-1)$ nontrivial.}
\label{DualATable}
\end{table}

The dual groups $\sch{\tilde G}^\Vee$ are consistent with the Iwahori-Hecke algebra isomorphisms found by Savin in \cite[Theorem 7.8]{SavinUnramified}.  In other words, the dual group $\sch{\tilde G}^\Vee$ coincides with the Langlands dual group of $\alg{SL}_{\ell + 1} / \alg{Z}_{[n]}$, where $\alg{Z}_{[n]}$ is the $n$-torsion subgroup of the center of $\alg{SL}_{\ell+1}$.  

The central 2-torsion elements $\tau_Q(-1)$ follow a somewhat predictable pattern.  For covering degree $2$, $\tau_Q(-1)$ is nontrivial for $\alg{SL}_2$, $\alg{SL}_6$, $\alg{SL}_{10}$, $\alg{SL}_{14}$, etc..  In covering degree $4$, $\tau_Q(-1)$ is nontrivial for $\alg{SL}_4$, $\alg{SL}_{12}$, $\alg{SL}_{20}$, etc..  In covering degree $6$, $\tau_Q(-1)$ is nontrivial for $\alg{SL}_2$, $\alg{SL}_6$, $\alg{SL}_{10}$, $\alg{SL}_{14}$, etc..  In covering degree $8$, $\tau_Q(-1)$ is nontrivial for $\alg{SL}_8$, $\alg{SL}_{24}$, etc..  In general, we suspect the following:
\begin{center}
$\tau_Q(-1)$ is nontrivial for a degree $2^e \cdot k$ ($k$ odd) cover of $\alg{SL}_m$ if and only if $m = 2^e \cdot j$ for $j$ odd.
\end{center} 

\subsubsection{$\alg{Spin}_{2 \ell + 1}$}
For $\alg{G} = \alg{Spin}_{2 \ell + 1}$, the Dynkin diagram has type B.
\begin{center}
\begin{tikzpicture}
\draw (1,0) -- (3.7,0);
\draw (5.3,0) -- (6,0);
\draw[dotted] (3.7,0) -- (5.3,0);
\draw (6,0.05) -- (7,0.05);
\draw (6,-0.05) -- (7,-0.05);
\draw (6.45,0.125) -- (6.55,0) -- (6.45,-0.125);
\foreach \i/\l in {1/1,2/2,3/3,6/{\ell-1}, 7/\ell}
{ 
\filldraw[fill=gray, draw=black] (\i,0) circle (0.1) node[below=2pt] {$\alpha_{\l}$};
}
\end{tikzpicture}
\end{center}
Let $\alg{\tilde G}$ be the cover of degree $n$, associated to the quadratic form taking the value $1$ at all short coroots.  Thus $Q(\alpha_i^\vee) = 2$ for all long coroots $1 \leq i \leq \ell - 1$.  If $n$ is odd, then $n_\alpha = n$ for all coroots $\alpha^\vee$.  If $n$ is even, then $n_i = n/2$ for $1 \leq i \leq \ell-1$ and $n_\ell = n$.  When $n$ is even, short coroots become long and long become short, after modification.  We find that the dual group is isogenous to $\sch{Sp}_{2 \ell}$ if $n$ is odd, and is isogenous to $\sch{Spin}_{2 \ell + 1}$ if $n$ is even.

The centers of $\sch{Sp}_{2 \ell}$ and $\sch{Spin}_{2 \ell + 1}$ are cyclic of order two.  Thus the dual group can be identified by the order of its center.

\begin{table}[!htbp]
\begin{tabular}{l|cccccc}
 & \multicolumn{6}{c}{Group $\alg{G}$} \\
$n$ \phantom{M} & $\sch{Spin}_7$ & $\sch{Spin}_9$ & $\sch{Spin}_{11}$ & $\sch{Spin}_{13}$ & $\sch{Spin}_{15}$ & $\sch{Spin}_{17}$  \\ \hline \hline
1 & $\sch{PGSp}_{6}$ & $\sch{PGSp}_{8}$ & $\sch{PGSp}_{10}$ & $\sch{PGSp}_{12}$ & $\sch{PGSp}_{14}$ & $\sch{PGSp}_{16}$ \\
2 &  $\sch{SO}_{7}$ & $\sch{Spin}_{9}$ & $\sch{SO}_{11}$ & ${}^\ast \sch{Spin}_{13}$ &$ \sch{SO}_{15}$ & $\sch{Spin}_{17}$ \\
3 & $\sch{PGSp}_{6}$ & $\sch{PGSp}_{8}$ & $\sch{PGSp}_{10}$ & $\sch{PGSp}_{12}$ & $\sch{PGSp}_{14}$ & $\sch{PGSp}_{16}$ \\
4 & ${}^\ast  \sch{Spin}_{7}$ & $\sch{Spin}_{9}$ & ${}^\ast \sch{Spin}_{11}$ & $\sch{Spin}_{13}$ & ${}^\ast \sch{Spin}_{15}$ & $\sch{Spin}_{17}$ \\
5 & $\sch{PGSp}_{6} $ & $\sch{PGSp}_{8} $ & $\sch{PGSp}_{10} $ & $\sch{PGSp}_{12}$ &$ \sch{PGSp}_{14}$ & $\sch{PGSp}_{16}$ \\
6 &  $\sch{SO}_{7}$ & $\sch{Spin}_{9}$ & $\sch{SO}_{11}$ & ${}^\ast \sch{Spin}_{13}$ & $\sch{SO}_{15}$ & $\sch{Spin}_{17}$  \\
\end{tabular}
\caption{Table of dual groups for degree $n$ covers of $\alg{Spin}_{2 \ell+1}$.  Groups marked with $\ast$ have $\tau_Q(-1)$ nontrivial.}
\label{DualBTable}
\end{table}
Table \ref{DualBTable} describes the dual groups.  Note that, in this case, the isogeny class of the dual group depends on the covering degree modulo $4$.  In covering degree $4k + 2$, we find that $\tau_Q(-1)$ is nontrivial for $\alg{Spin}_{8 j + 5}$ (corresponding to rank $4j + 2$) for all positive integers $j$.  In covering degree $4k$, we find that $\tau_Q(-1)$ is nontrivial for $\alg{Spin}_{4j + 3}$ for all positive integers $j$.

\subsubsection{$\alg{Sp}_{2 \ell}$}
For $\alg{G} = \alg{Sp}_{2 \ell}$, the Dynkin diagram has type C.
\begin{center}
\begin{tikzpicture}
\draw (1,0) -- (3.7,0);
\draw (5.3,0) -- (6,0);
\draw[dotted] (3.7,0) -- (5.3,0);
\draw (6,0.05) -- (7,0.05);
\draw (6,-0.05) -- (7,-0.05);
\draw (6.55,0.125) -- (6.45,0) -- (6.55,-0.125);
\foreach \i/\l in {1/1,2/2,3/3,6/{\ell-1}, 7/\ell}
{ 
\filldraw[fill=gray, draw=black] (\i,0) circle (0.1) node[below=2pt] {$\alpha_{\l}$};
}
\end{tikzpicture}
\end{center}
Let $\alg{\tilde G}$ be the cover of degree $n$, associated to the quadratic form taking the value $1$ at all short coroots.  As in type B, we find that short coroots become long, and long become short, after modification when $n$ is even.  We find that the dual group is isogenous to $\sch{Sp}_{2 \ell}$ if $n$ is even, and is isogenous to $\sch{Spin}_{2 \ell + 1}$ if $n$ is odd.  As before, the dual group can be identified by the order of its center.

\begin{table}[!htbp]
\begin{tabular}{l|ccc}
 & \multicolumn{3}{c}{Group $\alg{G}$} \\
$n$ \phantom{M} & $\alg{Sp}_6$ & $\alg{Sp}_8$ & $\alg{Sp}_{10}$ \\ \hline \hline
1 & $\sch{SO}_{7}$ & $\sch{SO}_{9}$ & $\sch{SO}_{11}$ \\
2 & ${}^\ast \sch{Sp}_{6}$ & ${}^\ast \sch{Sp}_{8}$ & ${}^\ast \sch{Sp}_{10}$ \\
3 &  $\sch{SO}_{7}$ & $\sch{SO}_{9}$ & $\sch{SO}_{11}$ \\
4 & $\sch{Sp}_{6}$ & $\sch{Sp}_{8}$ & $\sch{Sp}_{10}$ \\
5 &  $\sch{SO}_{7}$ & $\sch{SO}_{9}$ & $\sch{SO}_{11}$ \\
6 & ${}^\ast \sch{Sp}_{6}$ & ${}^\ast \sch{Sp}_{8}$ & ${}^\ast \sch{Sp}_{10}$ \\
\end{tabular}
\caption{Table of dual groups for degree $n$ covers of $\alg{Sp}_{2 \ell}$.  Groups marked with $\ast$ have $\tau_Q(-1)$ nontrivial.}
\label{DualCTable}
\end{table}
As Table \ref{DualCTable} illustrates, the dual group of the degree $n$ cover of $\alg{Sp}_{2 \ell}$ is the simply-connected Chevalley group $\sch{Sp}_{2 \ell}$ when $n$ is even, and the dual group is the adjoint group $\sch{SO}_{2 \ell + 1} = \sch{Spin}_{2 \ell + 1} / \amu_2$ when $n$ is odd.  The central 2-torsion element $\tau_Q(-1)$ is nontrivial when the covering degree is $4k + 2$ for some non-negative integer $k$.  This is consistent (in covering degree $2$) with expectations from the classical theta correspondence for metaplectic groups.

\subsubsection{$\alg{Spin}_{2 \ell}$}
For $\alg{G} = \alg{Spin}_{2 \ell}$, $\ell \geq 4$, the Dynkin diagram has type D.
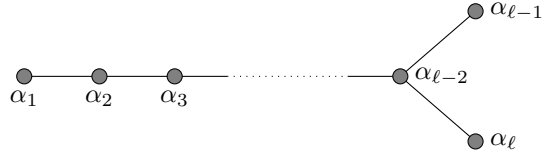
\begin{center}
\begin{tikzpicture}
\draw (1,0) -- (3.7,0);
\draw (5.3,0) -- (6,0);
\draw[dotted] (3.7,0) -- (5.3,0);
\draw (6,0) -- (7,0.866);
\draw (6,0) -- (7,-0.866);
\filldraw[fill=gray, draw=black] (6,0) circle (0.1) node[right=2pt] {$\alpha_{\ell-2}$};
\filldraw[fill=gray, draw=black] (7,0.866) circle (0.1) node[right=2pt] {$\alpha_{\ell-1}$};
\filldraw[fill=gray, draw=black] (7,-0.866) circle (0.1) node[right=2pt] {$\alpha_{\ell}$};
\foreach \i/\j/\l in {1/1,2/2,3/3}
{ 
\filldraw[fill=gray, draw=black] (\i,0) circle (0.1) node[below=2pt] {$\alpha_{\l}$};
}
\end{tikzpicture}
\end{center}
Let $\alg{\tilde G}$ be the cover of degree $n$, associated to the quadratic form taking the value $1$ at all coroots.  By the same methods as in type A, we find that the dual group is isogenous to $\sch{Spin}_{2 \ell}$.  

If $\ell$ is odd, then the center of $\sch{Spin}_{2 \ell}$ is a cyclic group of order $4$.  In this case, the dual group is determined by the order of its center.

If $\ell$ is even, then the center of $\sch{Spin}_{2 \ell}$ is isomorphic to $\amu_2 \times \amu_2$, and so the dual group is not a priori determined by the order of its center.  But fortunately, the order of the center of the dual group always equals $1$ or $4$ when $\ell$ is even, and this suffices to identify the dual group.

\begin{table}[!htbp]
\begin{tabular}{l|cccccc}
 & \multicolumn{6}{c}{Group $\alg{G}$} \\
$n$ \phantom{M} & $\alg{Spin}_8$ & $\alg{Spin}_{10}$ & $\alg{Spin}_{12}$ & $\alg{Spin}_{14}$ & $\alg{Spin}_{16}$ & $\alg{Spin}_{18}$ \\ \hline \hline
1 & $\sch{PGO}_8$ & $\sch{PGO}_{10}$ & $\sch{PGO}_{12}$ & $\sch{PGO}_{14}$ & $\sch{PGO}_{16}$ & $\sch{PGO}_{18}$  \\
2 & $\sch{Spin}_8$ & $\sch{SO}_{10}$ & ${}^\ast \sch{Spin}_{12}$ & $\sch{SO}_{14}$ & $\sch{Spin}_{16}$ & $\sch{SO}_{18}$ \\
3 & $\sch{PGO}_8$ & $\sch{PGO}_{10}$ & $\sch{PGO}_{12}$ & $\sch{PGO}_{14}$ & $\sch{PGO}_{16}$ & $\sch{PGO}_{18}$  \\
4 & $\sch{Spin}_8$ & ${}^\ast \sch{Spin}_{10}$ & $\sch{Spin}_{12}$ & ${}^\ast \sch{Spin}_{14}$ & $\sch{Spin}_{16}$ & ${}^\ast \sch{Spin}_{18}$\\
5 & $\sch{PGO}_8$ & $\sch{PGO}_{10}$ & $\sch{PGO}_{12}$ & $\sch{PGO}_{14}$ & $\sch{PGO}_{16}$ & $\sch{PGO}_{18}$ \\
6 & $\sch{Spin}_8$ & $\sch{SO}_{10}$ & ${}^\ast \sch{Spin}_{12}$ & $\sch{SO}_{14}$ & $\sch{Spin}_{16}$ & $\sch{SO}_{18}$ \\
\end{tabular}
\caption{Table of dual groups for degree $n$ covers of $\alg{Spin}_{2 \ell}$.  Groups marked with $\ast$ have $\tau_Q(-1)$ nontrivial.}
\label{DualDTable}
\end{table}
As Table \ref{DualDTable} illustrates, the dual group of an odd-degree cover of $\alg{Spin}_{2 \ell}$ coincides with the Langlands dual group of the linear group $\alg{Spin}_{2 \ell}$; this dual group is the adjoint group $\sch{PGO}_{2 \ell}$.  But the dual group of an even-degree cover of $\alg{Spin}_{2 \ell}$ depends on the parity of $\ell$ and the covering degree modulo $4$.  As in type A, these dual groups agree with expectations from the Hecke algebra isomorphisms of Savin \cite{SavinUnramified}.

When the covering degree is a multiple of $4$, the element $\tau_Q(-1)$ is nontrivial for $\alg{Spin}_{4j + 2}$ for all $j \geq 2$.  Since $\alg{Spin}_{4j + 2}$ has a unique central element of order two, this suffices to describe $\tau_Q$.  When the covering degree has the form $4k +2$, the element $\tau_Q(-1)$ is nontrivial for $\alg{Spin}_{8j + 4}$ for all $j \geq 1$.  The center of the group $\sch{Spin}_{8j + 4}$ is isomorphic to $\amu_2 \times \amu_2$, which has three distinct 2-torsion elements.  However, only one of these is invariant under the nontrivial outer automorphism of the pinned Chevalley group $\sch{Spin}_{8j + 4}$.  This one must be $\tau_Q(-1)$, since $Q$ is invariant under this outer automorphism.

\subsubsection{Exceptional groups}
Let $\alg{G}$ be a simply-connected split simple group of type $\text{E}_\ell$ (with $\ell = 6, 7, 8$), $\text{F}_4$, or $\text{G}_2$.  Let $\alg{\tilde G}$ be the cover of degree $n$, associated to the quadratic form taking the value $1$ at all short coroots (all coroots in type E).  As in types A and D, we find that the dual group is semisimple and isogenous to the Chevalley group of the same type as $\alg{G}$.  In types $\text{E}_8$, $\text{F}_4$, and $\text{G}_2$, the simply-connected group is centerless, and so $\sch{\tilde G}^\Vee$ coincides with the simply-connected Chevalley group of the same type.

The center of $\sch{E}_6$ has order $3$, and the center of $\sch{E}_7$ has order $2$.  Hence the dual group $\sch{\tilde G}^\Vee$ is determined by the order of its center.  The dual groups are listed in Table \ref{DualETable}.
\begin{table}[!htbp]
\begin{tabular}{l|ccccc}
 & \multicolumn{5}{c}{Group $\alg{G}$} \\
$n$ \phantom{M} & $\alg{E}_6$ & $\alg{E}_7$ & $\alg{E}_8$ & $\alg{F}_4$ & $\alg{G}_2$ \\ \hline
1 & $\sch{E}_6 / \amu_3$ & $\sch{E}_7 / \amu_2$ & $\sch{E}_8$ & $\sch{F}_4$ & $\sch{G}_2$ \\
2 & $\sch{E}_6 / \amu_3$ & ${}^\ast \sch{E}_7$ & $\sch{E}_8$ & $\sch{F}_4$ & $\sch{G}_2$ \\
3 & $\sch{E}_6$ & $\sch{E}_7 / \amu_2$ & $\sch{E}_8$ & $\sch{F}_4$ & $\sch{G}_2$ \\
4 & $\sch{E}_6 / \amu_3$ & $\sch{E}_7$ & $\sch{E}_8$ & $\sch{F}_4$ & $\sch{G}_2$  \\
5 & $\sch{E}_6 / \amu_3$ & $\sch{E}_7 / \amu_2$ & $\sch{E}_8$ & $\sch{F}_4$ & $\sch{G}_2$ \\
6 & $\sch{E}_6$ & ${}^\ast \sch{E}_7$ & $\sch{E}_8$ & $\sch{F}_4$ & $\sch{G}_2$ \\
\hline
\end{tabular}
\caption{Table of dual groups for degree $n$ covers of exceptional groups.  Groups marked with $\ast$ have $\tau_Q(-1)$ nontrivial.}
\label{DualETable}
\end{table}
In type E, these dual groups agree with expectations from \cite{SavinUnramified}.  The central 2-torsion element is nontrivial for $\alg{E}_7$, when the covering degree equals $4j + 2$ for some $j \geq 0$.

\subsubsection{$\alg{GL}_r$}

Suppose that $\alg{G}$ is split reductive, and the derived subgroup of $\alg{G}$ is simply-connected.  Let $\alg{T}$ be a split maximal torus in $\alg{G}$ with cocharacter lattice $Y$.  Then, for any Weyl-invariant quadratic form $Q \From Y \To \ZZ$, there exists a cover $\alg{\tilde G}$ with first Brylinski-Deligne invariant $Q$.

For example, when $\alg{G} = \alg{GL}_r$, there is a two-parameter family of such Weyl-invariant quadratic forms.  Write $\alg{T}$ for the standard maximal torus of diagonal matrices, and identify $Y = \ZZ^r$ in the usual way.  For any pair of integers $q,c$, there exists a unique Weyl-invariant quadratic form $Q_{q,c}$ satisfying
$$Q(1,-1,0,\ldots, 0) = q \text{ and } Q(1,0,\ldots, 0) = 1+c.$$
The $n$-fold covers $\widetilde{GL}_r^{(c)}$ studied by Kazhdan and Patterson \cite[\S 0.1]{KazhdanPatterson} can be constructed from Brylinski-Deligne extensions with first invariant $Q_{1,c}$.  The proof of following result is left to the reader.
\begin{proposition}
Let $\alg{\tilde G}$ be a degree $n$ cover of $\alg{G} = \alg{GL}_r$ with first Brylinski-Deligne invariant $Q_{1,c}$.  If $\GCD(n, 1+r + 2rc) = 1$, then $\sch{\tilde G}^\Vee$ is isomorphic to $\sch{GL}_r$.  If $\GCD(n,r) = 1$, then the derived subgroup of $\sch{\tilde G}^\Vee$ is isomorphic to $\sch{SL}_r$ and thus there exists an isogeny $\sch{\tilde G}^\Vee \To \sch{GL}_r$.
\end{proposition}
This may place the work of Kazhdan and Flicker \cite{KazhdanFlicker} in a functorial context.

\subsubsection{$\alg{GSp}_{2 r}$}

For $\alg{G} = \alg{GSp}_{2r}$, and a standard choice of split maximal torus and Borel subgroup, we write $e_0, \ldots, e_r$ for a basis of $Y$, $f_0, \ldots, f_r$ for the dual basis of $X$, and the simple roots and coroots are
$$\alpha_1 = f_1 - f_2, \ldots, \alpha_{r-1} = f_{r-1} - f_r, \quad \alpha_r = 2 f_r - f_0;$$
$$\alpha_1^\vee = e_1 - e_2, \ldots, \alpha_{r-1}^\vee = e_{r-1} - e_r, \quad \alpha_r^\vee = e_r.$$
The Weyl group is $S_r \ltimes \mu_2^r$, with $S_r$ acting by permutation of indices $1, \ldots, r$ (fixing $e_0$ and $f_0$), and elements $w_j$ (for $1 \leq j \leq r$) of order two which satisfy 
$$w_j(e_j) = - e_j, \quad w_j(e_i) = e_i \text{ for } i \neq j,0, \quad w_j(e_0) = e_0 + e_1.$$
Weyl-invariant quadratic forms on $Y$ are in bijection with pairs $(\kappa, \nu)$ of integers.  For any such pair, there is a unique Weyl-invariant quadratic form $Q_{\kappa, \nu}$ satisfying
$$Q_{\kappa, \nu}(e_0) = \kappa, \quad Q_{\kappa, \nu}(e_i) = \nu \text{ for } 1 \leq i \leq r.$$
The proof of following result is left to the reader.
\begin{proposition}
Let $\alg{\tilde G}$ be a degree $2$ cover of $\alg{G} = \alg{GSp}_{2r}$, with first Brylinski-Deligne invariant $Q_{0,1}$.  Then the dual group $\sch{\tilde G}^\Vee$ is isomorphic to $\sch{GSp}_{2r}$ if $r$ is odd, and to $\sch{PGSp}_{2r} \times \sch{G}_m$ if $r$ is even.
\end{proposition}
We find that that double-covers of $\alg{GSp}_{2r}$ behave differently depending on the parity of $r$; this phenomenon is consistent with the work of Szpruch \cite{Szpruch} on principal series.

\section{The gerbe associated to a cover}

In this section, we construct a gerbe $\gerb{E}_\epsilon(\alg{\tilde G})$ on $S_{\et}$ associated to a degree $n$ cover $\alg{\tilde G}$ of a quasisplit group $\alg{G}$, and a choice of injective character $\epsilon \From \mu_n \Into \CC^\times$.  Fix $\alg{\tilde G}$, $\alg{G}$, and $\epsilon$ throughout.  Also, choose a Borel subgroup containing a maximally split maximal torus $\alg{B} \supset \alg{T}$; we will see that our construction is independent of this choice (in a 2-categorical sense).

We make one assumption about our cover $\alg{\tilde G}$, which enables our construction and is essentially nonrestrictive.
\begin{assumption}[Odd $n$ implies even $Q$]
\label{OddnEvenQ}
If $n$ is odd, then we assume $Q \From \sheaf{Y} \To \ZZ$ takes only even values.
\end{assumption}

If $\alg{\tilde G}$ does not satisfy this assumption, i.e., $n$ is odd and $Q(y)$ is odd for some $y \in \sheaf{Y}$, then replace $\alg{\tilde G}$ by $(n+1) \dot \times \alg{\tilde G}$ (its Baer sum with itself $n+1$ times).  The first Brylinski-Deligne invariant becomes $(n+1) Q$, which is even-valued.  By Proposition \ref{DualGroupModn}, the dual group $\dgp{\tilde G^\vee}$ does not change since $Q \equiv (n+1) Q$ modulo $n$.  Moreover, the resulting extensions of groups over local or global fields, e.g., $\mu_n \Into \tilde G \Onto G$, remain the same (up to natural isomorphism).  Indeed, the Baer sum of $\tilde G$ with itself $n+1$ times is naturally isomorphic to the pushout via the $(n+1)^{\th}$ power map $\mu_n \To \mu_n$, which equals the identity map.   

We work with sheaves of abelian groups on $S_{\et}$, and great care is required to avoid confusion between those in the left column and the right column below.  Define
\begin{align*}
\sheaf{\hat T} = \shom(\sheaf{Y}_{Q,n}, \sheaf{G}_m), & \quad \sheaf{\tilde T^\vee} = \shom(\sheaf{Y}_{Q,n}, \CC^\times); \\
\sheaf{\hat T}_{\SC} = \shom(\sheaf{Y}_{Q,n}^{\SC}, \sheaf{G}_m), & \quad \sheaf{\tilde T}_{\SC}^\Vee = \shom(\sheaf{Y}_{Q,n}^{\SC}, \CC^\times); \\
\sheaf{\hat Z} = \shom(\sheaf{Y}_{Q,n} / \sheaf{Y}_{Q,n}^{\SC}, \sheaf{G}_m), & \quad \sheaf{\tilde Z^\vee} = \shom(\sheaf{Y}_{Q,n} / \sheaf{Y}_{Q,n}^{\SC}, \CC^\times).
\end{align*}
Here $\CC^\times$ denotes the constant sheaf on $S_{\et}$.  Thus, in the right column, we find the complex points of the dual groups,
$$\sheaf{\tilde T^\vee} = \dgp{\tilde T}^\Vee(\CC), \quad \sheaf{\tilde T}_{\SC}^\Vee = \dgp{\tilde T}_{\SC}^\Vee(\CC), \quad \sheaf{\tilde Z^\vee} = \dgp{\tilde Z}^\Vee(\CC).$$

Composing with $\epsilon$ defines homomorphisms of local systems of abelian groups,
\begin{align*}
\sheaf{\hat T}_{[n]} = \shom(\sheaf{Y}_{Q,n}, \mu_n) & \xrightarrow{\epsilon} \sheaf{\tilde T^\vee}; \\
\sheaf{\hat T}_{\SC, {[n]}} = \shom(\sheaf{Y}_{Q,n}^{\SC}, \mu_n) & \xrightarrow{\epsilon} \sheaf{\tilde T}_{\SC}^\Vee; \\
\sheaf{\hat Z}_{[n]} = \shom(\sheaf{Y}_{Q,n} / \sheaf{Y}_{Q,n}^{\SC}, \mu_n) & \xrightarrow{\epsilon} \sheaf{\tilde Z^\vee}.
\end{align*}

\subsection{The gerbe associated to a cover of a torus}

Associated to the cover $\alg{\tilde T} = (\alg{T}', n)$, the second Brylinski-Deligne invariant is a central extension of sheaves of groups on $S_{\et}$,
$$\sheaf{G}_m \Into \sheaf{D} \Onto \sheaf{Y}.$$
The commutator of this extension is given in \cite[Proposition 3.11]{B-D},
\begin{equation}
\label{CommForm}
\Comm(y_1, y_2) = (-1)^{n \beta_Q(y_1, y_2)}, \text{ for all } y_1, y_2 \in \sheaf{Y}.
\end{equation}
Pulling back via $\sheaf{Y}_{Q,n} \Into \sheaf{Y}$, we find an extension of sheaves of groups,
\begin{equation}
\label{BD2Qn}
\sheaf{G}_m \Into \sheaf{D}_{Q,n} \Onto \sheaf{Y}_{Q,n}.
\end{equation}

\begin{proposition}
$\sheaf{D}_{Q,n}$ is a commutative extension.
\end{proposition}
\proof
If $n$ is even and $y_1, y_2 \in \sheaf{Y}_{Q,n}$, then $\beta_Q(y_1, y_2) \in \ZZ$ and $n \beta_Q(y_1, y_2) \in 2 \ZZ$.  On the other hand, if $n$ is odd, Assumption \ref{OddnEvenQ} implies that 
$$n \beta_Q(y_1, y_2) = Q(y_1 + y_2) - Q(y_1) - Q(y_2) \in 2 \ZZ.$$
The commutator formula \eqref{CommForm} finishes the proof.
\qed

Let $\sspl(\sheaf{D}_{Q,n})$ denote the sheaf of {\em splittings} of the commutative extension \eqref{BD2Qn}.  In other words, $\sspl(\sheaf{D}_{Q,n})$ is the subsheaf of $\shom(\sheaf{Y}_{Q,n}, \sheaf{D}_{Q,n})$ consisting of homomorphisms which split \eqref{BD2Qn}.

Over any finite \'etale $U \To S$ splitting $\alg{T}$, $\sheaf{Y}_{Q,n}$ restricts to a constant sheaf of free abelian groups.  Thus $\sspl(\sheaf{D}_{Q,n})$ is a $\sheaf{\hat T}$-torsor on $S_{\et}$, which obtains a point over any such $U$.  The equivalence class of this torsor is determined by its cohomology class $\left[ \sspl(\sheaf{D}_{Q,n}) \right] \in H_{\et}^1(S, \sheaf{\hat T})$.

Consider the Kummer sequence, $\sheaf{\hat T}_{[n]} \Into \sheaf{\hat T}\xtwoheadrightarrow{n} \sheaf{\hat T}$.  Write $\Kum$ (for Kummer) for the connecting map in cohomology, $\Kum \From H_{\et}^1(S, \sheaf{\hat T}) \To H_{\et}^2(S, \sheaf{\hat T}_{[n]})$.  This map in cohomology corresopnds to the functor which sends a $\sheaf{\hat T}$-torsor to its gerbe of $n^{\th}$ roots (see \ref{AppendixGerbeRoots} for details).  We write
$\sqrt[n]{\sspl(\sheaf{D}_{Q,n})}$ for the gerbe of $n^{\th}$ roots of the $\sheaf{\hat T}$-torsor $\sspl(\sheaf{D}_{Q,n})$.  It is banded by the local system $\sheaf{\hat T}_{[n]}$ and its equivalence class satisfies
$$\left[ \sqrt[n]{\sspl(\sheaf{D}_{Q,n})} \right] = \Kum [\sspl(\sheaf{D}_{Q,n})].$$
Finally we push out via the homomorphism of local systems,
$$\epsilon \From \sheaf{\hat T}_{[n]} = \shom(\sheaf{Y}_{Q,n}, \mu_n) \To \shom(\sheaf{Y}_{Q,n}, \CC^\times) = \sheaf{\tilde T^\vee}.$$

\begin{definition}
The \defined{gerbe associated to the cover} $\alg{\tilde T}$ is defined by
$$\gerb{E}_\epsilon(\alg{\tilde T}) \defeq \epsilon_\ast \sqrt[n]{\sspl(\sheaf{D}_{Q,n})}.$$
It is a gerbe on $S_{\et}$ banded by the local system of abelian groups $\sheaf{\tilde T^\vee}$.
\end{definition}

\begin{example}
\label{SplitTorusGerbeNeutral}
Suppose that $\alg{T}$ is a split torus.  Then the exact sequence of sheaves $\sheaf{G}_m \Into \sheaf{D}_{Q,n} \Onto \sheaf{Y}_{Q,n}$ splits.  Indeed, $\sheaf{Y}_{Q,n}$ is a constant sheaf of free abelian groups, and Hilbert's Theorem 90 gives a short exact sequence
\begin{equation}
\label{SPointsD}
\sheaf{G}_m[S] \Into \sheaf{D}_{Q,n}[S] \Onto \sheaf{Y}_{Q,n}[S].
\end{equation}
Since $\sheaf{Y}_{Q,n}[S]$ is a free abelian group, this exact sequence splits, and any such splitting defines an $S$-point of the torsor $\sspl(\sheaf{D}_{Q,n})$.  An $S$-point of $\sspl(\sheaf{D}_{Q,n})$, in turn, neutralizes of the gerbe $\sqrt[n]{\sspl(\sheaf{D}_{Q,n})}$.

Thus when $\alg{T}$ is a split torus, the gerbe $\gerb{E}_\epsilon(\alg{\tilde T})$ is trivial.  Any splitting of the sequence \eqref{SPointsD} defines a neutralization of $\gerb{E}_\epsilon(\alg{\tilde T})$.
\end{example}  

\subsection{The gerbe of liftings}

Recall that $\sheaf{Y}_{Q,n}^{\SC}$ denotes the subgroup of $\sheaf{Y}_{Q,n}$ spanned by the modified coroots $\tilde \Phi^\vee$, and $\sheaf{\hat T}_{\SC} = \shom(\sheaf{Y}_{Q,n}^{\SC}, \sheaf{G}_m)$.  The inclusion $\sheaf{Y}_{Q,n}^{\SC} \Into \sheaf{Y}_{Q,n}$ corresponds to a surjective homomorphism,
$$p \From \sheaf{\hat T} \To \sheaf{\hat T}_{\SC}.$$
The extension $\sheaf{G}_m \Into \sheaf{D}_{Q,n} \Onto \sheaf{Y}_{Q,n}$ pulls back via $\sheaf{Y}_{Q,n}^{\SC} \Into \sheaf{Y}_{Q,n}$ to an extension,
$$\sheaf{G}_m \Into \sheaf{D}_{Q,n}^{\SC} \Onto \sheaf{Y}_{Q,n}^{\SC}.$$
A splitting of $\sheaf{D}_{Q,n}$ pulls back to a splitting of $\sheaf{D}_{Q,n}^{\SC}$, providing a map of torsors,
$$p^\ast \From \sspl(\sheaf{D}_{Q,n}) \To \sspl(\sheaf{D}_{Q,n}^{\SC}),$$
lying over $p \From \sheaf{\hat T} \To \sheaf{\hat T}_{\SC}$.  Taking $n^{\th}$ roots of torsors gives a functor of gerbes,
$$\sqrt[n]{p^\ast}\From \sqrt[n]{\sspl(\sheaf{D}_{Q,n})} \To \sqrt[n]{\sspl(\sheaf{D}_{Q,n}^{\SC})},$$
lying over $p \From \sheaf{\hat T}_{[n]} \To \sheaf{\hat T}_{\SC, [n]}$ (see \ref{AppendixGerbeRoots}).

Recall that $\sheaf{\tilde T^\vee} = \shom(\sheaf{Y}_{Q,n}, \CC^\times)$ and $\sheaf{\tilde T}_{\SC}^\Vee = \shom(\sheaf{Y}_{Q,n}^{\SC}, \CC^\times)$.  Define $\gerb{E}_\epsilon^{\SC}(\alg{\tilde T}) \defeq \epsilon_\ast \sqrt[n]{\sspl(\sheaf{D}_{Q,n}^{\SC})}$ by analogy to $\gerb{E}_\epsilon(\alg{\tilde T}) = \epsilon_\ast \sqrt[n]{\sspl(\sheaf{D}_{Q,n})}$.  Pushing out via $\epsilon$, the functor $\sqrt[n]{p^\ast}$ yields a functor of gerbes
$$\gerb{p} = \epsilon_\ast \sqrt[n]{p^\ast} \From \gerb{E}_\epsilon(\alg{\tilde T}) \To \gerb{E}_\epsilon^{\SC}(\alg{\tilde T}),$$
lying over the homomorphism $p \From \sheaf{\tilde T^\vee} \To \sheaf{\tilde T}_{\SC}^\Vee$.

In the next section, we define the {\em Whittaker torsor}, which gives an object $\whit$ neutralizing the gerbe $\gerb{E}_\epsilon^{\SC}(\alg{\tilde T})$.  We take this construction of $\whit$ for granted at the moment.

\begin{definition}
Define $\gerb{E}_\epsilon(\alg{\tilde G})$ to be the gerbe $\gerb{p}^{-1}(\whit)$ of liftings of $\whit$ via $\gerb{p}$ (see \ref{AppendixGerbeLiftings}).  In other words, $\gerb{E}_\epsilon(\alg{\tilde G})$ is the 
category of pairs $(\gerb{e}, j)$ where $\gerb{e}$ is an object of $\gerb{E}_\epsilon(\alg{\tilde T})$ and $j \From \gerb{p}(\gerb{e}) \To \whit$ is an isomorphism in $\gerb{E}_\epsilon^{\SC}(\alg{\tilde T})$. This is a gerbe on $S_{\et}$ banded by $\sheaf{\tilde Z^\vee} = \Ker(\sheaf{\tilde T^\vee} \xrightarrow{p} \sheaf{\tilde T}_{\SC}^\Vee)$.
\end{definition}

The cohomology classes of our gerbes now fit into a sequence
$$\begin{tikzcd}[row sep = 0em]
\left[ \gerb{E}_\epsilon(\alg{\tilde G}) \right] \arrow[mapsto]{r} & \left[ \gerb{E}_\epsilon(\alg{\tilde T}) \right] \arrow[mapsto]{r} & \left[ \gerb{E}_\epsilon^{\SC}(\alg{\tilde T}) \right] = 0 \\
\rotatebox{-90}{$\in$} & \rotatebox{-90}{$\in$}& \rotatebox{-90}{$\in$} \\
H_{\et}^2(S, \sheaf{\tilde Z^\vee}) \arrow{r} & H_{\et}^2(S, \sheaf{\tilde T^\vee}) \arrow{r} &H_{\et}^2(S, \sheaf{\tilde T_{\SC} ^\vee}) 
\end{tikzcd}$$

\begin{remark}
The construction of this gerbe relies on the (soon-to-be-defined) Whittaker torsor in a crucial way.  We view this as a good thing, since any putative Langlands correspondence should also connect the existence of Whittaker models to properties of the Langlands parameter (cf. \cite{VoganLLC}).
\end{remark}

\subsection{The Whittaker torsor}

Now we construct the object $\whit$ neutralizing the gerbe $\gerb{E}_\epsilon^{\SC}(\alg{\tilde T})$ over $S$.  Let $\alg{U}$ denote the unipotent radical of the Borel subgroup $\alg{B} \subset \alg{G}$, and let $\sheaf{U}$ be the sheaf of groups on $S_{\et}$ that it represents.  Let $\alg{G}_a$ denote the additive group scheme over $S$, and $\sheaf{G}_a$ the sheaf of groups on $S_{\et}$ that it represents.  Recall that $\Delta \subset \Phi$ denotes the subset of simple roots corresponding to the Borel subgroup $\alg{B}$.  

For $S' \To S$ finite \'etale and splitting $\alg{T}$, and $\alpha \in \Delta[S']$, write $\alg{U}_\alpha$ for the one-dimensional root subgroup of $\alg{U}_{S'}$ associated to $\phi$.  Let $\sheaf{U}_\alpha$ be the associated sheaf of abelian groups on $S_{\et}'$.  Write $\shom^\ast(\sheaf{U}_\alpha, \sheaf{G}_{a})$ for the sheaf (on $S_{\et}'$) of isomorphisms from $\sheaf{U}_\alpha$ to $\sheaf{G}_{a,S'}$.  The sheaf $\shom^\ast(\sheaf{U}_\alpha, \sheaf{G}_{a})$ naturally forms a $\sheaf{G}_m$-torsor on $S_{\et}'$, by the formula
$$[h \ast \xi](u) = h^{-1} \cdot \xi(u) \text{ for all } h \in \sheaf{G}_m, \xi \in \shom^\ast(\sheaf{U}_\alpha, \sheaf{G}_a).$$
\begin{definition}
The \defined{Whittaker torsor} is the subsheaf $\Whit \subset \shom(\sheaf{U}, \sheaf{G}_a)$ consisting of those homomorphisms which (locally on $S_{\et}$)  restrict to an isomorphism on every simple root subgroup.  The sheaf $\Whit$ is given the structure of a $\sheaf{\hat T}_{\SC}$-torsor as follows:  for a Galois cover $S' \To S$ splitting $\alg{T}$, we have
$$\sheaf{\hat T}_{\SC}[S'] = \shom(\sheaf{Y}_{Q,n}^{\SC}, \sheaf{G}_m)[S'] = \shom \left( \bigoplus_{\alpha \in \Delta[S']} \ZZ {\tilde \alpha}^\vee, \sheaf{G}_{m} \right) [S'] \ident \prod_{\alpha \in \Delta[S']} \sheaf{G}_{m}[S'].$$
Similarly, we can decompose the Whittaker sheaf
$$\Whit[S'] \ident \bigoplus_{\alpha \in \Delta[S']} \shom^\ast(\sheaf{U}_\alpha, \sheaf{G}_{a})[S'].$$
The $\sheaf{G}_m$-torsor structure on $\shom^\ast(\sheaf{U}_\alpha, \sheaf{G}_{a})$ yields (simple root by simple root) a $\sheaf{\hat T}_{\SC}$-torsor structure on $\Whit$.  Although we have defined the torsor structure locally on $S_{\et}$, the action descends since the $\Gal(S'/S)$-actions are compatible throughout.
\end{definition}

The third Brylinski-Deligne invariant of $\alg{\tilde G}$ is a homomorphism $f \From \sheaf{D}_Q \To \sheaf{D}$ of groups on $S_{\et}$.
$$\begin{tikzcd}
\sheaf{G}_m \inarrow{r} \arrow{d}{=} &  \sheaf{D}_Q \onarrow{r} \inarrow{d}{f} & \sheaf{Y}^{\SC} \inarrow{d} \\
\sheaf{G}_m \inarrow{r} & \sheaf{D} \onarrow{r} & \sheaf{Y}
\end{tikzcd}$$
Here $\sheaf{D}_Q$ is a sheaf on $S_{\et}$ which depends (up to unique isomorphism) only on the Weyl- and Galois-invariant quadratic form $Q \From \sheaf{Y}^{\SC} \To \ZZ$.  This is reviewed in \cite[\S 1.3, 3.3]{MWIntegral}, and characterized in \cite[\S 11]{B-D} when working over a field.

Consider a Galois cover $S' \To S$ splitting $\alg{T}$ as before.  For any $\eta \in \Whit[S']$, and any simple root $\alpha \in \Delta[S']$, there exists a unique element $e_{\eta, \alpha} \in \sheaf{U}_\alpha[S']$ such that $\eta(e_{\eta,\alpha}) = 1$.  From these, \cite[\S 11.2]{B-D} gives elements $[e_{\eta,\alpha}] \in \sheaf{D}_Q[S']$ lying over the simple coroots $\alpha^\vee \in \sheaf{Y}^{\SC}[S']$.
\begin{remark}
When $S = \Spec(F)$ this follows directly from \cite[\S 11.2]{B-D}.  When $S = \Spec(\OO)$, $S' = \Spec(\OO')$, $\eta \in \Whit[S']$, and $F'$ is the fraction field of $\OO'$, we find elements $e_{\eta,\alpha} \in \sheaf{U}_\alpha[F']$; as $\eta$ gives an isomorphism from $\sheaf{U}_\alpha$ to $\sheaf{G}_a$ (as sheaves of groups on $\OO_{\et}'$), it follows that $e_{\eta,\alpha} \in \sheaf{U}_\alpha[\OO']$ as well.  The map $e \mapsto [e]$ of \cite[\S 11.1]{B-D} similarly makes sense over $\OO'$ as well as it does over a field; since we assume $\alg{G}$ is a reductive group over $\OO$, split over $\OO'$, every root $SL_2$ over $F'$ arises from one over $\OO'$.  Thus the results of \cite[\S 11.2]{B-D} apply in the setting of $S = \Spec(\OO)$ as well as in the setting of a field.
\end{remark}

Using the elements $[e_{\eta,\alpha}] \in \sheaf{D}_Q[S']$ lying over the simple coroots $\alpha^\vee$, define
$$\omega(\eta)(\tilde \alpha^\vee) \defeq f([e_{\eta,\alpha}])^{n_\alpha} \in \sheaf{D}_{Q,n}^{\SC}[S'], \text{ lying over } \tilde \alpha^\vee = n_\alpha \alpha^\vee \in \sheaf{Y}_{Q,n}^{\SC}[S'].$$
The map $\omega(\eta) \From \tilde \alpha^\vee \mapsto f([e_{\eta, \alpha}])^{n_\alpha}$ extends uniquely to a splitting of the sequence
\begin{equation}
\begin{tikzcd}
\sheaf{G}_m[S'] \inarrow{r} & \sheaf{D}_{Q,n}^{\SC}[S'] \onarrow{r} & \sheaf{Y}_{Q,n}^{\SC}[S'] \arrow[bend right=20]{l}[swap]{\omega(\eta)}.
\end{tikzcd}
\end{equation}
As $(\sheaf{Y}_{Q,n}^{\SC})_{S'}$ is a constant sheaf, this gives an element $\omega(\eta) \in \sspl(\sheaf{D}_{Q,n}^{\SC})[S']$.  Allowing $\eta$ to vary, and appyling Galois descent (cf.~\cite[Proposition 11.7]{B-D}), we find a map of sheaves on $S_{\et}$,
$$\omega \From \Whit \To \sspl(\sheaf{D}_{Q,n}^{\SC}).$$
To summarize, $\omega$ is the map that sends a nondegenerate character $\eta$ of $\alg{U}$ to the splitting $\omega(\eta)$, which (locally on $S_{\et}$) sends each modified simple coroot $\tilde \alpha^\vee$ to the element $f([e_{\eta, \alpha}])^{n_\alpha}$ of $\sheaf{D}_{Q,n}$. 

Both $\Whit$ and $\sspl(\sheaf{D}_{Q,n}^{\SC})$ are $\sheaf{\hat T}_{\SC}$-torsors, and the following proposition describes how $\omega$ interacts with the torsor structure.
\begin{proposition}
Let $\nu \From \sheaf{\hat T}_{\SC} \To \sheaf{\hat T}_{\SC}$ be the homomorphism corresponding to the unique homomorphism $\sheaf{Y}_{Q,n}^{\SC} \To \sheaf{Y}_{Q,n}^{\SC}$ which sends $\tilde \alpha^\vee$ to $- n_\alpha Q(\alpha^\vee) \tilde \alpha^\vee$ for all simple roots $\alpha$.  Then $\omega$ lies over $\nu$, i.e., the following diagram commutes.
$$\begin{tikzcd}
\sheaf{\hat T_{\SC}} \times \Whit \arrow{r}{\ast} \arrow{d}{\nu \times \omega} & \Whit \arrow{d}{\omega} \\
\sheaf{\hat T_{\SC}} \times \sspl(\sheaf{D}_{Q,n}^{\SC}) \arrow{r}{\ast} & \sspl(\sheaf{D}_{Q,n}^{\SC})
\end{tikzcd}$$
\end{proposition}
\proof
We must trace through the action of $\sheaf{\hat T}_{\SC} = \shom(\sheaf{Y}_{Q,n}^{\SC}, \sheaf{G}_m)$; we work over a finite \'etale cover of $S$ over which $\alg{T}$ splits in what follows.  Then, for any simple root $\alpha \in \Delta$, and any element $h \in \sheaf{G}_m$, there exists a unique element $h_\alpha \in \sheaf{\hat T}_{\SC}$ such that for all $\beta \in \Delta$,
$$h_\alpha(\tilde \beta^\vee) = \begin{cases} 1 & \text{ if } \beta \neq \alpha; \\ h & \text{ if } \beta = \alpha.  \end{cases}$$
If $\eta \in \Whit$ then $[h_\alpha \ast \eta](e_{h_\alpha \ast \eta, \alpha}) = 1$ and so $\eta(e_{h_\alpha \ast \eta, \alpha}) = h$.  Therefore,
$$e_{h_\alpha \ast \eta, \beta} = \begin{cases} e_{\eta, \beta} & \text{ if } \beta \neq \alpha; \\ h \ast  e_{\eta, \alpha} & \text{ if } \beta = \alpha. \end{cases}$$
If $\beta \neq \alpha$, then $\omega(h_\alpha \ast \eta)(\tilde \beta^\vee) = f([e_{h_\alpha \ast \eta, \beta}])^{n_\beta} = f([e_{\eta, \beta}])^{n_\beta} = \omega(\eta)(\tilde \beta^\vee)$.  On the other hand, in the case $\beta = \alpha$ we compute using \cite[Equation (11.2.1)]{B-D},
\begin{align*}
\omega(h_\alpha \ast \eta)(\tilde \alpha^\vee) = f([e_{h_\alpha \ast \eta, \alpha}])^{n_\alpha} &= f \left( [h \ast e_{\eta, \alpha} ] \right)^{n_\alpha} \\
&= f \left( h^{-Q(\alpha^\vee)} \cdot [e_{\eta, \alpha} ] \right)^{n_\alpha} \\
&= h^{ - n_\alpha Q(\phi_\alpha^\vee) } \cdot f \left( [e_{\eta, \alpha} ] \right)^{n_\alpha} = h^{ - n_\alpha Q(\phi_\alpha^\vee) } \cdot \omega(\eta)(\tilde \alpha^\vee).
\end{align*}
This computation demonstrates that the diagram commutes as desired.
\qed

%\begin{remark}
%Equation (11.2.1) of \cite{B-D} has a negative sign in the exponent.  Our Whittaker torsor is dual to (a product of the) torsors $U_\alpha^\ast$ discussed in \cite[\S 11.2]{B-D}, which cancels the negative sign in the exponent in the computation above.
%\end{remark}

Now let $\mu \From \sheaf{\hat T}_{\SC} \Onto \sheaf{\hat T}_{\SC}$ be the homomorphism corresponding to the unique homomorphism $\sheaf{Y}_{Q,n}^{\SC} \Into \sheaf{Y}_{Q,n}^{\SC}$ which sends $\tilde \alpha^\vee$ to $- m_\alpha \tilde \alpha^\vee$ for all $\alpha \in \Delta$.  As $Q(\alpha^\vee) n_\alpha = m_\alpha \cdot n$, we find that $\nu = n \circ \mu$, where $n$ denotes the $n^{\th}$-power map.

Let $\mu_\ast \Whit$ denote the pushout of the $\sheaf{\hat T}_{\SC}$-torsor $\Whit$, via $\mu$.  Since $\nu$ factors through $\mu$, we find that $\omega \From \Whit \To \sspl(\sheaf{D}_{Q,n}^{\SC})$ factors uniquely through $\bar \omega \From \mu_\ast \Whit \To \sspl(\sheaf{D}_{Q,n}^{\SC})$, making the following diagram commute.
$$\begin{tikzcd}
\sheaf{\hat T_{\SC}} \times \left( \mu_\ast \Whit \right) \arrow{r}{\ast} \arrow{d}{n \times \bar \omega} & \mu_\ast \Whit \arrow{d}{\bar \omega}  \\
\sheaf{\hat T_{\SC}} \times \sspl(\sheaf{D}_{Q,n}^{\SC}) \arrow{r}{\ast} & \sspl(\sheaf{D}_{Q,n}^{\SC})
\end{tikzcd}$$

The pair $(\mu_\ast \Whit, \bar \omega)$ is therefore an object of the category $\sqrt[n]{\sspl(\sheaf{D}_{Q,n}^{\SC})}[S]$; it {\em neutralizes} the gerbe $\sqrt[n]{\sspl(\sheaf{D}_{Q,n}^{\SC})}$.  In particular,
$$\left[ \sqrt[n]{\sspl(\sheaf{D}_{Q,n}^{\SC})} \right] = 0.$$   
Write $\whit = (\mu_\ast \Whit, \bar \omega)$ for this object.  Pushing out via $\epsilon$, we view $\whit$ as an $S$-object of $\gerb{E}_\epsilon^{\SC}(\alg{\tilde T})$.  This completes the construction of the gerbe $\gerb{E}_\epsilon(\alg{\tilde G}) = \gerb{p}^{-1}(\gerb{w})$ associated to the cover $\alg{\tilde G}$ and character $\epsilon$.

\begin{example}
Suppose that $\sheaf{Y}_{Q,n} / \sheaf{Y}_{Q,n}^{\SC}$ is torsion-free and a constant sheaf (equivalently, the center of $\dgp{\tilde G^\vee}$ is connected and constant as a sheaf on $S_{\et}$).  Then the following short exact sequence splits:
$$\sheaf{Y}_{Q,n}^{\SC} \Into \sheaf{Y}_{Q,n} \Onto \sheaf{Y}_{Q,n} /    \sheaf{Y}_{Q,n}^{\SC}.$$

Given such a splitting, write $\sheaf{Y}_{Q,n}^{\cent} \subset \sheaf{Y}_{Q,n}$ for the image of $\sheaf{Y}_{Q,n} / \sheaf{Y}_{Q,n}^{\SC}$ via the splitting.  The identification $\sheaf{Y}_{Q,n} = \sheaf{Y}_{Q,n}^{\SC} \oplus \sheaf{Y}_{Q,n}^{\cent}$ corresponds to an isomorphism $\sheaf{\hat T} \xrightarrow{\sim} \sheaf{\hat T}_{\SC} \times \sheaf{\hat Z}$.   Let $\sheaf{D}_{Q,n}^{\cent}$ be the pullback of $\sheaf{D}_{Q,n}$ to $\sheaf{Y}_{Q,n}^{\cent}$.  From Example \ref{SplitTorusGerbeNeutral}, the short exact sequence $\sheaf{G}_m \Into \sheaf{D}_{Q,n}^{\cent} \Onto \sheaf{Y}_{Q,n}^{\cent}$ splits, providing an object of $\sqrt[n]{\sspl(\sheaf{D}_{Q,n}^{{\cent}})}$.  

Chasing diagrams gives a map of object sets,
$$\sqrt[n]{\sspl(\sheaf{D}_{Q,n}^{\SC})}[S] \times \sqrt[n]{\sspl(\sheaf{D}_{Q,n}^{{\cent}})}[S] \To \sqrt[n]{\sspl(\sheaf{D}_{Q,n})}[S].$$
A splitting of $\sheaf{D}_{Q,n}^{\cent}$ gives an object of $\sqrt[n]{\sspl(\sheaf{D}_{Q,n}^{{\cent}})}[S]$ and $\whit$ provides an object of $\sqrt[n]{\sspl(\sheaf{D}_{Q,n}^{\SC})}[S]$.  Hence the gerbe $\gerb{E}_\epsilon(\alg{\tilde G})$ is neutral when $\sheaf{Y}_{Q,n} / \sheaf{Y}_{Q,n}^{\SC}$ is torsion-free and a constant sheaf.
\end{example}

\subsection{Well-aligned functoriality}
\label{WAFGerbe}

Consider a well-aligned homomorphism $\tilde \iota \From \alg{\tilde G}_1 \To \alg{\tilde G}_2$ of covers, each endowed with Borel subgroup and maximally split maximal torus, i.e., a morphism in the category $\Cat{WAC}_S$.  Fix $\epsilon$ as before.  We have constructed gerbes $\gerb{E}_\epsilon(\alg{\tilde G}_1)$ and $\gerb{E}_\epsilon(\alg{\tilde G}_2)$ associated to $\alg{\tilde G}_1$ and $\alg{\tilde G}_2$, banded by $\sheaf{\tilde Z}_1^\Vee$ and $\sheaf{\tilde Z}_2^\Vee$, respectively.  We have constructed a homomorphism of dual groups $\iota^\vee \From \dgp{\tilde G}_2^\Vee \To \dgp{\tilde G}_1^\Vee$ in Section \ref{WAFDualGroup}, which (after taking $\CC$-points) restricts to $\iota^\vee \From \sheaf{\tilde Z}_2^\Vee \To \sheaf{\tilde Z}_1^\Vee$.  Here we construct a functor of gerbes $\gerb{i} \From \gerb{E}_\epsilon(\alg{\tilde G}_2) \To \gerb{E}_\epsilon(\alg{\tilde G}_1)$, lying over $\iota^\vee \From \sheaf{\tilde Z}_{2}^\Vee \To \sheaf{\tilde Z}_{1}^\Vee$.

Well-alignedness give a commutative diagram in which the first row is the pullback of the second.
$$\begin{tikzcd}
\alg{K}_2 \inarrow{r} \arrow{d}{=} & \alg{T}_1' \onarrow{r} \arrow{d}{\iota'} & \alg{T}_1 \arrow{d}{\iota} \\
\alg{K}_2 \inarrow{r} & \alg{T}_2' \onarrow{r} & \alg{T}_2
\end{tikzcd}$$
This gives a commutative diagram for the second Brylinski-Deligne invariants.  After pulling back to $\sheaf{Y}_{1,Q_1,n}$ and $\sheaf{Y}_{2,Q_2,n}$, we get a commutative diagram of sheaves of abelian groups on $S_{\et}$, in which the top row is the pullback of the bottom.
$$\begin{tikzcd}
\sheaf{G}_m \inarrow{r} \arrow{d}{=} & \sheaf{D}_{1,Q_1,n} \onarrow{r} \arrow{d}{\iota'} & \sheaf{Y}_{1,Q_1,n} \arrow{d}{\iota} \\
\sheaf{G}_m \inarrow{r} & \sheaf{D}_{2,Q_2,n} \onarrow{r}& \sheaf{Y}_{2,Q_2,n}
\end{tikzcd}$$
(Assumption \ref{OddnEvenQ} is in effect, so $Q_1$ and $Q_2$ are even-valued if $n$ is odd.)

We have homomorphisms of sheaves of abelian groups,
$$\iota \From \sheaf{Y}_{1,Q_1,n} \To \sheaf{Y}_{2,Q_2,n}, \quad \sheaf{Y}_{1,Q_1,n}^{\SC} \To \sheaf{Y}_{2,Q_2,n}^{\SC}, \quad \frac{\sheaf{Y}_{1,Q_1,n}}{\sheaf{Y}_{1,Q_1,n}^{\SC}} \To \frac{\sheaf{Y}_{2,Q_2,n}}{\sheaf{Y}_{2,Q_2,n}^{\SC}}.$$
Applying $\shom(\bullet, \sheaf{G}_m)$ yields homomorphisms of sheaves of abelian groups,
$$\hat \iota \From \sheaf{\hat T}_2 \To \sheaf{\hat T}_1, \quad \sheaf{\hat T}_{\SC,2} \To \sheaf{\hat T}_{\SC,1}, \quad \sheaf{\hat Z}_2 \To \sheaf{\hat Z}_1.$$

A splitting of $\sheaf{D}_{2,Q_2,n}$ pulls back to a splitting of $\sheaf{D}_{1,Q_1,n}$ (call this pullback map $\iota^\ast$), giving a commutative diagram of sheaves on $S_{\et}$.
$$\begin{tikzcd}
\sheaf{\hat T}_2 \times \sspl(\sheaf{D}_{2,Q_2,n}) \arrow{r}{\ast} \arrow{d}{\hat \iota \times \iota^\ast} & \sspl(\sheaf{D}_{2,Q_2,n}) \arrow{d}{\iota^\ast} \\
\sheaf{\hat T}_1 \times \sspl(\sheaf{D}_{1,Q_1,n}) \arrow{r}{\ast} & \sspl(\sheaf{D}_{1,Q_1,n})
\end{tikzcd}$$
This defines a functor of gerbes 
$$\sqrt[n]{\iota^\ast} \From \sqrt[n]{\sspl(\sheaf{D}_{2,Q_2,n})} \To \sqrt[n]{\sspl(\sheaf{D}_{1,Q_1,n})},$$
lying over $\hat \iota \From \sheaf{\hat T}_{2,[n]} \To \sheaf{\hat T}_{1,[n]}$.  Pushing out via $\epsilon$ yields a functor of gerbes,
$$\gerb{i} \From \gerb{E}_\epsilon(\alg{\tilde T}_2) \To \gerb{E}_\epsilon(\alg{\tilde T}_1).$$
The same process applies to $\iota \From \sheaf{Y}_{1,Q_1,n}^{\SC} \To \sheaf{Y}_{2,Q_2,n}^{\SC}$, giving a functor of gerbes, $\gerb{i}^{\SC} \From \gerb{E}_\epsilon^{\SC}(\alg{\tilde T}_2) \To \gerb{E}_\epsilon^{\SC}(\alg{\tilde T}_1)$.  By pulling back in stages, we find a square of gerbes and functors, and a natural isomorphism $\Fun{S} \From \gerb{p}_1 \circ \gerb{i} \xRightarrow{\sim} \gerb{i}^{\SC} \circ \gerb{p}_2$ making the diagram 2-commute.
$$\begin{tikzcd}
 \gerb{E}_\epsilon(\alg{\tilde T}_2) \arrow{r}{\gerb{i}} \arrow{d}{\gerb{p}_2} & \gerb{E}_\epsilon(\alg{\tilde T}_1) \arrow{d}{\gerb{p}_1} \\
\gerb{E}_\epsilon^{\SC}(\alg{\tilde T}_2)  \arrow{r}{\gerb{i}^{\SC}} & \gerb{E}_\epsilon^{\SC}(\alg{\tilde T}_1)
\end{tikzcd}$$

If $\alg{U}_1$ and $\alg{U}_2$ are the unipotent radicals of $\alg{B}_1$ and $\alg{B}_2$, respectively, then pulling back gives a map $\iota^\ast \From \Whit_2 \To \Whit_1$.  Condition (1) of well-alignedness states that $\Ker(\iota)$ is contained in the center of $\alg{G}_1$, and so simple root subgroups in $\alg{U}_1$ map isomorphically to simple root subgroups in $\alg{U}_2$.  Compatibility of quadratic forms $Q_1$ and $Q_2$ implies compatibility in the constants $n_\alpha$ and $m_\alpha$ for simple roots, and so we find a commutative diagram
$$\begin{tikzcd}
(\mu_2)_\ast \Whit_2 \arrow{r}{\bar \omega_2} \arrow{d}{\iota^\ast} & \sspl(\sheaf{D}_{2,Q_2,n}^{\SC}) \arrow{d}{\iota^\ast} \\
(\mu_1)_\ast \Whit_1 \arrow{r}{\bar \omega_1} & \sspl(\sheaf{D}_{1,Q_1,n}^{\SC})  
\end{tikzcd}$$
Write $\whit_1$ for the object $((\mu_1)_\ast \Whit_1, \bar \omega_1)$ of $\gerb{E}_\epsilon^{\SC}(\alg{\tilde T}_1)$, and similarly $\whit_2$ for the object $((\mu_2)_\ast \Whit_2, \bar \omega_2)$ of $\gerb{E}_\epsilon^{\SC}(\alg{\tilde T}_2)$.  The commutative diagram above gives an isomorphism $f \From \gerb{i}^{\SC}( \whit_2) \xrightarrow{\sim} \whit_1$ in $\gerb{E}_\epsilon^{\SC}(\alg{\tilde T}_1)$.

If $(\gerbob{e},j)$ is an object of $\gerb{E}_\epsilon(\alg{\tilde G}_2) = \gerb{p}_2^{-1}(\whit_2)$, i.e., $\gerbob{e}$ is an object of $\gerb{E}_\epsilon(\alg{\tilde T}_2)$ and $j \From \gerb{p}_2(\gerbob{e}) \To \whit_2$ is an isomorphism, then we find a sequence of isomorphisms,
$$\gerb{i}(j) \From \gerb{p}_1(\gerb{i}(\gerbob{e})) \xrightarrow{\Fun{S}} \gerb{i}^{\SC}(\gerb{p}_2(\gerbob{e})) \xrightarrow{j} \gerb{i}^{\SC}(\whit_2) \xrightarrow{f} \whit_1.$$
In this way $(\gerb{i}(\gerbob{e}),\gerb{i}(j))$ becomes an object of $\gerb{E}_\epsilon(\alg{\tilde G}_1) = \gerb{p}_1^{-1}(\whit_1)$.  This extends to a functor of gerbes, $\gerb{i} \From \gerb{E}_\epsilon(\alg{\tilde G}_2) \To \gerb{E}_\epsilon(\alg{\tilde G}_1)$, lying over $\iota^\vee \From \sheaf{\tilde Z}_2^\Vee \To \sheaf{\tilde Z}_1^\Vee$.  

Given a pair of well-aligned homomorphisms,
$$\alg{\tilde G}_1 \xrightarrow{\tilde \iota_1} \alg{\tilde G}_2 \xrightarrow{\tilde \iota_2} \alg{\tilde G}_2,$$
with $\tilde \iota_3 = \tilde \iota_2 \circ \tilde \iota_1$, we find three functors of gerbes.
$$\begin{tikzcd}
\gerb{E}_\epsilon(\alg{\tilde G}_3) \arrow{r}[swap]{\gerb{i}_2} \arrow[bend left=20]{rr}{\gerb{i}_3} & \gerb{E}_\epsilon(\alg{\tilde G}_2) \arrow{r}[swap]{\gerb{i}_1} & \gerb{E}_\epsilon(\alg{\tilde G}_1).
\end{tikzcd}$$
lying over three homomorphisms of sheaves,
$$\begin{tikzcd}
\sheaf{\tilde Z}_3^\Vee \arrow{r}[swap]{\iota_2^\vee} \arrow[bend left=20]{rr}{\iota_3^\vee} & \sheaf{\tilde Z}_2^\Vee \arrow{r}[swap]{\iota_1^\vee} & \sheaf{\tilde Z}_1^\Vee.
\end{tikzcd}$$

The functoriality of pullbacks, pushouts, taking $n^{\th}$ roots of torsors, etc., yields a \textbf{natural isomorphism} of functors:
\begin{equation}
\label{CompositionWellAligned}
N(\iota_1, \iota_2) \From \gerb{i}_3 \xRightarrow{\sim} \gerb{i}_1 \circ \gerb{i}_2.
\end{equation}

We can summarize these results using the language of ``weak functors'' \cite[Tag 003G]{Stacks}.  We have given all the data for such a weak functor, and a reader who wishes to check commutativity of a few diagrams may verify the following. 
\begin{proposition}
The construction of the associated gerbe defines a weak functor from $\Cat{WAC}_S$ (the category of triples $(\alg{\tilde G}, \alg{B}, \alg{T})$ and well-aligned homomorphisms),  to the 2-category of gerbes on $S_{\et}$, functors of gerbes, and natural isomorphisms thereof.
\end{proposition}

\subsection{Well-definedness}

The construction of the associated gerbe $\gerb{E}_\epsilon(\alg{\tilde G})$ depended on the choice of Borel subgroup $\alg{B}$ and maximally split maximal torus $\alg{T}$.  Now we demonstrate that $\gerb{E}_\epsilon(\alg{\tilde G})$ is well-defined independently of these choices, in a suitable 2-categorical sense.

Consider another choice $\alg{B}_0 \supset \alg{T}_0$.  Our constructions, with these two choices of tori and Borel subgroups, yield two dual groups $\dgp{\tilde G}_0^\Vee$ and $\dgp{\tilde G}^\Vee$ and two gerbes $\gerb{E}_\epsilon(\alg{\tilde G})$ and $\gerb{E}_{0,\epsilon}(\alg{\tilde G})$.
    
Let $\sheaf{Y}_0$ and $\sheaf{Y}$ be the cocharacter lattices of $\alg{T}_0$ and $\alg{T}$, respectively, and $\Phi_0^\vee, \Phi^\vee$ the coroots therein.  The Borel subgroups provide systems of simple coroots $\Delta_0^\vee$ and $\Delta^\vee$, respectively (sheaves of sets on $S_{\et}$).  Similarly we have character lattices $\sheaf{X}_0, \sheaf{X}$, roots $\Phi_0, \Phi$, and simple roots $\Delta_0, \Delta$.  The cover $\alg{\tilde G}$ yields quadratic forms $Q \From \sheaf{Y} \To \ZZ$ and $Q_0 \From \sheaf{Y}_0 \To \ZZ$.  The second Brylinski-Deligne invariant gives extensions $\sheaf{G}_m \Into \sheaf{D} \Onto \sheaf{Y}$ and $\sheaf{G}_m \Into \sheaf{D}_0 \Onto \sheaf{Y}_0$.  

By Proposition \ref{BTConj}, there exists $g \in \alg{G}(S)$ such that $\Int(g) \alg{T}_0 = \alg{T}$ and $\Int(g) \alg{B}_0 = \alg{B}$.  This automorphism $\Int(g)$ lifts to an automorphism of $\alg{G}'$.  This defines a well-aligned isomorphism of covers.
$$\begin{tikzcd}
\alg{K}_{2} \inarrow{r} \arrow{d}{=} & \alg{G}' \onarrow{r} \arrow{d}{\Int(g)'} & \alg{G} \arrow{d}{\Int(g)} & & \alg{B}_{0} \arrow{d} & \alg{T}_{0} \arrow{d}  \\
\alg{K}_{2} \inarrow{r} & \alg{G}' \onarrow{r} & \alg{G} & & \alg{B} & \alg{T}
\end{tikzcd}$$
This well-aligned isomorphism of covers yields an equivalence of gerbes, 
$$\gerb{Int}(g) \From \gerb{E}_\epsilon(\alg{\tilde G}) \To \gerb{E}_{0,\epsilon}(\alg{\tilde G}),$$ 
lying over $\Int(g)^\vee \From \sheaf{\tilde Z^\vee} \To \sheaf{\tilde Z}_0^\Vee$.

Suppose that $g' \in \alg{G}(S)$ also satisfies $\Int(g') \alg{T}_0 = \alg{T}$ and $\Int(g') \alg{B}_0 = \alg{B}$.  As in the proof of Theorem \ref{WellDefinedDualGroup}, $g' = t g$ for a unique $t \in \alg{T}(S)$ and $\Int(g') = \Int(t) \Int(g)$.  

This gives a natural isomorphism of functors,
$$N(g', g) \From \gerb{Int}(g') \xRightarrow{\sim} \gerb{Int}(g) \circ  \gerb{Int}(t)$$
Our upcoming Proposition \ref{IntTFunctors} will provide a natural isomorphism $\gerb{Int}(t) \xRightarrow{\sim} \gerb{Id}$.  Assuming this for the moment, we find a natural isomorphism $\gerb{Int}(g') \xRightarrow{\sim} \gerb{Int}(g)$.  This demonstrates that the gerbe $\gerb{E}_\epsilon(\alg{\tilde G})$ is well-defined in the following 2-categorical sense:  
\begin{enumerate}
\item
for each pair $i = (\alg{B}, \alg{T})$ consisting of a Borel subgroup of $\alg{G}$ and a maximally split maximal torus therein, we have constructed a gerbe $\gerb{E}_\epsilon^i(\alg{\tilde G})$;
\item
for each pair $i = (\alg{B}, \alg{T})$, $j = (\alg{B}_0, \alg{T}_0)$, we have constructed a family $P(i,j)$ of gerbe equivalences $\gerb{Int}(g) \From \gerb{E}_\epsilon^i(\alg{\tilde G}) \To \gerb{E}_\epsilon^j(\alg{\tilde G})$, indexed by those $g$ which conjugate $i$ to $j$;
\item
for any two elements $g, g'$ which conjugate $i$ to $j$, there is a distinguished natural isomorphism of gerbe equivalences from $\gerb{Int}(g')$ to $\gerb{Int}(g)$.
\end{enumerate}
Once we define the natural isomorphism $\gerb{Int}(t) \xRightarrow{\sim} \gerb{Id}$, we will have defined $\gerb{E}_\epsilon(\alg{\tilde G})$ uniquely up to equivalence, the equivalences being defined uniquely up to unique natural isomorphism (we learned this notion of well-definedness from reading various works of James Milne).  Defining the natural isomorphism requires some computation and is the subject of the section below.

\subsubsection{The isomorphism $\gerb{Int}(t) \xRightarrow{\sim} \gerb{Id}$}
For what follows, define $\delta_Q \From \sheaf{Y} \To n^{-1} \sheaf{X}$ to be the unique homomorphism satisfying
\begin{equation}
\label{DefinitionOfDelta}
\langle \delta_Q(y_1), y_2 \rangle = \beta_Q(y_1, y_2) \text{ for all } y_1, y_2 \in \sheaf{Y}.
\end{equation}
The constants $m_\phi$ and $n_\phi$ arise in the following useful result.
\begin{lemma}
\label{DeltaM}
For all $\phi \in \Phi$, we have $\delta_Q(n_\phi \phi^\vee) = \delta_Q(\tilde \phi^\vee) = m_\phi \phi$.
\end{lemma}
\proof
For all $\phi \in \Phi$, and all $y \in \sheaf{Y}$, we have
\begin{equation}
\label{DeltaBeta}
\langle \delta_Q(\tilde \phi^\vee), y \rangle = \beta_Q(n_\phi \phi^\vee, y).
\end{equation}
Weyl-invariance of the quadratic form (applying the root reflection for $\phi$) implies
$$\beta_Q(\phi^\vee, y) = \beta_Q(-\phi^\vee, y - \langle \phi, y \rangle \phi^\vee) = - \beta_Q(\phi^\vee, y) + \beta_Q(\phi^\vee, \phi^\vee) \langle \phi, y \rangle.$$
Adding $\beta_Q(\phi^\vee, y)$ to both sides and dividing by two, we find that
$$\beta_Q(\phi^\vee, y) = n^{-1} Q(\phi^\vee) \langle \phi, y \rangle.$$
Substituting into \eqref{DeltaBeta} yields
$$\langle \delta_Q(\tilde \phi^\vee), y \rangle = n_\phi  \beta_Q( \phi^\vee, y)= \langle n_\phi n^{-1} Q(\phi^\vee) \cdot \phi, y \rangle.$$
Since this holds for all $y \in \sheaf{Y}$, we have
$$\delta_Q(\tilde \phi^\vee) = n_\phi n^{-1} Q(\phi^\vee) \cdot \phi= m_\phi \cdot \phi.$$
\qed

Consider the homomorphisms of sheaves of abelian groups on $S_{\et}$,
$$\sheaf{T} \xrightarrow{\delta_Q} \sheaf{\hat T} \xtwoheadrightarrow{p} \sheaf{\hat T}_{\SC},$$
obtained by applying $\shom(\bullet, \sheaf{G}_m)$ to the homomorphisms
$$\sheaf{Y}_{Q,n}^{\SC} \Into \sheaf{Y}_{Q,n} \xrightarrow{\delta_Q} \sheaf{X}.$$
An object of $\gerb{E}_\epsilon(\alg{\tilde G})$ is a triple $(\sheaf{H}, h, j)$ where
\begin{itemize}
\item
$\sheaf{H} \xrightarrow{h} \sspl(\sheaf{D}_{Q,n})$ is an $n^{\th}$ root of the $\sheaf{\hat T}$-torsor $\sheaf{D}_{Q,n}$.
\item
$j \From p_\ast \sheaf{H} \To \mu_\ast \Whit$ is an isomorphism in the gerbe $\gerb{E}_\epsilon^{\SC}(\alg{\tilde T}) = \epsilon_\ast \sqrt[n]{\sspl(\sheaf{D}_{Q,n}^{\SC})}$.  Thus $j = j_0 \wedge \tau^\vee$ where $j_0 \From p_\ast \sheaf{H} \To \mu_\ast \Whit$ is an isomorphism in the gerbe $\sqrt[n]{\sspl(\sheaf{D}_{Q,n}^{\SC})}$ and $\tau^\vee \in \sheaf{\tilde T}_{\SC}^\Vee$.
\end{itemize}

Given $\sheaf{H} \xrightarrow{h} \sspl(\sheaf{D}_{Q,n})$ and $\hat t \in \sheaf{\hat T}$, write $\hat t \circ h \From \sheaf{H} \To \sspl(\sheaf{D}_{Q,n})$ for the map obtained by composing $h$ with the automorphism $\hat t$ of the $\sheaf{\hat T}$-torsor $\sspl(\sheaf{D}_{Q,n})$.

Similarly, given $j = j_0 \wedge \tau^\vee$, with $p_\ast \sheaf{H} \xrightarrow{j_0} \mu_\ast \Whit$, $\tau^\vee \in \sheaf{\tilde T}_{\SC}^{\Vee}$, and given $\hat t_{\SC} \in \sheaf{\hat T}_{\SC}$, write $\hat t_{\SC} \circ j$ for $(\hat t_{\SC} \circ j_0) \wedge \tau^\vee$.  Here $\hat t_{\SC} \circ j_0 \From p_\ast \sheaf{H} \To \mu_\ast \Whit$ is the map obtained by composing $j$ with the automorphism $\hat t_{\SC}$ of the $\sheaf{\hat T}_{\SC}$-torsor $\mu_\ast \Whit$.

The following result describes the functor $\gerb{Int}(t)$ explicitly.
\begin{lemma}
For all $t \in \alg{T}(S)$, the equivalence of gerbes $\gerb{Int}(t) \From \gerb{E}_\epsilon(\alg{\tilde G}) \To \gerb{E}_\epsilon(\alg{\tilde G})$ sends an object $(\sheaf{H}, h,j)$ to the object $(\sheaf{H}, \delta_Q(t)^{-n} \circ h, p(\delta_Q(t))^{-1} \circ j)$.
\end{lemma}  
\proof
We work locally on $S_{\et}$ throughout the proof.  The action $\Int(t)$ of $t$ on the extension $\sheaf{G}_m \Into \sheaf{D}_{Q,n} \Onto \sheaf{Y}_{Q,n}$ is given by \cite[Equation 11.11.1]{B-D}, which states (in different notation) that
$$\Int(t) d = d \cdot \delta_Q(y)(t)^{-n} = d \cdot y \left( \delta_Q(t)^{-n} \right),$$
for all $d \in \sheaf{D}_{Q,n}$ mapping to $y \in \sheaf{Y}_{Q,n}$.  We find that $\Int(t)$ is the automorphism of the extension $\sheaf{D}_{Q,n}$ determined by the element $\delta_Q(t)^{-n} \in \sheaf{\hat T} = \shom(\sheaf{Y}_{Q,n}, \sheaf{G}_m)$.  Hence the functor $\gerb{Int}(t) \From \gerb{E}_\epsilon(\alg{\tilde T}) \To \gerb{E}_\epsilon(\alg{\tilde T})$ sends $(\sheaf{H}, h)$ to $(\sheaf{H}, \delta_Q(t)^{-n} \circ h)$.  It remains to see how $\gerb{Int}(t)$ affects the third term in a triple $(\sheaf{H}, h, j)$.  

Conjugation by $t$ gives a map of torsors $\Whit \To \Whit$.  Let $j_0'$ be the map which makes the following triangle commute:
$$\begin{tikzcd}
p_\ast \sheaf{H} \arrow{r}{j_0} \arrow{dr}[swap]{j_0'} & \mu_\ast \Whit \arrow{d}{\Int(t)} \\
\phantom{a} & \mu_\ast \Whit
\end{tikzcd}$$
We have $\gerb{Int}(t) (\sheaf{H}, h, j_0 \wedge \tau^\vee) = (\sheaf{H}, \delta_Q(t)^{-n} \circ h, j_0' \wedge \tau^\vee)$ and must describe $j_0'$.

Take an element $\eta= (\eta_\alpha : \alpha \in \Delta) \in \Whit$ (decomposed with respect to the basis of simple roots $\Delta$).  Conjugation by $t$ yields a new element $\Int(t) \eta$ satisfying
$$\Int(t) \eta(u) = \eta(t^{-1} u t), \quad \text{ for all } u \in \sheaf{U}.$$
Decomposing along the simple root subgroups,
$$\Int(t) \eta = (\alpha(t) \ast \eta_\alpha : \alpha \in \Delta) \in \Whit,$$
where $[a \ast \eta_\alpha](e) = \eta_\alpha(a^{-1} e)$ for any $a \in \sheaf{G}_m$.

Comparing to the action of $\sheaf{\hat T}_{\SC}$ on $\Whit$, we find that
$$\Int(t) \eta= \theta(t) \ast \eta,$$
where $\theta \From \sheaf{T} \To \sheaf{\hat T}_{\SC}$ denotes the homomorphism dual to the homomorphism of character lattices sending $\tilde \alpha^\vee \in \sheaf{Y}_{Q,n}^{\SC}$ to $\alpha \in \sheaf{X}$.  

In the quotient $\sheaf{\hat T}_{\SC}$-torsor $\mu_\ast \Whit$, the action is given by $\Int(t) \bar \eta = \mu(\theta(t)) \cdot \bar \eta$, for all $\bar \eta \in \mu_\ast \Whit$.  More explicitly, $\mu \circ \theta \From \sheaf{\hat T}_{\SC} \To \sheaf{T}$ is the homomorphism dual to the map of character lattices sending $\tilde \alpha^\vee$ to $-m_{\alpha} \alpha$.  By Lemma \ref{DeltaM}, we have $\mu(\theta(t)) = p(\delta_Q(t))^{-1}$ (recall that $p$ corresponds to the inclusion $\sheaf{Y}_{Q,n}^{\SC} \Into \sheaf{Y}_{Q,n}$).  It follows that $j_0' = \Int(t) \circ j_0 = p(\delta_Q(t))^{-1} \circ j_0$.  This yields the result:
$$\gerb{Int}(t) (\sheaf{H}, h, j) = (\sheaf{H}, \delta_Q(t)^{-n} \circ h, p(\delta_Q(t))^{-1} \circ j).$$
\qed

\begin{proposition}
\label{IntTFunctors}
Let $(\sheaf{H}, h, j)$ be an object of the category $\gerb{E}_\epsilon(\alg{\tilde G})$.  Then, for all $t \in \alg{T}(S)$, the morphism $\delta_Q(t)^{-1} \From \sheaf{H} \To \sheaf{H}$ defines an isomorphism from $\gerb{Int}(t) (\sheaf{H}, h, j)$ to $(\sheaf{H}, h, j)$.  As objects vary, this defines a natural isomorphism $\Fun{A}(t) \From \gerb{Int}(t) \Rightarrow \gerb{Id}$ of functors from $\gerb{E}_\epsilon(\alg{\tilde G})$ to itself.  For a pair of elements $t_1, t_2 \in \alg{T}(S)$, we have a commutative diagram of functors and natural isomorphisms.
$$\begin{tikzcd}
\gerb{Int}(t_1 t_2) \arrow[Rightarrow]{r}{=} \arrow[Rightarrow]{d}{\Fun{A}(t_1 t_2)} & \gerb{Int}(t_1) \circ \gerb{Int}(t_2) \arrow[Rightarrow]{d}{\Fun{A}(t_1) \circ \Fun{A}(t_2)} \\
\gerb{Id} \arrow[Rightarrow]{r}{=} & \gerb{Id}
\end{tikzcd}$$
\end{proposition}
\proof
As $\sheaf{H} \xrightarrow{h} \sspl(\sheaf{D}_{Q,n})$ is an $n^{\th}$ root of the torsor $\sspl(\sheaf{D}_{Q,n})$, we have $\delta_Q(t)^{-n} \circ h = h \circ \delta_Q(t)^{-1}$.  As $j = j_0 \wedge \tau$, with $j_0$ an isomorphism in the gerbe $\sqrt[n]{\sspl(\sheaf{D}_{Q,n}^{\SC})}$, we have $p(\delta_Q(t)^{-1}) \circ j = j \circ p(\delta_Q(t)^{-1})$.  

Hence the isomorphism $\delta_Q(t)^{-1} \From \sheaf{H} \To \sheaf{H}$ defines an isomorphism
\begin{align*}
\gerb{Int}(t) (\sheaf{H}, h, j) &= (\sheaf{H}, \delta_Q(t)^{-n} \circ h, p(\delta_Q(t)^{-1}) \circ j) \\
&= (\sheaf{H}, h \circ \delta_Q(t)^{-1}, j \circ p(\delta_Q(t)^{-1}) ) \\
& \xrightarrow{\delta_Q(t)^{-1}} (\sheaf{H}, h, j).
\end{align*}

As their definition depends only on $t$, these isomorphisms $\delta_Q(t)^{-1}$ form a natural isomorphism $\Fun{A}(t) \From \gerb{Int}(t) \Rightarrow \gerb{Id}$.  As $\delta_Q(t_1 t_2)^{-1} = \delta_Q(t_1)^{-1} \delta_Q(t_2)^{-1}$, we find an {\em equality} of functors $\gerb{Int}(t_1 t_2) = \gerb{Int}(t_1) \circ \gerb{Int}(t_2)$, and the commutative diagram of the proposition.  
\qed

\subsection{Change of base scheme}
\label{BaseChangeGerbe}

Let $\alg{\tilde G}$ be a degree $n$ cover of a quasisplit group $\alg{G} \supset \alg{B} \supset \alg{T}$ over $S$, as before.  Let $\gamma \From S_0 \To S$ be a morphism of schemes as in Section \ref{BaseChangeDualGroup}.  Then $\gamma$ gives rise to a pullback functor $\gamma^\ast$ from sheaves on $S_{\et}$ to sheaves on $S_{0,\et}$, and from gerbes on $S_{\et}$ to gerbes on $S_{0, \et}$.  

Assuming that $\amu_n(S)$ and $\amu_n(S_0)$ are cyclic groups of order $n$, we may identify these groups via $\gamma \From \amu_n(S) \To \amu_n(S_0)$, and a character $\epsilon \From \mu_n \Into \CC^\times$ corresponds to a character $\epsilon_0 \From \amu_n(S_0) \Into \CC^\times$.

Pullback via $\gamma$ defines a degree $n$ cover $\alg{\tilde G}_0$ of a quasisplit group $\alg{G}_0 \supset \alg{B}_0 \supset \alg{T}_0$ over $S_0$.  We have constructed gerbes $\gerb{E}_\epsilon(\alg{\tilde G})$ and $\gerb{E}_{\epsilon_0}(\alg{\tilde G}_0)$ associated to this data, banded by $\sheaf{\tilde Z}^\Vee$ and $\sheaf{\tilde Z}_0^\Vee$, respectively.  We also consider the pullback gerbe $\gamma^\ast \gerb{E}_\epsilon(\alg{\tilde G})$, banded by $\gamma^\ast \sheaf{\tilde Z}^\Vee$, and recall from Section \ref{BaseChangeDualGroup} that there is a natural isomorphism $N^\vee \From \gamma^\ast \sheaf{\tilde Z}^\Vee \To \sheaf{\tilde Z}_0^\Vee$.

Similarly, we find isomorphisms of sheaves on $S_{0,\et}$,
\begin{enumerate}
\item $\hat N \From \gamma^\ast \sheaf{\hat T} \xrightarrow{\sim} \sheaf{\hat T}_0$ (from $N \From \gamma^\ast \sheaf{Y}_{Q,n} \xrightarrow{\sim} \sheaf{Y}_{0,Q_0,n}$);
\item $N^\vee \From \gamma^\ast \sheaf{\tilde T}_{\SC}^\Vee \xrightarrow{\sim} \sheaf{\tilde T}_{0,\SC}^\Vee$;
\item $N \From \gamma^\ast \sheaf{D}_{Q,n} \xrightarrow{\sim} \sheaf{D}_{0,Q_0,n}$ (the construction of the second Brylinski-Deligne invariant is compatible with pullbacks); 
\item $N \From \gamma^\ast \Whit \xrightarrow{\sim} \Whit_0$ (since $\alg{B}_0 = \gamma^\ast \alg{B}$).
\end{enumerate}

From these observed isomorphisms, we find that the pullback of \'etale sheaves from $S_{\et}$ to $S_{0,\et}$ defines a functor 
$$\gerb{N}' \From \gerb{E}_\epsilon(\alg{\tilde G})  \To \gamma_\ast \gerb{E}_{\epsilon_0}(\alg{\tilde G}_0)),$$
given on objects by $\gerb{N}'(\sheaf{H}, h, j) = (\gamma^\ast \sheaf{H}, \gamma^\ast h, \gamma^\ast j)$.  For example, if $\sheaf{H}$ is a $\sheaf{\hat T}$-torsor, then $\gamma^\ast \sheaf{H}$, a priori a $\gamma^\ast \sheaf{\hat T}$-torsor, becomes a $\sheaf{\hat T}_0$-torsor via $\hat N$.  

From \cite[Chapitre V, Proposition 3.1.8]{Giraud}, such an equivalence of gerbes $\gerb{N}'$ determines a unique, up to unique natural isomorphism, equivalence of gerbes,
$$\gerb{N} \From \gamma^\ast \gerb{E}_\epsilon(\alg{\tilde G}) \To \gerb{E}_{\epsilon_0}(\alg{\tilde G}_0),$$
lying over the natural isomorphism of bands $N^\vee \From \gamma^\ast \sheaf{\tilde Z}^\Vee \To \sheaf{\tilde Z}_0^\Vee$.  In this way, the construction of the gerbe associated to a cover is compatible with change of base scheme.

 \subsection{Parabolic subgroups}
\label{ParabolicGerbe}

Return to the degree $n$ cover $\alg{\tilde G}$ of a quasisplit group $\alg{G} \supset \alg{B} \supset \alg{T}$ over $S$.  Let $\alg{P} \subset \alg{G}$ be a parabolic subgroup defined over $S$, containing $\alg{B}$.  As in Section \ref{ParabolicDualGroup}, we consider a Levi decomposition $\alg{P} = \alg{M} \alg{N}$ and resulting cover $\alg{\tilde M}$.  Fix $\epsilon \From \mu_n \Into \CC^\times$ as before.

Consider the gerbes $\gerb{E}_\epsilon(\alg{\tilde G})$ and $\gerb{E}_\epsilon(\alg{\tilde M})$.  An object of $\gerb{E}_\epsilon(\alg{\tilde G})$ is a triple $(\sheaf{H}, h, j)$, where $(\sheaf{H}, h)$ is an $n^{\th}$ root of the torsor $\sspl(\sheaf{D}_{Q,n})$, $j = j_0 \wedge \tau$, $\tau \in \sheaf{\tilde T}_{\SC}^\Vee$, and $j_0 \From p_\ast\sheaf{H} \To \mu_\ast \Whit$ is an isomorphism of $\sheaf{\hat T}_{\SC}$-torsors.

Restriction of characters via $\alg{U} \supset \alg{U}_{\alg{M}} = \alg{U} \cap \alg{M}$ provides a homomorphism of sheaves from $\Whit$ to $\Whit_{\alg{M}}$ (the Whittaker torsor for $\alg{M}$).  The inclusion $\sheaf{Y}_{\alg{M}}^{\SC} \subset \sheaf{Y}^{\SC}$ of coroot lattices provides a homomorphism of sheaves of abelian groups $\sheaf{\hat T}_{\SC} \Onto \sheaf{\hat T}_{\alg{M},\SC}$, where 
$$\sheaf{\hat T}_{\SC} = \shom(\sheaf{Y}^{\SC}, \sheaf{G}_m), \quad \sheaf{\hat T}_{\alg{M},\SC} = \shom(\sheaf{Y}_{\alg{M}}^{\SC}, \sheaf{G}_m).$$
Similarly, it provides a homomorphism $\sheaf{\tilde T}_{\SC}^\Vee \Onto \sheaf{\tilde T}_{\alg{M}, \SC}^\Vee$, where 
$$\sheaf{\tilde T}_{\SC}^\Vee = \shom(\sheaf{Y}^{\SC}, \CC^\times), \quad \sheaf{\tilde T}_{\alg{M},\SC}^\Vee = \shom(\sheaf{Y}_{\alg{M}}^{\SC}, \CC^\times).$$
Define $\tau_{\alg{M}}^\vee$ to be the image of $\tau^\vee$ under this homomorphism.

The constants defining $\mu \From \sheaf{\hat T}_{\SC} \To \sheaf{\hat T}_{\SC}$ are the same as those defining the corresponding map, $\mu_{\alg{M}} \From \sheaf{\hat T}_{\alg{M},\SC} \To \sheaf{\hat T}_{\alg{M},\SC}$.  We find a commutative diagram.
$$\begin{tikzcd}
 \sheaf{\hat T}_{\SC} \times \mu_\ast \Whit  \arrow{r}{\ast} \arrow{d} & \mu_\ast \Whit \arrow{d} \\
  \sheaf{\hat T}_{\alg{M},\SC} \times (\mu_{\alg{M}})_\ast  \Whit_{\alg{M}} \arrow{r}{\ast} &  (\mu_{\alg{M}})_\ast \Whit_{\alg{M}} 
\end{tikzcd}$$

Composing $j_0 \From p_\ast \sheaf{H} \To \mu_\ast \Whit$ with the map to $(\mu_{\alg{M}})_\ast \Whit_{\alg{M}}$ defines a map
$$j_{\alg{M},0} \From \sheaf{H} \To (\mu_{\alg{M}})_\ast \Whit_{\alg{M}}.$$
This defines a functor of gerbes,
$$\gerb{i} \From \gerb{E}_\epsilon(\alg{\tilde G}) \To \gerb{E}_\epsilon(\alg{\tilde M}), \quad (\sheaf{H}, h, j_0 \wedge \tau^\vee) \mapsto (\sheaf{H}, h, j_{\alg{M},0} \wedge \tau_{\alg{M}}^\vee),$$
lying over $\iota^\vee \From \sheaf{\tilde Z}^\Vee \Into \sheaf{\tilde Z}_{\alg{M}}^\Vee$.

\section{The metaGalois group}
\label{MetaGaloisSection}
We now specialize to three classes of base scheme of arithmetic interest.
\begin{description}
\item[Global]  $S = \Spec(F)$ for a global field $F$;
\item[Local]  $S = \Spec(F)$ for a local field $F$;
\item[Local integral]  $S = \Spec(\OO)$ for the ring of integers $\OO$ in a nonarchimedean local field $F$.
\end{description}

Choose a separable closure $\bar F / F$ in all three cases, and write $\Gal_F = \Gal(\bar F / F)$ for the resulting absolute Galois group.  In the local integral case, the separable closure $\bar F / F$ provides a geometric base point $\bar s$ for $\Spec(\OO)$ as well, and define $\Gal_\OO = \pi_1^{\et}(\Spec(\OO), \bar s)$.  This is a profinite group, topologically generated by a geometric Frobenius element $\Fr$; thus we write $\Gal_\OO = \langle \Fr \rangle_{\prof}$.

When $S = \Spec(F)$ for a local or global field, or $S = \Spec(\OO)$ for the ring of integers in a nonarchimedean local field, write $\Gal_S$ for $\Gal_F$ or $\Gal_\OO$ accordingly.

\subsection{Construction of the metaGalois group}
The metaGalois group will be a profinite group fitting into a central extension,
$$\mu_2 \Into \mGal_S \Onto \Gal_S.$$
When $F$ has characteristic two, we {\em define} the metaGalois group $\mGal_F$ to be the trivial extension $\Gal_F \times \mu_2$.  The metaGalois group $\mGal_\OO$ will {\em not be defined} when $\OO$ has residual characteristic $2$, reflecting the idea that metaGalois representations cannot be ``unramified at $2$'' (though one might propose an alternative notion of ``minimally ramified'').

\subsubsection{Local fields}

When $F$ is a local field (with $2 \neq 0$) the quadratic Hilbert symbol defines a symmetric nondegenerate $\ZZ$-bilinear form
$$\Hilb_2 \From F^\times_{/2} \times F^\times_{/2} \To \mu_2.$$
The abelianized Galois group $\Gal_F^{\ab}$ is defined, up to unique isomorphism, from $F$ alone.  When $F$ is nonarchimedean, we normalize the valuation so that $\val(F^\times) = \ZZ$, and we normalize the reciprocity map of local class field theory to send a geometric Frobenius element to an element of valuation $1$.  The reciprocity map gives a surjective homomorphism
$$\rec_{F/2} \From \Gal_F^{\ab} \Onto F^\times_{/2}.$$
Composing the Hilbert symbol with $\rec_{F/2}$ defines a function
$$h \From \Gal_F^{\ab} \times \Gal_F^{\ab} \To \mu_2,$$
and it is straightforward to verify that $h$ is a (bimultiplicative) continuous symmetric $2$-cocycle.  This incarnates a commutative extension
$$\mu_2 \Into \mGal_F^{\ab} \Onto \Gal_F^{\ab}$$
of profinite groups.  Concretely, $\mGal_F^{\ab} = \Gal_F^{\ab} \times \mu_2$ as sets, and 
$$(\gamma_1, \epsilon_1) \cdot (\gamma_2, \epsilon_2) \defeq \left( \gamma_1 \gamma_2, \epsilon_1 \epsilon_2 \cdot h(\gamma_1, \gamma_2) \right).$$

The pullback of this extension to $\Gal_F$ will be called the \defined{metaGalois group of $F$}, written $\mGal_F$.  It is a central extension of $\Gal_F$ by $\mu_2$,
$$\mu_2 \Into \mGal_F \Onto \Gal_F.$$

\subsubsection{The local integral case}

Suppose that $F$ is a nonarchimedean \textbf{nondyadic} (i.e., $\val(2) = 0$) local field.  Then the quadratic Hilbert symbol satisfies
$$\Hilb_2(u,v) = 1 \text{ for all } u,v \in \OO^\times.$$
Write $\Inertia \subset \Gal_F$ for the inertial subgroup, so that $\rec_F(\Inertia) = \OO^\times$.  The cocycle $h$ is trivial when restricted to $\Inertia \times \Inertia$.  Thus $\gamma \mapsto (\gamma, 1)$ gives a canonical splitting $\sigma^\circ \From \Inertia \Into \mGal_F$.  The natural map $\Gal_F \Onto \Gal_\OO$ identifies $\Gal_\OO$ with $\Gal_F / \Inertia$.  Define $\mGal_\OO = \mGal_F / \sigma^\circ(\Inertia)$ to obtain a commutative diagram with exact rows.
$$\begin{tikzcd}
\mu_2 \arrow{d}{=} \inarrow{r} & \mGal_F \arrow{r} \onarrow{d} & \Gal_F \onarrow{d} \\
\mu_2 \inarrow{r} & \mGal_\OO \arrow{r} & \Gal_\OO
\end{tikzcd}$$
We call $\mGal_\OO$ the \defined{metaGalois group of $\OO$}.  If $\gamma \in \Gal_F$ lifts $\Frob$, then every element of $\mGal_\OO$ is equal to $(\gamma^{\hat n}, \pm 1)$ (mod $\sigma^\circ(\Inertia)$) for some $\hat n \in \hat \ZZ$.  In this way, a Frobenius lift provides an isomorphism from $\mGal_\OO$ to the group with underlying set $\langle \gamma \rangle_{\prof} \times \mu_2$ and multiplication given by
$$(\gamma^{\hat n_1}, \epsilon_1) \cdot (\gamma^{\hat n_2}, \epsilon_2) = \left( \gamma^{\hat n_1 + \hat n_2}, \epsilon_1 \epsilon_2 \cdot (-1)^{ \hat n_1 \hat n_2 (q-1)/2} \right),$$
where $q$ is the cardinality of the residue field of $\OO$.

\subsubsection{Global fields}

When $F$ is a global field (with $2 \neq 0$ as before), the Hilbert symbol defines a symmetric $\ZZ$-bilinear form,
$$\Hilb_2 \From \AA^\times_{/2} \times \AA^\times_{/2} \To \mu_2,$$
obtained as the product of local Hilbert symbols.  This defines a continuous symmetric 2-cocycle, from which we get a commutative extension,
\begin{equation}
\label{ExtensionAdelic}
\mu_2 \Into \widetilde{\AA^\times_{/2}} \Onto \AA^\times_{/2}.
\end{equation}
Global quadratic reciprocity for the Hilbert symbol ($\Hilb_2(u,v) = 1$ for all $u,v \in F^\times$) provides a canonical splitting $\sigma_F \From F^\times_{/2} \Into \widetilde{\AA^\times_{/2}}$.  Taking the quotient yields a commutative extension,
\begin{equation}
\label{QuotientExtension}
\mu_2 \Into \frac{\widetilde{\AA^\times_{/2}}}{\sigma_F(F^\times_{/2})} \Onto \frac{ \AA^\times_{/2} }{F^\times_{/2}}.
\end{equation}
The global reciprocity map of class field theory gives an surjective homomorphism,
$$\rec_{F/2} \From \Gal_F^{\ab} \Onto \AA^\times_{/2} / F^\times_{/2} \ident (\AA^\times / F^\times)_{/2}.$$
Pulling back \eqref{QuotientExtension} via $\rec_{F/2}$ yields a commutative extension,
$$\mu_2 \Into \mGal_F^{\ab} \Onto \Gal_F^{\ab}.$$
Pulling back via $\Gal_F \Onto \Gal_F^{\ab}$ defines the \defined{metaGalois group of $F$},
$$\mu_2 \Into \mGal_F \Onto \Gal_F.$$

\subsubsection{Compatibilities}
\label{CompatibilityMetaGalois}
If $v \in \VV$ is a place of a global field $F$, then an embedding $\bar F \Into \bar F_v$ of separable closures determines an injective homomorphism $\iota_v \From \Gal_{F_v} \Into \Gal_F$.  As the global Hilbert symbol is the product of local ones, we find a homomorphism $\tilde \iota_v$ realizing $\widetilde{\Gal}_{F_v}$ as the pullback of the extension $\widetilde{\Gal}_F$.
$$\begin{tikzcd}
\mu_2 \inarrow{r} \arrow{d}{=} & \mGal_{F_v} \arrow{r} \inarrow{d}{\tilde \iota_v} & \Gal_{F_v} \inarrow{d}{\iota_v} \\
\mu_2 \inarrow{r} & \mGal_F \arrow{r} & \Gal_F
\end{tikzcd}$$

For local and global fields $F$, a choice of separable closure $\bar F / F$ entered the construction of the metaGalois group.  Suppose that $\bar F_0$ is another separable closure of $F$.  Every $F$-algebra isomorphism $\iota \From \bar F_0 \xrightarrow{\sim} \bar F$ yields an isomorphism $\iota \From \Gal(\bar F_0 / F) \xrightarrow{\sim} \Gal(\bar F / F)$.  The resulting isomorphism $\Gal(\bar F_0 / F)^{\ab} \xrightarrow{\sim} \Gal(\bar F / F)^{\ab}$ does not depend on $\iota$.  

The separable closure $\bar F_0$ yields a cocycle $h_0 \From \Gal(\bar F_0 / F)^{\ab} \times \Gal(\bar F_0 / F)^{\ab} \To \mu_2$, and thus a metaGalois group $\mGal(\bar F_0 / F)$.  Since the defining cocycles $h_0$ and $h$ factor through abelianized Galois groups, the isomorphism $\iota \From \Gal(\bar F_0 / F) \To \Gal(\bar F / F)$ lifts canonically to an isomorphism of metaGalois groups.
$$\begin{tikzcd}
\mu_2 \arrow{d}{=} \inarrow{r} & \mGal(\bar F_0 / F)\onarrow{r} \arrow{d}{\tilde \iota} & \Gal(\bar F_0 / F)\arrow{d}{\iota} \\
\mu_2 \inarrow{r} & \mGal(\bar F / F) \onarrow{r} & \Gal(\bar F / F)
\end{tikzcd}$$

\subsection{The Brauer class}
When $F$ is a local or global field, the metaGalois group is an extension $\mu_2 \Into \mGal_F \Onto \Gal_F$.  As such, it has a cohomology class in the Brauer group.
$$\left[ \mGal_F \right] \in H_{\et}^2(F, \mu_2) = \Br(F)_{[2]}.$$

This Brauer class is often trivial -- the metaGalois group often splits, though it rarely has a canonical splitting.
\begin{proposition}
\label{OddRCTrivialBr}
Suppose that $F$ is a \textbf{nondyadic} (i.e., $\val(2) = 0$) nonarchimedean local field.  Then $\left[ \mGal_F \right]$ is the trivial class.
\end{proposition}
\proof
The projection $\mGal_F \To \mGal_\OO$ identifies the metaGalois group of $F$ with the pullback of the metaGalois group of $\OO$.
$$\begin{tikzcd}
\mu_2 \arrow{r} \arrow{d}{=} & \mGal_F \arrow{r} \arrow{d} & \Gal_F \arrow{d} \\
\mu_2 \arrow{r} & \mGal_\OO \arrow{r} & \Gal_\OO = \langle \Fr \rangle_{\prof}
\end{tikzcd}$$
But every extension of $\hat \ZZ$ by $\mu_2$ splits (though not canonically); hence the metaGalois group splits and its Brauer class is trivial.
\qed

\begin{proposition}
\label{RealNontrivBr}
Over $\RR$, the metaGalois group is a nonsplit extension, so $\left[ \mGal_\RR \right]$ is the unique nontrivial class in the Brauer group $\Br(\RR)$.
\end{proposition}
\proof
Let $\sigma$ denote complex conjugation, $\Gal_\RR = \Gal(\CC / \RR) = \{ \Id, \sigma \}$.  The metaGalois group is a cyclic group of order $4$ sitting in an extension
$$\mu_2 \Into \mGal_\RR \Onto \Gal_\RR.$$
Indeed, the cocycle $h$ satisfies $h(\sigma, \sigma) = \Hilb_2(-1,-1) = -1$.  Thus
$$(\sigma,1) \cdot (\sigma,1) = (\Id, -1) \in \mGal_\RR.$$
Hence $(\sigma,1)$ is an element of order $4$ and $\left[ \mGal_\RR \right]$ is nontrivial.
\qed

\begin{proposition}
Let $F_2$ be a \textbf{dyadic} nonarchimedean local field of characteristic zero.  Then $\left[ \mGal_{F_2} \right]$ is trivial if $[F_2 : \QQ_2]$ is even, and is nontrivial if $[F_2 : \QQ_2]$ is odd.  
\end{proposition}
\proof
Let $d = [F_2 : \QQ_2]$.  By approximation, there exists a global field $F$ such that $F \otimes_\QQ \QQ_2$ is isomorphic to $F_2$ as an $F$-algebra.  The global metaGalois group $\mGal_F$ has a Brauer class $\beta_F$ with local components $\beta_{F,v}$ satisfying
\begin{itemize}
\item
$\beta_{F,v}$ is trivial when $F_v$ has odd residual characteristic (by Proposition \ref{OddRCTrivialBr});
\item
$\beta_{F,2} = \left[ \mGal_{F_2} \right]$ at the unique place of even residual characteristic;
\item
$\beta_{F,v}$ is nontrivial at all real places (by Proposition \ref{RealNontrivBr});
\item
(Parity condition) $\beta_{F,v}$ is nontrivial at a set of places of even cardinality.
\end{itemize}
We have $d = [F_2 : \QQ_2] = [F : \QQ] = r_1 + 2 r_2$, where $r_1$ is the number of real places, and $r_2$ the number of complex places.  It follows that $d$ is even if and only if $r_1$ is even.  The parity condition on the global Brauer class implies that $r_1$ is even if and only if $\beta_{F,2}$ is the trivial class.
\qed

\begin{corollary}
Let $F$ be a global field, with $2 \neq 0$ in $F$.  Then the Brauer class of $\mGal_F$ is that of the unique quaternion algebra which is ramified at all real places and all dyadic places of odd degree over $\QQ_2$.  
\end{corollary}
\proof
This follows directly from the previous three propositions, and the local-global compatibility of the metaGalois group.
\qed

In particular, the Brauer class $\left[ \mGal_\QQ \right]$ is that of the quaternion algebra $\frac{(-1,-1)}{\QQ}$ ramified only at $2$ and $\infty$.  If $F$ is a global field of characteristic $p \neq 2$, then $\left[ \mGal_F \right]$ is the trivial class.

\subsection{Splitting by additive characters}

The metaGalois group may be a nonsplit extension of $\Gal_S$ by $\mu_2$, and even when it splits, it rarely splits canonically.  However, an additive character suffices to split the metaGalois group after pushing out via $\mu_2 \Into\mu_4$.  In the three cases of interest, define a $\sheaf{G}_m[S]$-torsor $\Psi_S$ as follows.
\begin{itemize}
\item
When $F$ is local, let $\Psi_F$ be the set of nontrivial continuous homomorphisms from $F$ to $\CC^\times$.  If $u \in F^\times$, $\psi \in \Psi_F$, write $[u \ast \psi](x) = \psi(u^{-1} x)$.  In this way, $\Psi_F$ is a $F^\times$-torsor.
\item
When $F$ is global, let $\Psi_F$  be the set of nontrivial continuous homomorphisms from $\AA / F$ to $\CC^\times$.  If $u \in F^\times$, $\psi \in \Psi_F$, write $[u \ast \psi](x) = \psi(u^{-1} x)$.  In this way, $\Psi_F$ is a $F^\times$-torsor.
\item
When $F$ is local nonarchimedean, with ring of integers $\OO$, let $\Psi_\OO$  be the set of nontrivial continuous homomorphisms from $F / \OO$ to $\CC^\times$.  If $u \in \OO^\times$, $\psi \in \Psi_\OO$, write $[u \ast \psi](x) = \psi(u^{-1} x)$.  In this way, $\Psi_\OO$ is a $\OO^\times$-torsor.
\end{itemize}

Define here $\mu_4 = \amu_4(\CC) = \{ 1,-1, i, -i \}$.  When $F$ is a local field (with $2 \neq 0$), and $\psi \in \Psi_F$, the \defined{Weil index} is a function $\weil_F(\bullet, \psi) \From F^\times_{/2} \To \mu_4$ which satisfies
\begin{equation}
\label{WeilHilbert}
\frac{\weil_F(uv, \psi)}{\weil_F(u, \psi) \weil_F(v,\psi)} = \Hilb_2(u,v).
\end{equation}
Our $\weil_F(u, \psi)$ is defined to be $\gamma(u x^2) / \gamma(x^2)$ in Weil's notation from \cite[\S 29]{WeilActa} and is written $\gamma_F(u, \psi)$ in \cite[\S A.3]{RangaRao} and elsewhere.  From \cite[Proposition A.11]{RangaRao}, the local Weil indices are trivial on $\OO^\times$ at all nondyadic places. 

When $F$ is a global field and $\psi \in \Psi_F$, the \defined{Weil index} is the function
$$\weil_F(\bullet, \psi) \From \AA^\times_{/2} \To \mu_4$$
defined as the product of local Weil indices.  The global Weil index is trivial on $F^\times$ by \cite[\S II.30, Proposition 5]{WeilActa}.  As the global Hilbert symbol $\Hilb_2 \From \AA^\times \times \AA^\times \To \mu_2$ is defined as the product of local Hilbert symbols, the formula \eqref{WeilHilbert} holds in the global setting too.  

Write $\mGal_S^{(4)}$ for the pushout of $\mGal_S$ via the inclusion $\mu_2 \Into \mu_4$ (when $S = \Spec(F)$ or $S = \Spec(\OO)$ as usual).
$$\begin{tikzcd}[column sep = 4em]
\mu_2 \inarrow{r} \inarrow{d}{\iota} & \mGal_S \onarrow{r} \inarrow{d} & \Gal_S \arrow{d}{=} \\
\mu_4 \inarrow{r} & \mGal_S^{(4)} \onarrow{r} & \Gal_S \arrow[bend right=20, dashed]{l}[swap]{s(\psi)?}
\end{tikzcd}$$
The splittings of $\mGal_S^{(4)}$, if they exist, form a $\Hom(\Gal_S, \mu_4)$-torsor.  In what follows, if $u \in \sheaf{G}_m[S]$, define $\chi_u \From \Gal_S \To \mu_2$ to be the quadratic character associated to the \'etale extension $F[\sqrt{u}]$ (in the local or global case) or $\OO[\sqrt{u}]$ (in the nondyadic local integral case).

\begin{proposition}
\label{WeilSplitsMetaGalois}
For each additive character $\psi \in \Psi_S$, the Weil index provides a splitting $s(\psi) \From \Gal_S \To \mGal_S^{(4)}$.  Moreover, this system of splittings satisfies
$$s(u \ast \psi) = \chi_u \ast s(\psi) \text{ for all } u \in \sheaf{G}_m[S].$$
\end{proposition}
The splittings $s(\psi)$ are described in three cases below.

\subsubsection{Local fields}
When $F$ is a local field, the pushout $\mGal_S^{(4)}$ can be identified with the product $\Gal_F \times \mu_4$ as a set, with multiplication given by
$$(\gamma_1, \zeta_1) \cdot (\gamma_2, \zeta_2) = \left( \gamma_1 \gamma_2, \zeta_1 \zeta_2 \cdot \Hilb_2(\rec_{F/2}(\gamma_1), \rec_{F/2}(\gamma_2) ) \right).$$
For $\psi \in \Psi_F$, \eqref{WeilHilbert} provides a splitting $s(\psi) \From \Gal_F \To \mGal_F^{(4)}$,
$$s(\psi)(\gamma) = \left( \gamma, \weil_F( \rec_{F/2}(\gamma), \psi) \right), \text{ for all } \gamma \in \Gal_F.$$

If $u \in F^\times$, then \cite[Corollary A.5]{RangaRao} states that $\weil_F(a, u \ast \psi) = \Hilb_2(a,u) \cdot \weil_F(a, \psi)$.  Since $\Hilb_2(\rec_{F/2} \gamma, u) = \chi_u(\gamma)$, we find $s(u \ast \psi) = \chi_u \ast s(\psi)$.

\subsubsection{The local integral case}

When $F$ is a nonarchimedean, nondyadic local field, the local Weil index is trivial on $\OO^\times$.  Given a character $\psi \in \Psi_\OO \subset \Psi_F$, the splitting $s(\psi) \From \Gal_F \To \mGal_F^{(4)}$ coincides with the canonical splitting $\sigma^\circ$ on intertia,
$$s(\psi) (\gamma) = \sigma^\circ(\gamma) = (\gamma, 1), \text{ for all } \gamma \in \Inertia.$$
It follows that $s(\psi)$ descends to a splitting of $\mGal_\OO^{(4)}$ at nondyadic places.
$$\begin{tikzcd}[column sep = 4em]
\mu_4 \inarrow{r} & \mGal_\OO^{(4)} \onarrow{r} & \Gal_\OO. \arrow[bend right=20]{l}[swap]{s(\psi)}
\end{tikzcd}$$

If $u \in \OO^\times$, write $\bar u$ for its image in the residue field $\FF_q$.  As before, we have $s(u \ast \psi) = \chi_u \cdot s(\psi)$.  But now, the quadratic character $\chi_u$ is restricted to $\Gal_\OO = \langle \Frob \rangle_{\prof}$; we have
$$\chi_u(\Frob) = \bar u^{(q-1) / 2} \in \mu_2.$$
In other words, $\chi_u$ is the character of $\Gal_\OO$ which sends $\Frob$ to the Legendre symbol of the reduction of $u$.  

\subsubsection{Global fields}

In the global setting, pushing out via $\mu_2 \Into \mu_4$ gives a short exact sequence
\begin{equation}
\label{Adelic4}
\mu_4 \Into \widetilde{\AA^\times_{/2}}^{(4)} \Onto \AA^\times_{/2}.
\end{equation}
The middle term is given by $\widetilde{\AA^\times_{/2}}^{(4)} = \AA^\times_{/2} \times \mu_4$ as a set, with multiplication given by
$$(u_1, \zeta_1) \cdot (u_2, \zeta_2) = (u_1 u_2, \zeta_1 \zeta_2 \cdot \Hilb_2(u_1, u_2)).$$
A character $\psi \in \Psi_F$ provides a splitting of the extension \eqref{Adelic4},
$$s_\AA(\psi)(u) = (u, \weil_F(u, \psi) ) \text{ for all } u \in \AA^\times_{/2}.$$
Since $\weil_F(u, \psi) = 1$ for all $u \in F^\times$, this splitting restricts to the canonical splitting $\sigma_F \From F_{/2}^\times \To \widetilde{\AA_{/2}^\times}$.  Thus $s_\AA(\psi)$ descends and pulls back to a splitting $s(\psi) \From \Gal_F \To  \mGal_F^{(4)}$.  If $u \in F^\times$, then our local results and local-global compatibility imply that $s(u \ast \psi) = \chi_u \cdot s(\psi)$.

\subsection{Restriction}

Suppose that $F' / F$ is a finite separable extension with $F' \subset \bar F$.  In the local integral case, suppose that $F' / F$ is unramified and let $\OO'$ be the ring of integers in $F'$.  Write $S' = \Spec(F')$ in the cases of local or global fields, and write $S' = \Spec(\OO')$ in the local integral case.  We have defined metaGalois groups for $S$ and $S'$.
\begin{equation}
\label{MGalInc}
\begin{tikzcd}
\mu_2 \inarrow{r} \arrow{d}{=} & \mGal_{S'} \onarrow{r} \arrow[dashed]{d}{?}& \Gal_{S'} \inarrow{d} \\
\mu_2 \inarrow{r} & \mGal_S \onarrow{r} & \Gal_S
\end{tikzcd}
\end{equation}

The inclusion $F' \subset \bar F$ gives an inclusion of Galois groups $\Gal_{S'} \Into \Gal_S$, but a natural inclusion of metaGalois groups is not obvious.  In particular, the cocycle defining $\mGal_S$ does not restrict to the cocycle defining $\mGal_{S'}$.

Fortunately, a beautiful insight of Wee Teck Gan gives such an inclusion of metaGalois groups, using a ``lifting theorem'' of Edward Bender \cite{Bender}.  We explain this insight here.

In the case of local fields, consider a nonzero element $u \in F'$, and the ``trace form'' (cf. \cite{SerreTrace}) $F' \To F$ given by $x \mapsto \Tr_{F' / F}(u x^2)$.  Viewing this as a quadratic form on a finite-dimensional $F$-vector space $F'$, it has a Hasse-Witt invariant (an element of $\{ \pm 1 \}$).  Define
$$\HW(u) = \frac{ \text{Hasse-Witt invariant of } x \mapsto \Tr_{F'/F}(u x^2) }{ \text{Hasse-Witt invariant of } x \mapsto \Tr_{F'/F}(x^2) }.$$
This function depends only on the square class of $u$.

Bender's theorem \cite[Theorem 1]{Bender} states that
$$\Hilb_{F',2}(u,v) = \frac{ \HW(u) \HW(v) }{ \HW(uv) } \cdot \Hilb_{F,2} ( \Norm_{F'/F} u, \Norm_{F'/F} v).$$
Let $\iota \From \Gal_{F'} \Into \Gal_F$ be the canonical inclusion, so that $\rec_{F}( \iota(\gamma)) = \Norm_{F'/F} \rec_F(\gamma)$ for all $\gamma \in \Gal_{F'}$.
\begin{proposition}
Let $F$ be a local field (with $2 \neq 0$ as usual).  Then the function $\tilde \iota \From \mGal_{F'} \Into \mGal_F$, given by
$$\tilde \iota (\gamma, \pm 1) = \left( \iota(\gamma), \pm \HW( \rec_{F'} \gamma) \right)$$
is a group homomorphism completing the commutative diagram \eqref{MGalInc}.
\end{proposition}
\proof
Consider any $\gamma_1, \gamma_2 \in \Gal_{F'}$ and define $u_1 \defeq \rec_{F'}(\gamma_1)$, $u_2 \defeq \rec_{F'}(\gamma_2)$.  Thus  $\rec_{F}( \iota(\gamma_1)) = \Norm_{F'/F} u_1$ and $\rec_{F}( \iota(\gamma_2)) = \Norm_{F'/F} u_2$.  For all $\epsilon_1, \epsilon_2 \in \{ \pm 1 \}$, we compute
\begin{align*}
\tilde \iota \left( (\gamma_1, \epsilon_1) \cdot (\gamma_2, \epsilon_2) \right) &= \tilde \iota \left( \gamma_1 \gamma_2, \epsilon_1 \epsilon_2 \Hilb_{F',2} \left( \rec_{F'}(\gamma_1), \rec_{F'}(\gamma_2) \right) \right)  \\
& = \left( \iota(\gamma_1 \gamma_2), \epsilon_1 \epsilon_2 \Hilb_{F',2} (u_1, u_2) \cdot \HW(u_1 u_2) \right) \\
& = \left( \iota(\gamma_1) \iota(\gamma_2), \epsilon_1 \epsilon_2 \Hilb_{F,2} (\Norm_{F' / F} u_1, \Norm_{F' / F} u_2) \HW(u_1) \HW(u_2) \right) \\
& = \left( \iota(\gamma_1), \epsilon_1 \HW(u_1) \right) \cdot \left( \iota(\gamma_2), \epsilon_2 \HW(u_2) \right) \\
&= \tilde \iota( \gamma_1, \epsilon_1) \cdot \tilde \iota(\gamma_2, \epsilon_2).
\end{align*}
\qed

In the local integral case, when $\OO$ is the ring of integers in a nondyadic nonarchimedean field,  $\HW(u) = 1$ for all $u \in \OO^\times$.  From this it follows that $\tilde \iota \From \mGal_{F'} \Into \mGal_{F}$ descends to an injective homomorphism.
$$\begin{tikzcd}
\mu_2 \inarrow{r} \arrow{d}{=} & \mGal_{\OO'} \onarrow{r} \arrow{d}{\tilde \iota}& \Gal_{\OO'} \inarrow{d}{\iota} \\
\mu_2 \inarrow{r} & \mGal_\OO \onarrow{r} & \Gal_\OO
\end{tikzcd}$$

In the global case, when $F$ is a number field, we note that $\prod_{v} \HW_v(u) = 1$ for all $u \in F^\times$ (here $\HW_v$ denotes the invariant as the place $v$).  From this it follows that the injective homomorphisms $\tilde \iota_v \From \mGal_{F_v'} \Into \mGal_{F_v}$ yield a injective homomorphism globally.
$$\begin{tikzcd}
\mu_2 \inarrow{r} \arrow{d}{=} & \mGal_{F'} \onarrow{r} \arrow{d}{\tilde \iota}& \Gal_{F'} \inarrow{d}{\iota} \\
\mu_2 \inarrow{r} & \mGal_F \onarrow{r} & \Gal_F
\end{tikzcd}$$

Taken together, these inclusions $\tilde \iota \From \mGal_{S'} \Into \mGal_S$ allow one to canonically ``restrict'' metaGalois representations (representations of $\mGal_S$).  

\section{L-groups, parameters, L-functions}

\subsection{L-groups}

We use the term ``L-group'' to refer to a broad class of extensions of Galois groups by complex reductive groups.  Unlike Langlands, Vogan, and others, we do not assume that our L-groups are endowed with a conjugacy class of splittings.  Our L-groups are more closely related to the ``weak E-groups'' of \cite[Definition 3.24]{VoganLLC}.  But we maintain the letter ``L'' since our L-groups are still connected to L-functions.  

The other difference between our L-groups and those in the literature is that (for reasons which will become clear) we consider our L-groups as objects of a 2-category.  A base scheme $S$, of the three classes discussed in the previous section, will be fixed.
\begin{definition}
An \defined{L-group} is a pair $(G^\vee, {}^\EL G)$, where $G^\vee$ is a complex linear algebraic group (not necessarily connected) and ${}^\EL G$ is an extension of locally compact groups
$$G^\vee \Into {}^\EL G \Onto \Gal_S,$$
for which the conjugation action of any element of ${}^\EL G$ on $G^\vee$ is complex-algebraic.
\end{definition}

\begin{remark}
For complex linear algebraic groups, we do not distinguish between the underlying variety and its $\CC$-points.  Thus we say $G^\vee$ is a complex linear algebraic group, and also view $G^\vee$ as a locally compact group.
\end{remark}

Of course, Langlands' L-group ${}^\EL G = \Gal_F \ltimes G^\vee$ (associated to a reductive group $\alg{G}$ over a field $F$) is an example.  When $W$ is a finite-dimensional complex vector space, the direct product $\Gal_S \times GL(W)$ is an L-group.  Since we don't assume $G^\vee$ to be connected, our metaGalois group $\mGal_S$ is another example. 

\begin{definition}
Given two L-groups,
$$G_1^\vee \Into {}^\EL G_1 \Onto \Gal_S, \quad G_2^\vee \Into {}^\EL G_2 \Onto \Gal_S,$$
an \defined{L-morphism} ${}^\EL \rho \From {}^\EL G_1 \To {}^\EL G_2$ will mean a continuous group homomorphism lying over $\Id \From \Gal_S \To \Gal_S$, which restricts to a complex algebraic homomorphism $\rho^\vee \From G_1^\vee \To G_2^\vee.$  An \defined{L-equivalence} will mean an invertible  L-morphism.  
\end{definition}

In other words, an L-morphism fits into a commutative diagram, with the middle column continuous and the left column complex-algebraic.
$$\begin{tikzcd}
G_1^\vee \inarrow{r} \arrow{d}{\rho^\vee}& {}^\EL G_1 \onarrow{r} \arrow{d}{{}^\EL \rho} & \Gal_S \arrow{d}{=} \\
G_2^\vee \inarrow{r} & {}^\EL G_2 \onarrow{r} & \Gal_S
\end{tikzcd}
$$

An \defined{L-representation} of ${}^\EL G$ will mean a pair $({}^\EL \rho, W)$ where $W$ is a finite-dimensional complex vector space and ${}^\EL \rho \From {}^\EL G \To \Gal_S \times GL(W)$ is an L-morphism.  If we project further, from $\Gal_S \times GL(W)$ to $GL(W)$, then an L-representation yields a continuous homomorphism $\rho \From {}^\EL G \To GL(W)$, whose restriction to $G^\vee$ is complex algebraic.  Such a pair $(\rho, W)$ determines an L-representation $({}^\EL \rho, W)$, and so we often abuse notation slightly and say $(\rho, W)$ is an L-representation of ${}^\EL G$.

\begin{definition}
Given two L-morphisms ${}^\EL \rho, {}^\EL \rho' \From {}^\EL G_1 \To {}^\EL G_2$, a natural isomorphism ${}^\EL \rho \xRightarrow{\sim} {}^\EL \rho'$ will mean an element $a \in Z_2^\vee$ (the center of $G_2^\vee$) such that
$${}^\EL \rho'(g) = a \cdot {}^\EL \rho(g) \cdot a^{-1} \text{ for all } g \in {}^\EL G_1.$$
In particular, note that $\rho$ and $\rho'$ coincide on $G_1^\vee$ when they are naturally isomorphic.

This defines the \defined{2-category of L-groups}, L-morphisms, and natural isomorphisms of L-morphisms.
\end{definition}

In many cases of interest (e.g., when ${}^\EL G_2$ arises as the L-group of a split reductive group), the only natural isomorphism is the identity.  However, in some  nonsplit cases, e.g., ${}^\EL G_2 = \Gal_S \ltimes SL_3(\CC)$, the Langlands L-group of a quasisplit $\alg{G} = \alg{PGU}_3$, a nontrivial element $a \in Z_2^\vee$ does not lie in the center of ${}^\EL G_2$.  Such an element $a$ may determine a nonidentity natural isomorphism.

\subsection{Parameters}

Write $\Weil_S$ for the \defined{Weil group}.  When $S = \Spec(F)$ for a local or global field, this Weil group $\Weil_S$ is $\Weil_F$ defined as in \cite{ArtinTate}; when $S = \Spec(\OO)$, we define $\Weil_S$ to be the free cyclic group $\langle \Frob \rangle \isom \ZZ$ generated by a geometric Frobenius $\Frob$.  In all cases, the Weil group is endowed with a continuous homomorphism $\Weil_S \To \Gal_S$ with dense image.

Let $G^\vee \Into {}^\EL G \Onto \Gal_S$ be an L-group.  A \defined{Weil parameter} is a continuous homomorphism $\phi \From \Weil_S \To {}^\EL G$ lying over $\Weil_S \To \Gal_S$, such that $\phi(w)$ is semisimple for all $w \in \Weil_S$ (see \cite[\S 8.2]{BorelCorvallis}).  The reader may follow \cite{BorelCorvallis} and \cite{GrossReeder} to define Weil-Deligne parameters in this general context, when working over a local field.

Write $\Par(\Weil_S, {}^\EL G)$ for the set of ${}^\EL G$-valued Weil parameters.  It is endowed with an action of $G^\vee$ by conjugation:  if $g \in G^\vee$ and $\phi$ is a parameter, then define
$${}^ g \phi(w) = \phi(g^{-1} w g).$$
Two parameters are called \defined{equivalent} if they are in the same $G^\vee$-orbit (this makes sense in the local or global context).

Composition with an L-morphism ${}^\EL \rho \From {}^\EL G_1 \To {}^\EL G_2$ defines a map,
$${}^\EL \rho \From \Par(\Weil_S, {}^\EL G_1) \To \Par(\Weil_S, {}^\EL G_2).$$
Moreover, this map is equivariant, in the sense that for all $g_1 \in G_1^\vee$ and all parameters $\phi \in \Par({}^\EL G_1)$, we have
$${}^\EL \rho \left( {}^{g_1} \phi \right) = {}^{\rho^\vee(g_1)} \left( {}^\EL \rho(\phi) \right).$$
Thus the L-morphism $\rho$ descends to a well-defined map of equivalence classes
$${}^\EL \rho \From \frac{ \Par(\Weil_S, {}^\EL G_1)}{G_1^\vee-\text{conjugation} } \To \frac{\Par(\Weil_S, {}^\EL G_2)}{G_2^\vee-\text{conjugation} }.$$

Next, consider a natural isomorphism of L-morphisms $\rho \xRightarrow{\sim} \rho'$, with $\rho, \rho' \From {}^\EL G_1 \To {}^\EL G_2$.  We find two maps of parameter spaces,
$${}^\EL \rho, {}^\EL \rho' \From \Par(\Weil_S, {}^\EL G_1) \To \Par(\Weil_S, {}^\EL G_2),$$ 
and an element $a \in Z_2^\vee$ such that ${}^\EL \rho'$ is obtained from ${}^\EL \rho$ by conjugation by $a$.

It follows that ${}^\EL \rho$ and ${}^\EL \rho'$ induce the {\em same} map on equivalence classes,
$${}^\EL \rho = {}^\EL \rho' \From \frac{ \Par(\Weil_S, {}^\EL G_1)}{G_1^\vee-\text{conjugation} } \To \frac{\Par(\Weil_S, {}^\EL G_2)}{G_2^\vee-\text{conjugation} }.$$

Suppose that an L-group ${}^\EL G$ is defined up to L-equivalence, and the L-equivalence defined up to unique natural isomorphism.  Then the set of equivalence classes of parameters
$$\frac{\Par(\Weil_S, {}^\EL G)}{G^\vee-\text{conjugation}}$$    
is uniquely defined up to unique isomorphism.

Refinements of the Langlands parameterization suggest that one should look not only at equivalence classes of (Weil or Weil-Deligne) parameters, but also irreducible representations of the component group of the centralizer of a parameter.  Or, following Vogan \cite{VoganLLC}, one can look at $G^\vee$-equivariant perverse sheaves on a suitable variety of parameters.  The fact that conjugation by $a \in Z_2^\vee$ commutes with the conjugation action of $G_2^\vee$ implies that conjugation by $a$ preserves not only the equivalence class of a Weil parameter for ${}^\EL G_2$, but also the equivalence class of such a refined parameter.  If an L-group is defined up to L-equivalence, and the L-equivalence defined up to unique natural isomorphism, then the set of equivalence classes of refined parameters is uniquely defined up to unique isomorphism.

\subsection{L-functions}

Let $G^\vee \Into {}^\EL G \Onto \Gal_S$ be an L-group, and $\phi \From \Weil_S \To {}^\EL G$ a Weil parameter (or we may take $\phi$ to be a Weil-Deligne parameter in the local case).  Let $(\rho, W)$ be an L-representation of ${}^\EL G$.  Then
$$\rho \circ \phi \From \Weil_S \To GL(W)$$
is a Weil representation (or Weil-Deligne representation in the local case).  As such we obtain an L-function (as defined by Weil and discussed in \cite[\S 3.3]{TateCorvallis}),
$$L(\phi, \rho, s) \defeq L(\rho \circ \phi, s).$$
Choosing an additive character $\psi$ as well gives an $\epsilon$-factor (see \cite[\S 3.4]{TateCorvallis}, based on work of Langlands and Deligne),
$$\epsilon(\phi, \rho, \psi, s) \defeq \epsilon(\rho \circ \phi, \psi, s).$$

In the local integral case $S = \Spec(\OO)$, we have $\Weil_S = \langle \Frob \rangle$, and we define the L-functions and $\epsilon$-factors to be those coming from the unramified representation of $\Weil_F$ by pullback.  

In the setting of Langlands L-groups, a zoo of L-representations arises from complex algebraic representations of $G^\vee$, yielding well-known ``standard'' L-functions, symmetric power and exterior power L-functions, etc..  

But in our very broad setting, we limit our discussion to {\em adjoint} L-functions, as these play an important role in representation theory and their definition is ``internal.''  Consider any L-morphism $\rho \From {}^\EL H \To {}^\EL G$ of L-groups.
$$\begin{tikzcd}
H^\vee \inarrow{r} \arrow{d}{\rho^\vee} & {}^\EL H \onarrow{r} \arrow{d}{ {}^\EL \rho} & \Gal_S \arrow{d}{=}  \\
G^\vee \inarrow{r} & {}^\EL G \onarrow{r} & \Gal_S
\end{tikzcd}$$
For example, we might consider the case where $H^\vee$ is a Levi subgroup of $G^\vee$ (as arises in the Langlands-Shahidi method, \cite{Shahidi}).

Let $\Lie{g}^\vee$ be the complex Lie algebra of $G^\vee$.  The homomorphism ${}^\EL \rho$ followed by conjugation gives an adjoint representation:
$$Ad_\rho \From {}^\EL H \To GL\left( \Lie{g}^\vee  \right).$$
Suppose we have a decomposition of $\Lie{g}^\vee$ as a representation of ${}^\EL H$,
\begin{equation}
\label{AdSummands}
\Lie{g}^\vee = \bigoplus_{i = 0}^h \Lie{g}_{i}^\vee.
\end{equation}
For example, when $H^\vee$ is a Levi subgroup of a parabolic $P^\vee \subset G^\vee$, we may decompose $\Lie{g}^\vee$ into $\Lie{h}^\vee$ and the steps in the nilradical of the Lie algebra of $P^\vee$ and its opposite.   

A decomposition \eqref{AdSummands} gives representations $Ad_{i} \From {}^\EL H \To GL\left( \Lie{g}_{i}^\vee  \right)$.  When $\phi \From \Weil_S \To {}^\EL H$ is a Weil parameter, we obtain L-functions
$$L(\phi, Ad_i, s) \defeq L(Ad_i \circ \phi, s).$$
In particular, when ${}^\EL H = {}^\EL G$, and $\rho = \Id$, we write $Ad$ for the adjoint representation of ${}^\EL G$ on $\Lie{g}^\vee$.  This yields \defined{the adjoint L-function} $L(\phi, Ad, s)$ for any Weil parameter $\phi \From \Weil_S \To {}^\EL G$.  When $H^\vee$ is a Levi subgroup of a parabolic in $G^\vee$, and $Ad_i$ arises from a step in the nilradical of the parabolic, we call $L(\phi, Ad_i, s)$ a \defined{Langlands-Shahidi L-function}.

\begin{remark}
The importance of such L-functions for covering groups is suggested by recent work of D. Szpruch \cite{Szp}, who demonstrates that the Langlands-Shahidi construction of L-functions carries over to the metaplectic group.  But it is not clear how to extend the Langlands-Shahidi method to other covering groups, where uniqueness of Whittaker models often fails.  The thesis work of Gao Fan \cite{GaoThesis} takes some promising steps in this direction.  The general machinery of adjoint L-functions also suggests an analogue, for covering groups, of the Hiraga-Ichino-Ikeda conjecture \cite[Conjecture 1.4]{HII} on formal degrees (see Ichino-Lapid-Mao \cite{IchinoLapidMao}).  It is also supported by the simpler observation that theta correspondence for the metaplectic group $\widetilde{Sp}_{2n}$ provides a definition of adjoint L-functions independently of choices of additive characters.
\end{remark}

\subsection{The L-group of a cover}

Now we define the L-group of a cover.  Let $\alg{\tilde G}$ be a degree $n$ cover of a quasisplit group $\alg{G}$ over $S$.  Fix an injective character $\epsilon \From \mu_n \Into \CC^\times$.  Choose a separable closure $\bar F / F$, yielding a geometric base point $\bar s \To S$ and the absolute Galois group $\Gal_S = \pi_1^{\et}(S, \bar s)$.

Recall the constructions of the previous three sections.
\begin{itemize}
\item
$\dgp{\tilde G^\vee}$ denotes the dual group of $\alg{\tilde G}$, a local system on $S_{\et}$ of pinned reductive groups over $\ZZ$,  with center $\dgp{\tilde Z^\vee}$.  It is endowed with a homomorphism $\tau_Q \From \amu_2 \To \dgp{\tilde Z^\vee}$.
\item
$\gerb{E}_\epsilon(\alg{\tilde G})$ is the gerbe associated to $\alg{\tilde G}$, a gerbe on $S_{\et}$ banded by $\sheaf{\tilde Z^\vee} = \dgp{\tilde Z^\vee}(\CC)$.  
\item
$\mu_2 \Into \mGal_S \Onto \Gal_S$ is the metaGalois group.
\end{itemize}

Define $\tilde Z^\vee = \sheaf{\tilde Z}_{\bar s}^\Vee = \dgp{\tilde Z}_{\bar s}^\Vee(\CC)$.  This is the center of the complex reductive group $\tilde G^\vee = \sheaf{\tilde G}_{\bar s}^\Vee = \dgp{\tilde G}_{\bar s}^\Vee(\CC)$.

Pushing out $\mGal_S$ via $\tau_Q \From \mu_2 \To \tilde Z^\vee$ defines an L-group,
\begin{equation}
\label{Twist1}
\tilde Z^\vee \Into (\tau_Q)_\ast \mGal_S \Onto \Gal_S.
\end{equation}
From Theorem \ref{AppxFundGpDefined}, the fundamental group of the gerbe $\gerb{E}_\epsilon(\alg{\tilde G})$ is an L-group, well-defined up to L-equivalence, and the L-equivalence well-defined up to unique natural isomorphism,
\begin{equation}
\label{Twist2}
\tilde Z^\vee \Into \pi_1^{\et}(\gerb{E}_\epsilon(\alg{\tilde G}), \bar s) \Onto \Gal_S.
\end{equation}

\begin{remark}
The extensions \eqref{Twist1} and \eqref{Twist2} play the role of the first and second twist in \cite{MWCrelle}.  In fact \eqref{Twist1} is canonically isomorphic to the first twist in the split case; the extension \eqref{Twist2} may not coincide with the second twist under some circumstances, and the construction here is more general than \cite{MWCrelle} in both cases.
\end{remark}

The Baer sum of \eqref{Twist1} and \eqref{Twist2} is an L-group which will be called ${}^\EL \tilde Z$,
$$\tilde Z^\vee \Into {}^\EL \tilde Z \Onto \Gal_S.$$
\defined{The L-group} of $\alg{\tilde G}$ is defined to be the pushout of ${}^\EL \tilde Z$ via the inclusion $\tilde Z^\vee \Into \tilde G^\vee$,
$$\tilde G^\vee \Into {}^\EL \tilde G \Onto \Gal_S.$$
More explicitly, this pushout is in the $\Gal_S$-equivariant sense.  In other words,
$${}^\EL \tilde G = \frac{\tilde G^\vee \rtimes {}^\EL \tilde Z}{ \langle (z, z^{-1}) : z \in \tilde Z^\vee \rangle },$$
where the semidirect product action ${}^\EL \tilde Z \To \Aut(\tilde G^\vee)$ is given by projection ${}^\EL \tilde Z \To \Gal_S$ followed by the action $\Gal_S \To \Aut(\tilde G^\vee)$ (the Galois group acts by pinned automorphisms on the dual group). 

By construction, ${}^\EL G$ is well-defined by $\alg{\tilde G}$ and $\epsilon$ up to L-equivalence, and the equivalence well-defined up to unique natural isomorphism.

\subsubsection{Local-global compatibility}

Suppose that $\gamma \From F \Into F_v$ is the inclusion of a global field $F$ into its localization at a place.  Let $\alg{\tilde G}$ be a degree $n$ cover of a quasisplit group $\alg{G}$ over $F$, and let $\epsilon \From \mu_n \Into \CC^\times$ be an injective character.  Let $\bar F \Into \bar F_v$ be an inclusion of separable closures, inducing an inclusion $\Gal_{F_v} \Into \Gal_F$ of absolute Galois groups.  Write $S = \Spec(F)$ and $S_v = \Spec(F_v)$, and $\bar s \To S$ and $\bar s_v \To S_v$ for the geometric base points arising from $\bar F \Into \bar F_v$.

Write $\alg{\tilde G}_v$ for the pullback of $\alg{\tilde G}$ via $\Spec(F_v) \To \Spec(F)$.  Similarly, write $Q_v$ for its Brylinski-Deligne invariant.  The results of Section \ref{BaseChangeDualGroup} identify $\tilde G^\vee = \dgp{\tilde G}_{\bar s}^\Vee(\CC)$ with the corresponding dual group for $\alg{\tilde G}_v$ (relative to the separable closure $\bar F_v$).  Thus we simply write $\tilde G^\vee$ for their dual groups and $\tilde Z^\vee$ for the centers thereof.  

The results of Sections \ref{BaseChangeDualGroup} and \ref{CompatibilityMetaGalois} together provide an L-morphism, unique up to unique natural isomorphism.
\begin{equation}
\label{LGTwist1}
\begin{tikzcd}
\tilde Z^\vee \inarrow{r} \arrow{d}{=} & (\tau_{Q_v})_\ast {}^\EL \mGal_{F_v} \onarrow{r} \inarrow{d} & \Gal_{F_v} \inarrow{d} \\
\tilde Z^\vee \inarrow{r} & (\tau_Q)_\ast {}^\EL \mGal_F \onarrow{r} & \Gal_F
\end{tikzcd}
\end{equation}

The results of Section \ref{BaseChangeGerbe} and following Theorem \ref{AppxFundGpDefined} give an L-morphism, unique up to unique natural isomorphism.
\begin{equation}
\label{LGTwist2}
\begin{tikzcd}
\tilde Z^\vee \inarrow{r} \arrow{d}{=} & \pi_1^{\et}(\gerb{E}_\epsilon(\alg{\tilde G}_v), \bar s_v) \onarrow{r} \inarrow{d} & \Gal_{F_v} \inarrow{d} \\
\tilde Z^\vee \inarrow{r} & \pi_1^{\et}(\gerb{E}_\epsilon(\alg{\tilde G}), \bar s) \onarrow{r} & \Gal_F
\end{tikzcd}
\end{equation}

Applying the Baer sum to \eqref{LGTwist1} and \eqref{LGTwist2}, and pushing out via $\tilde Z^\vee \Into \tilde G^\vee$, yields an L-morphism, unique up to unique natural isomorphism.
$$\begin{tikzcd}
\tilde G^\vee \inarrow{r} \arrow{d}{=} & {}^\EL \tilde G_v \onarrow{r} \inarrow{d} & \Gal_{F_v} \inarrow{d} \\
\tilde G^\vee \inarrow{r} & {}^\EL \tilde G \onarrow{r} & \Gal_F
\end{tikzcd}$$

\subsubsection{Parabolic subgroups}

Return to a degree $n$ cover $\alg{\tilde G}$ of a quasisplit group $\alg{G}$ over $S$, and let $\alg{P} \subset \alg{G}$ be a parabolic subgroup defined over $S$.  As before, consider a Levi decomposition $\alg{P} = \alg{M} \alg{N}$ and the resulting cover $\alg{\tilde M}$.  Fix $\epsilon \From \mu_n \Into \CC^\times$.

Compatibility of the dual groups from Section \ref{ParabolicDualGroup} gives inclusions
$$ \tilde Z^\vee \Into \tilde Z_M^\vee \Into \tilde M^\vee \Into \tilde G^\vee,$$
where $\tilde Z^\vee$ denotes the center of $\tilde G^\vee$, and $\tilde Z_M^\vee$ denotes the center of $\tilde M^\vee$.  As these inclusions are compatible with the 2-torsion elements in $\tilde Z^\vee$ and $\tilde Z_M^\vee$, we find an L-morphism.
\begin{equation}
\label{PSTwist1}
\begin{tikzcd}
\tilde Z^\vee \inarrow{r}  \inarrow{d} &  (\tau_Q)_\ast {}^\EL \mGal_F \onarrow{r} \inarrow{d}  & \Gal_S \arrow{d}{=} \\
\tilde Z_M^\vee \inarrow{r} &  (\tau_{Q_{\alg{M}}})_\ast {}^\EL \mGal_F \onarrow{r}  & \Gal_S \\
\end{tikzcd}
\end{equation}
Section \ref{ParabolicGerbe} provided a functor of gerbes $\gerb{i} \From \gerb{E}_\epsilon(\alg{\tilde G}) \To \gerb{E}_\epsilon(\alg{\tilde M})$, lying over $\sheaf{\tilde Z}^\Vee \Into \sheaf{\tilde Z}_{\alg{M}}^\Vee$.  The \'etale fundamental groups give an L-morphism.
\begin{equation}
\label{PSTwist2}
\begin{tikzcd}
\tilde Z^\vee \inarrow{r} \inarrow{d} & \pi_1^{\et}(\gerb{E}_\epsilon(\alg{\tilde G}), \bar s) \onarrow{r} \inarrow{d} & \Gal_S \arrow{d}{=} \\
\tilde Z_M^\vee \inarrow{r} & \pi_1^{\et}(\gerb{E}_\epsilon(\alg{\tilde M}), \bar s) \onarrow{r}  & \Gal_S \\
\end{tikzcd}
\end{equation}

Applying the Baer sum to \eqref{PSTwist1} and \eqref{PSTwist2} yields an L-morphism, unique up to unique natural isomorphism.
$$\begin{tikzcd}
\tilde Z^\vee \inarrow{r} \inarrow{d} & {}^\EL \tilde Z \onarrow{r} \inarrow{d} & \Gal_S \arrow{d}{=} \\
\tilde Z_M^\vee \inarrow{r} & {}^\EL \tilde Z_M \onarrow{r} & \Gal_S \\
\end{tikzcd}$$
The universal property of pushouts yields an L-morphism.
$$\begin{tikzcd}
\tilde M^\vee \inarrow{r} \inarrow{d}& {}^\EL \tilde M \onarrow{r} \inarrow{d} & \Gal_S  \arrow{d}{=} \\
\tilde G^\vee \inarrow{r}  & {}^\EL \tilde G \onarrow{r} & \Gal_S
\end{tikzcd}$$
This L-morphism is well-defined up to conjugation by $\tilde Z_M^\vee$.

\appendix

\section{Torsors, gerbes, and fundamental groups}
\label{GroupsTorsors}
Let $S$ be a connected scheme, and $S_{\et}$ the \'etale site.  Our treatment of sheaves on $S_{\et}$ follows \cite[\S II]{SGA4half}.  Recall that a geometric point of $S$ is a morphism of schemes $\bar s \From \Spec(\bar F) \To S$, where $\bar F$ is a separably closed field.  

An \defined{open \'etale neighborhood} of $\bar s$ is an \'etale morphism $U \To S$ endowed with a lift $\bar u \From \Spec(\bar F) \To U$ of the geometric point $\bar s$.  If $\bar s$ is a geometric point, we write $\pi_1^{\et}(S, \bar s)$ for the \'etale fundamental group.  When $\bar s$ is fixed, we define $\Gal_S = \pi_1^{\et}(S, \bar s)$.

\subsection{Local systems on $S_{\et}$}

\begin{definition}
A \defined{local system on $S_{\et}$} is a locally constant sheaf $\sheaf{J}$ of sets on $S_{\et}$.
\end{definition}
When $\sheaf{J}$ is a local system on $S_{\et}$ and $U \To S$ is \'etale, we write $\sheaf{J}[U]$ for the set of sections over $U$ and we write $\sheaf{J}_U$ for the local system on $U_{\et}$ obtained by restriction.  If $\bar s$ is a geometric point of $S$, then the fibre $\sheaf{J}_{\bar s}$ is the inductive limit  $\limdir_U \sheaf{J}[U]$, over open \'etale neighborhoods of $\bar s$.  By local constancy, $\sheaf{J}_{\bar s} = \sheaf{J}[U]$ for some such open \'etale neighborhood.  Often in this paper we work locally on $S_{\et}$ and abuse notation a bit by writing $j \in \sheaf{J}$ rather than $j \in \sheaf{J}[U]$ (for an  \'etale $U \To S$).  

More generally, if $\Cat{C}$ is a category, then one may work with local systems on $S_{\et}$ of objects of $\Cat{C}$, or ``$\Cat{C}$-valued local systems.''  If $\Cat{C} \To \Cat{D}$ is a functor, then one finds a corresponding functor from the category of $\Cat{C}$-valued local systems to the category of $\Cat{D}$-valued local systems. Fibres of such $\Cat{C}$-valued local systems over geometric points make sense in this generality, by local constancy.

\begin{example}
\label{LocSpec}
Let $\sheaf{M}$ be a local system on $S_{\et}$ of finitely-generated abelian groups, and let $R$ be a commutative ring.  Then $\Spec(R[\sheaf{M}])$ will denote the local system on $S_{\et}$ of affine group schemes over $R$ given by
$$\Spec(R[\sheaf{M}])[U] = \Spec(R[\sheaf{M}[U]]).$$
\end{example}
We will work with local systems of groups, local systems of affine group schemes over $\ZZ$, local systems of root data, etc..

\subsection{Torsors on $S_{\et}$}
\begin{definition}
Let $\sheaf{G}$ be a sheaf of groups on $S_{\et}$.  A \defined{$\sheaf{G}$-torsor} is a sheaf of sets $\sheaf{V}$ on $S_{\et}$, endowed with an action $\ast \From \sheaf{G} \times \sheaf{V} \To \sheaf{V}$, such that
$$\sheaf{G} \times \sheaf{V} \To \sheaf{V} \times \sheaf{V}, \quad (g,v) \mapsto (g \ast v, v)$$
is an isomorphism of sheaves of sets on $S_{\et}$.  Morphisms of $\sheaf{G}$-torsors are morphisms of sheaves on $S_{\et}$ which intertwine the $\sheaf{G}$-action.  The \defined{category of $\sheaf{G}$-torsors} will be denoted $\Cat{Tors}(\sheaf{G})$.
\end{definition}

If $\sheaf{V}$ is a $\sheaf{G}$-torsor, we write $[ \sheaf{V} ]$ for its isomorphism class.  The isomorphism classes of $\sheaf{G}$-torsors form a pointed set denoted $H_{\et}^1(S, \sheaf{G})$.  

%Similarly, the isomorphism classes of $\dgp{G}$-torsors form a pointed set denoted $H^1(S_{\et}, \dgp{G})$.
%Suppose that $f \From \sheaf{G} \To \sheaf{H}$ is a homomorphism of sheaves of groups on $S_{\et}$.  If $\sheaf{V}$ is a $\sheaf{G}$-torsor, then we write $f_\ast \sheaf{V}$ for the pushout,
%$$f_\ast \sheaf{V} \defeq \frac{ \sheaf{V} \times \sheaf{H} }{ (g \ast v, h) \sim (v, f(g) \cdot h) }.$$
%This operation of \defined{pushing out torsors} defines a functor,
%$$f_\ast \From \Cat{Tors}(\sheaf{G}) \To \Cat{Tors}(\sheaf{H}).$$
%Similarly, such a pushout functor can be defined for $\dgp{G}$-torsors when $\dgp{G} \To \dgp{H}$ is a homomorphism of sheaves of affine group $R$-schemes on $S_{\et}$.

The category of torsors has more structure in the abelian case.  If $\sheaf{A}$ is a sheaf of {\em abelian} groups on $S_{\et}$, then the category $\Cat{Tors}(\sheaf{A})$ inherits a monoidal structure.  Namely, if $\sheaf{V}_1$ and $\sheaf{V}_2$ are two $\sheaf{A}$-torsors, define
$$\sheaf{V}_1 \Baer \sheaf{V}_2 = \frac{\sheaf{V}_1 \times \sheaf{V}_2}{ (a \ast v_1, v_2) \sim (v_1, a \ast v_2)}.$$
With this monoidal structure, the trivial torsor $\sheaf{A}$ as zero object, and obvious isomorphisms for commutativity and associativity and unit, the category $\Cat{Tors}(\sheaf{A})$ becomes a Picard groupoid (i.e., a strictly commutative Picard category, in the terminology of \cite[Expos\'e XVIII \S 1.4]{SGA4T3}.  The pointed set of isomorphism classes $H_{\et}^1(S, \sheaf{A})$ becomes an abelian group, with
$$[\sheaf{V}_1] + [\sheaf{V}_2] \defeq \left[ \sheaf{V}_1 \Baer \sheaf{V}_2 \right].$$
The group $H_{\et}^1(S, \sheaf{A})$ is identified with the \'etale cohomology with coefficients $\sheaf{A}$.

Suppose that $f \From \sheaf{A} \To \sheaf{G}$ is a homomorphism of sheaves of groups on $S_{\et}$, with $\sheaf{A}$ abelian, and $f$ central (i.e., $f$ factors through the inclusion of the center $\sheaf{Z} \Into \sheaf{G}$).  If $\sheaf{V}$ is an $\sheaf{A}$-torsor, then we write $f_\ast \sheaf{V}$ for the pushout,
$$f_\ast \sheaf{V} \defeq \frac{ \sheaf{V} \times \sheaf{G} }{ (a \ast v, g) \sim (v, f(a) \cdot g) }.$$
This operation of \defined{pushing out torsors} defines a functor,
$$f_\ast \From \Cat{Tors}(\sheaf{A}) \To \Cat{Tors}(\sheaf{G}).$$

Suppose that $c \From \sheaf{A}_1 \To \sheaf{A}_2$ is a homomorphism of sheaves of abelian groups on $S_{\et}$.  Let $\sheaf{V}_1$ be an $\sheaf{A}_1$-torsor and $\sheaf{V}_2$ an $\sheaf{A}_2$-torsor.  A map of torsors $\pi \From \sheaf{V}_1 \To \sheaf{V}_2$ \defined{lying over} $c$ means a morphism of sheaves of sets on $S_{\et}$ satisfying
$$\pi(a_1 \ast v_1) = c(a_1) \ast \pi(v_1), \text{ for all } a_1 \in \sheaf{A}_1, v_1 \in \sheaf{V}_1.$$
Such a map factors uniquely through $c_\ast \sheaf{V}_1$.  

A short exact sequence of sheaves of abelian groups on $S_{\et}$,
\begin{equation}
\label{ABC}
\sheaf{A} \xhookrightarrow{\alpha}\sheaf{B} \xtwoheadrightarrow{\beta} \sheaf{C},
\end{equation}
yields two more constructions of torsors.

First, the sequence yields a boundary map in cohomology, $\partial \From H_{\et}^0(S, \sheaf{C}) \To H_{\et}^1(S, \sheaf{A})$.  There is a corresponding map from global sections of $\sheaf{C}$ to objects of the category of $\sheaf{A}$-torsors as follows.

Begin with $c \in \sheaf{C}[S]$ and write $[c]$ to consider it as an element of $H_{\et}^0(S, \sheaf{C})$.  For any \'etale $U \To S$, write $c_U \in \sheaf{C}[U]$ for the restriction of $c$ to $U$.  Define $\partial c$ to be the sheaf on $S_{\et}$ whose sections are given by
$$\partial c[U] = \{ b \in \sheaf{B}[U] :  \beta(b) = c_U \}.$$
The sheaf $\partial c$ is naturally an $\sheaf{A}$-torsor; the equivalence class $[\partial c] \in H_{\et}^1(S, \sheaf{A})$ coincides with $\partial [c]$.  The sheaf $\partial c$ is called the \defined{torsor of liftings} of $c$ via $\beta$.

Next, write $\shom(\sheaf{C}, \sheaf{B})$ for the sheaf of homomorphisms (``sheaf-hom'') from $\sheaf{C}$ to $\sheaf{B}$.  This is a sheaf of abelian groups on $S_{\et}$, and there is a subsheaf of sets $\sspl(\sheaf{B})$ consisting of those homomorphisms which split the extension \eqref{ABC}.  This subsheaf $\sspl(\sheaf{B})$ is naturally a $\shom(\sheaf{C}, \sheaf{A})$-torsor, called the \defined{torsor of splittings}.

\subsection{Gerbes on $S_{\et}$}

Here we introduce a class of gerbes on $S_{\et}$.  In what follows, let $\sheaf{A}$ be a sheaf of {\em abelian} groups on $S_{\et}$.
\begin{definition}
A \defined{gerbe on $S_{\et}$ banded by $\sheaf{A}$} is a (strict) stack $\gerb{E}$ on $S_{\et}$ of groupoids such that $\gerb{E}$ is locally nonempty, locally connected, and banded by $\sheaf{A}$.
\end{definition}
We unravel this definition here, beginning with the data.
\begin{description}
\item[a (strict) stack $\gerb{E}$ on $S_{\et}$ of groupoids]  
For each \'etale $U \To S$, we have a (possibly empty) groupoid $\gerb{E}[U]$.  For $\gamma \From U' \To U$, a morphism of schemes \'etale over $S$, we have a pullback functor $\gamma^\ast \From \gerb{E}[U] \To \gerb{E}[U']$.  
\item[banded by $\sheaf{A}$]
For every object $x$ of $\gerb{E}[U]$, there is given an isomorphism $\sheaf{A}[U] \To \Aut(x)$ (written $\alpha \mapsto \alpha_x$).
\end{description}
This data satisfies additional axioms:
\begin{description}
\item[(strict) stack axioms]
For each pair $\gamma \From U' \To U$ and $\delta \From U'' \To U'$, we require {\em equality} of functors $\delta^\ast \circ \gamma^\ast = (\gamma \delta)^\ast$.  (``Strictness'' refers to the requirement of equality rather than extra data of natural isomorphisms).  Descent for objects and morphisms is effective.  
\item[locally nonempty]
There exists a finite \'etale $U \To S$ such that $\gerb{E}[U]$ is nonempty.
\item[locally connected]
For any \'etale $U \To S$ and pair of objects $x,y \in \gerb{E}[U]$, there exists a finite \'etale $\gamma \From U' \To U$ such that $\gamma^\ast x$ is isomorphic to $\gamma^\ast y$ in $\gerb{E}[U']$.
\item[banding]
Given a morphism $f \From x \To y$ in $\gerb{E}[U]$, and $\alpha \in \sheaf{A}[U]$, $\alpha_y \circ f = f \circ \alpha_x$.  Also, given $\gamma \From U' \To U$ \'etale, $\alpha_{\gamma^\ast x} = \gamma^\ast \alpha_x \in \Aut(\gamma^\ast x)$.   
\end{description}

\begin{remark}
We will not require gerbes banded by nonabelian groups -- commutativity greatly simplifies the theory.  For a fuller treatment of gerbes, one can consult the original book of Giraud \cite{Giraud}, work of Breen \cite{BreenNotes}, \cite{BreenOriginal}, the book of Brylinski \cite[Chapter V]{BrylBook}, the article of Deligne \cite{DelModere}, the introduction of Garland and Patnaik \cite{Garl}, and the Stacks Project \cite{Stacks}, among others.  We don't keep track of universes along the way, while Giraud \cite{Giraud} is careful about set-theoretic subtleties.  Our ``strictness'' assumption is typically referred to as an assumption that the fibred category $\gerb{E} \To S_{\et}$ is ``split'' (see \cite{Vistoli}).  The strictness assumption is not so restrictive, since every stack is equivalent to a strict stack (cf. \cite[Theorem 3.45]{Vistoli}).  See also \cite[\S I.1]{Giraud}.
\end{remark}

If $\gerb{E}$ is a gerbe on $S_{\et}$ banded by $\sheaf{A}$, and $U \To S$ is \'etale, then we write $\gerb{E}_U$ for its restriction to $U_{\et}$; this is a gerbe on $U_{\et}$ banded by $\sheaf{A}_U$.

Our strictness assumption allows us to easily define the \defined{fibre} of a gerbe $\gerb{E}$ at a geometric point $\bar s$.  This is the category $\gerb{E}_{\bar s}$ whose object set is the direct limit $\limdir_{U} \gerb{E}[U]$ of object sets, indexed by open \'etale neighborhoods of $\bar s$.  Write $A = \sheaf{A}_{\bar s}$ for the fibre of $\sheaf{A}$ over $\bar s$.  If $\bar z_1, \bar z_2$ are objects of $\gerb{E}_{\bar s}$, then $\Hom(\bar z_1, \bar z_2)$ is naturally an $A$-torsor.

If $x,y \in \gerb{E}[U]$, then for $\gamma \From U' \To U$ define
$$\shom(x,y)[U'] = \Hom(\gamma^\ast x, \gamma^\ast y).$$
In this way, we construct a sheaf $\shom(x,y)$ of sets on $U_{\et}$.  In particular, we find a sheaf of groups $\sAut(x)$ on $U_{\et}$.  The banding provides an isomorphism of sheaves of groups, $\sheaf{A}_U \xrightarrow{\sim} \sAut(x)$, and $\shom(x,y)$ becomes an $\sheaf{A}_U$-torsor.

\subsubsection{Functors of gerbes}

Suppose that $c \From \sheaf{A}_1 \To \sheaf{A}_2$ is a homomorphism of sheaves of abelian groups on $S_{\et}$, and $\gerb{E}_1, \gerb{E}_2$ are gerbes on $S_{\et}$ banded by $\sheaf{A}_1, \sheaf{A}_2$, respectively.  A \defined{functor of gerbes} $\phi \From \gerb{E}_1 \To \gerb{E}_2$, \defined{lying over $c$}, is a (strict) functor of stacks lying over $c$.  This entails the following.
\begin{description}
\item[(strict) functor of stacks]
For each \'etale $U \To S$, a functor of categories $\phi[U] \From \gerb{E}_1[U] \To \gerb{E}_2[U]$.  For every $\gamma \From U' \To U$, with pullback functors $\gamma_1^\ast$ in $\gerb{E}_1$ and $\gamma_2^\ast$ in $\gerb{E}_2$, the ``strictness'' condition requires an {\em equality} of functors, $\gamma_2^\ast \circ \phi[U] = \phi[U'] \circ \gamma_1^\ast$.
\item[lying over $c$]
The ``lying over $c$'' condition requires that, for each $\alpha_1 \in \sheaf{A}_1[U]$ with $\alpha_2 = c(\alpha_1)$, and object $x_1 \in \gerb{E}_1[U]$ with $x_2 = \phi[U] x_1$, we have
$$(\alpha_2)_{x_2} = \phi[U] \left( (\alpha_1)_{x_1} \right) \in \Aut(x_2).$$
\end{description} 

Gerbes on $S_{\et}$ banded by a fixed sheaf of abelian groups $\sheaf{A}$ form a 2-category; if $\gerb{E}_1$ and $\gerb{E}_2$ are two such gerbes banded by the same $\sheaf{A}$, an \defined{equivalence of gerbes} $\phi \From \gerb{E}_1 \To \gerb{E}_2$ is a functor of gerbes lying over $\Id \From \sheaf{A} \To \sheaf{A}$.  Given two such equivalences of gerbes $\phi, \phi' \From \gerb{E}_1 \To \gerb{E}_2$, a natural isomorphism $\phi \xRightarrow{\sim} \phi'$ consists of natural isomorphisms of functors $\phi[U] \Rightarrow \phi'[U]$ for each $U$, compatible with pullback.  This defines a 2-category of gerbes banded by $\sheaf{A}$, equivalences, and natural isomorphisms of equivalences. 

%Given an equivalence of gerbes $\phi \From \gerb{E}_0 \To \gerb{E}$, there is a pair $(\phi^{-1}, \nu)$ in which $\phi^{-1} \From \gerb{E} \To \gerb{E}_0$ is an equivalence of gerbes and $\nu \From \phi^{-1} \circ \phi \Rightarrow \Id$ is a natural isomorphism.  The pair $(\phi^{-1}, \nu)$ is unique up to unique natural isomorphism; given another such pair $(\phi_\ast^{-1}, \nu_\ast)$, there exists a unique natural isomorphism $\rho \From \phi^{-1} \Rightarrow \phi_\ast^{-1}$ such that the following diagram of functors and natural isomorphisms commutes.
%$$\begin{tikzcd}
%\phi^{-1} \circ \phi \arrow[Rightarrow]{r}{\rho} \arrow[Rightarrow]{d}{\nu} & \phi_\ast^{-1} \circ \phi \arrow[Rightarrow]{d}{\nu_\ast} \\
%\Id \arrow[Rightarrow]{r}{\Id} & \Id
%\end{tikzcd}$$

Given two gerbes $\gerb{E}_1, \gerb{E}_2$ banded by $\sheaf{A}$, one may ``contract'' them to form another gerbe $\gerb{E}_1 \Baer \gerb{E}_2$ banded by $\sheaf{A}$.    The family of categories of torsors, $\stors(\sheaf{A})$, given by $\stors(\sheaf{A})[U] = \Cat{Tors}(\sheaf{A}_U)$ (for each \'etale $U \To S$), with pullbacks given by restriction of sheaves, forms the \defined{neutral} $\sheaf{A}$-gerbe on $S_{\et}$.  $\gerb{E} \Baer \stors(\sheaf{A})$ is equivalent to $\gerb{E}$, for any gerbe $\gerb{E}$ banded by $\sheaf{A}$ (and the equivalence is determined up to unique natural isomorphism).

Suppose that $\gerb{E}$ is a gerbe on $S_{\et}$ banded by $\sheaf{A}$, and $x$ is an object of $\gerb{E}[U]$ for some \'etale $U \To S$.  Then, for $\gamma \From U' \To U$ \'etale, and $y \in \gerb{E}[U']$, we have a $\sheaf{A}_{U'}$-torsor $\shom(\gamma^\ast x, y)$.  This map $y \mapsto \shom(\gamma^\ast x, y)$ extends to an equivalence of gerbes from $\gerb{E}_U$ to $\stors(\sheaf{A})_U$.  In this way, we say that $x$ \defined{neutralizes} the gerbe $\gerb{E}$ over $U$.

If $\gerb{E}$ is a gerbe banded by $\sheaf{A}$, we write $\left[ \gerb{E} \right]$ for its equivalence class.  The set of such equivalence classes is denoted $H_{\et}^2(S, \sheaf{A})$.  This forms an abelian group, with zero corresponding to the neutral gerbe of $\sheaf{A}$-torsors, and addition arising from contraction.  From \cite{Giraud}, we identify $H_{\et}^2(S, \sheaf{A})$ with the \'etale cohomology of $S$ with coeffiencts $\sheaf{A}$.

\subsubsection{Pushouts}

Given a gerbe $\gerb{E}_1$ banded by $\sheaf{A}_1$, and $c \From \sheaf{A}_1 \To \sheaf{A}_2$ as above, one may construct a gerbe $c_\ast \gerb{E}_1$ banded by $\sheaf{A}_2$ called the \defined{pushout} of $\gerb{E}_1$ by $c$.  Any functor of gerbes $\phi \From \gerb{E}_1 \To \gerb{E}_2$ lying over $c$ factors through a functor $c_\ast \gerb{E}_1 \To \gerb{E}_2$ (the functor being determined uniquely up to unique natural isomorphism, see \cite[\S 5.3]{DelModere}).  The objects of $c_\ast \gerb{E}$ are the same as those of $\gerb{E}$.  But, given two such objects $x,y \in \gerb{E}[U]$, the morphism set $\Hom_{c_\ast \gerb{E}}(x,y)$ is defined as the pushout of torsors,
$$\Hom_{c_\ast \gerb{E}}(x,y) = c_\ast \Hom_{\gerb{E}}(x,y).$$
If $x,y \in \gerb{E}[U]$, $j_0 \in \Hom_{\gerb{E}}(x,y)$, $\alpha \in \sheaf{A}_2[U]$, we write $f \wedge \alpha$ for the resulting morphism from $x$ to $y$ in $c_\ast \gerb{E}$.

The pushout of gerbes corresponds to the map in cohomology,
$$H_{\et}^2(S, \sheaf{A}_1) \To H_{\et}^2(S, \sheaf{A}_2), \quad [\gerb{E}_1] \mapsto c_\ast [\gerb{E}_1].$$
See \cite[Chapitre IV, \S 3.3, 3.4]{Giraud} for details.

\subsubsection{The gerbe of liftings}

\label{AppendixGerbeLiftings}

If $\sheaf{A} \xhookrightarrow{\alpha}\sheaf{B} \xtwoheadrightarrow{\beta} \sheaf{C}$ is a short exact sequence of sheaves of abelian groups on $S_{\et}$, then the sequence of cohomology groups,
$$H_{\et}^2(S, \sheaf{A}) \xrightarrow{\alpha} H_{\et}^2(S, \sheaf{B}) \xrightarrow{\beta} H_{\et}^2(S, \sheaf{C})$$
is also exact.  The analogous construction with gerbes is the following:  suppose that $\gerb{E}$ is a gerbe banded by $\sheaf{B}$, $\gerb{F}$ is a gerbe banded by $\sheaf{C}$, and $\gerb{p} \From \gerb{E} \To \gerb{F}$ is a functor of gerbes lying over $\sheaf{B} \xrightarrow{\beta} \sheaf{C}$.  If $\gerb{z}$ is an $S$-object of $\gerb{F}$ (neutralizing $\gerb{F}$, so that $0 = [\gerb{F}] \in H_{\et}^2(S, \sheaf{C})$), then cohomology suggests that $\gerb{E}$ arises as the pushout of a gerbe banded by $\sheaf{A}$.

Indeed, we define the gerbe $\gerb{p}^{-1}(\gerb{z})$ as follows:  the objects of $\gerb{p}^{-1}(\gerb{z})[U]$ are pairs $(\gerb{y}, j)$ where $\gerb{y}$ is an object of $\gerb{E}$[U], and $j \From \gerb{p}(\gerb{y}) \To \gerb{z}$ is an isomorphism in $\gerb{F}[U]$.  The morphisms in $\gerb{p}^{-1}(\gerb{z})$ are those in $\gerb{E}$ which are compatible with the isomorphisms to $\gerb{z}$.  The gerbe $\gerb{p}^{-1}(\gerb{z})$ will be called the \defined{gerbe of liftings} of $\gerb{z}$ via $p$.

\subsubsection{The gerbe of $n^{\th}$ roots}
\label{AppendixGerbeRoots}
Suppose that $\sheaf{C}$ is a sheaf of abelian groups on $S_{\et}$, and the multiplication-by-$n$ homomorphism is surjective from $\sheaf{C}$ to itself.  An important example of a gerbe of liftings arises from the Kummer sequence
$$\sheaf{C}_{[n]} \Into \sheaf{C} \xtwoheadrightarrow{n} \sheaf{C}.$$
Pushing out via multiplication by $n$ gives to a functor of gerbes,
$$n_\ast \From \stors(\sheaf{C}) \To \stors(\sheaf{C}),$$
lying over $\sheaf{C} \xtwoheadrightarrow{n} \sheaf{C}$.  Given a $\sheaf{C}$-torsor $\sheaf{V}$, the gerbe of liftings of $\sheaf{V}$ via $n_\ast$ will be called the \defined{gerbe of $n^{\th}$ roots}, denoted $\sqrt[n]{\sheaf{V}}$.  It is a gerbe on $S_{\et}$ banded by $\sheaf{C}_{[n]}$.  The map which sends a $\sheaf{C}$-torsor to its gerbe of $n^{\th}$ roots corresponds to the Kummer coboundary map $\Kum \From H_{\et}^1(S, \sheaf{C}) \To H_{\et}^2(S, \sheaf{C}_{[n]})$.

Explicitly, an object of $\sqrt[n]{\sheaf{V}}$ is a pair $(\sheaf{H}, h)$ where $\sheaf{H}$ is a $\sheaf{C}$-torsor, and $h \From \sheaf{H} \To \sheaf{V}$ is a morphism of sheaves making the following diagram commute.
$$\begin{tikzcd}
\sheaf{C} \times \sheaf{H} \arrow{r}{\ast} \arrow{d}{n \times h} & \sheaf{H} \arrow{d}{h} \\
\sheaf{C} \times \sheaf{V} \arrow{r}{\ast} & \sheaf{V}
\end{tikzcd}$$

The construction of the gerbe of $n^{\th}$ roots is itself functorial.  Consider a homomorphism of sheaves of abelian groups, $c \From \sheaf{C}_1 \To \sheaf{C}_2$, and assume that $\sheaf{C}_1 \xrightarrow{n} \sheaf{C}_1$ and $\sheaf{C}_2 \xrightarrow{n} \sheaf{C}_2$ are surjective.  Suppose that $\sheaf{V}_1$ is a $\sheaf{C}_1$-torsor, and $\sheaf{V}_2$ is a $\sheaf{C}_2$-torsor.  Suppose that $f \From \sheaf{V}_1 \To \sheaf{V}_2$ is a morphism of torsors lying over the homomorphism $c \From \sheaf{C}_1 \To \sheaf{C}_2$.  Then we find a functor of gerbes $\sqrt[n]{f} \From \sqrt[n]{\sheaf{V}_1} \To \sqrt[n]{\sheaf{V}_2}$ lying over the homomorphism of bands $\sheaf{C}_{1,[n]} \To \sheaf{C}_{2,[n]}$.

\subsection{Fundamental group}
Now let $\sheaf{A}$ be a \textbf{local system} of abelian groups on $S_{\et}$.  Let $\gerb{E}$ be a gerbe on $S_{\et}$ banded by $\sheaf{A}$.  Let $\bar F$ be a separably closed field, let $\bar s \From \Spec(\bar F) \To S$ be a geometric point, and recall that $\Gal_S = \pi_1^{\et}(S,\bar s)$ denotes the \'etale fundamental group.  

Suppose that $U \To S$ is a Galois cover, and $\bar u \From \Spec(\bar F) \To U$ lifts the geometric point $\bar s$.  Writing $\Gal_U = \pi_1^{\et}(U, \bar u)$, we find a short exact sequence
$$\Gal_U \Into \Gal_S \Onto \Gal(U / S).$$
If $\gamma \in \Gal_S$ we write $\gamma_U$ for its image in $\Gal(U / S)$.  

Suppose moreover that $\sheaf{A}_U$ is a constant sheaf and $\gerb{E}[U]$ is a nonempty groupoid (i.e., $\gerb{E}$ is neutral over $U$).  Write $A = \sheaf{A}[U]$ for the resulting abelian group.  Then $A$ is endowed with an action of $\Gal_S$ that factors through the finite quotient $\Gal(U/S)$.  An object $z \in \gerb{E}[U]$ will be called a \defined{base point} for the gerbe $\gerb{E}$ (over $U$).  The banding identifies $A$ with the automorphism group of $z$.

Without loss of generality, pulling back to a larger Galois cover if necessary, we may assume that $\Hom(z, \gamma_U^\ast z)$ is nonempty for all $\gamma \in \Gal_S$.  In this way, the base point $z \in \gerb{E}[U]$ and $\gamma \in \Gal_S$ define an $A$-torsor,
$$\Aut_\gamma(z) \defeq \Hom(z, \gamma_U^\ast z).$$
We write $\gamma^\ast$ instead of $\gamma_U^\ast$, when there is little choice of confusion.  

Define the \defined{\'etale fundamental group of the gerbe} $\gerb{E}$, at the base point $z$, by
$$\pi_1(\gerb{E}, z) = \bigsqcup_{\gamma \in \Gal_S} \Aut_\gamma(z).$$
The group structure is given, for $\gamma_1, \gamma_2 \in \Gal_S$, by the following sequence.
\begin{align*}
\Aut_{\gamma_1}(z) \times \Aut_{\gamma_2}(z) &= \Hom(z, \gamma_1^\ast z) \times \Hom(z, \gamma_2^\ast z)  \\
& \xrightarrow{\gamma_2^\ast \times \Id} \Hom(\gamma_2^\ast z, \gamma_2^\ast \gamma_1^\ast z) \times \Hom(z, \gamma_2^\ast z) \\
& \xrightarrow{\circ} \Hom(z, \gamma_2^\ast  \gamma_1^\ast z) \\
& \xrightarrow{=} \Hom(z, (\gamma_1 \gamma_2)^\ast z) = \Aut_{\gamma_1 \gamma_2}(z)
\end{align*}
As $\Aut_{\Id}(z) = \Aut(z)$, the isomorphism $A \xrightarrow{\sim} \Aut(z)$, $\alpha \mapsto \alpha_z$, gives an extension of groups,
\begin{equation}
\label{FundGpGerbe}
A \Into \pi_1^{\et}(\gerb{E}, z) \Onto \Gal_S.
\end{equation}
If $\gamma \in \Gal_U \subset \Gal_S$, then $\gamma_U = \Id$ and so $\Aut_\gamma(z) = \Aut_{\Id}(z)$.  In this way, we find a splitting $\Gal_U \Into \pi_1(\gerb{E}, z)$.  In other words, $\pi_1^{\et}(\gerb{E}, z)$ arises as the pullback of an extension of $\Gal(U/S)$ by $A$.  The conjugation action of $\Gal_S$ on $A$, in the extension \eqref{FundGpGerbe}, coincides with the canonical action of $\Gal(U/S)$ on $A = \sheaf{A}[U]$.

The sequence \eqref{FundGpGerbe} describes the fundamental group of a gerbe (with base point) as an extension of $\Gal_S$ by $A$.  Here we analyze how this fundamental group depends on the choice of base point, and how it behaves under equivalence of gerbes.  

Consider a further Galois cover $\delta \From U' \To U$ and geometric base point $\bar u'$ lifting $\bar u$.  By constancy of $\sheaf{A}_U$, we identify $A = \sheaf{A}[U] = \sheaf{A}[U']$.  For all $\gamma \in \Gal_S$, we have $\gamma_U \circ \delta = \delta \circ \gamma_{U'}$.  This defines an isomorphism of $A$-torsors, $\delta^\ast \From \Aut_\gamma(z) \xrightarrow{\sim} \Aut_\gamma(\delta^\ast z)$, using the sequence below.
\begin{align*}
\Aut_\gamma(z) = \Hom(z, \gamma_U^\ast z) & \xrightarrow{\delta^\ast} \Hom( \delta^\ast z, \delta^\ast \gamma_U^\ast z) \\
& \xrightarrow{=} \Hom \left( \delta^\ast z, (\gamma_U \delta)^\ast z \right) \\
& \xrightarrow{=} \Hom \left( \delta^\ast z, (\delta \gamma_{U'})^\ast z \right) \\
& \xrightarrow{=} \Hom \left( \delta^\ast z, \gamma_{U'}^\ast \delta^\ast z \right)  = \Aut_\gamma(\delta^\ast z).
\end{align*}
Putting these isomorphisms together, we find an isomorphism of extensions.
$$\begin{tikzcd}
A \inarrow{r} \arrow{d}{=} & \pi_1^{\et}(\gerb{E}, z) \onarrow{r} \arrow{d}{\iota_\delta} & \Gal_S \arrow{d}{=} \\
A \inarrow{r} & \pi_1^{\et}(\gerb{E}, \delta^\ast z) \onarrow{r} & \Gal_S
\end{tikzcd}$$

A further cover $\delta' \From U'' \To U'$, with $\delta'' = \delta \circ \delta' \From U'' \To U$, gives a commutative diagram in the category of extensions of $\Gal_S$ by $A$.
$$\begin{tikzcd}
\pi_1^{\et}(\gerb{E}, z) \arrow{r}[swap]{\iota_\delta} \arrow[bend left=20]{rr}{\iota_{\delta''}} & \pi_1^{\et}(\gerb{E}, \delta^\ast z)  \arrow{r}[swap]{\iota_{\delta'}} & \pi_1^{\et}(\gerb{E}, (\delta')^\ast \delta^\ast z) = \pi_1^{\et}(\gerb{E}, (\delta'')^\ast z)  \end{tikzcd}$$

Define $\bar z \in \gerb{E}_{\bar s}$ to be the image of the base point $z$ in the direct limit.  We call $\bar z$ a \defined{geometric base point} for the gerbe $\gerb{E}$.  Define
$$\pi_1^{\et}(\gerb{E}, \bar z) = \limdir_{U'} \pi_1^{\et}(\gerb{E}, \delta^\ast z),$$
the direct limit over Galois covers $\delta \From (U', \bar u') \To (U, \bar u)$, via the isomorphisms $\iota_\delta$ described above.  This gives an extension of groups, depending (up to unique isomorphism) only on the {\em geometric} base point $\bar z \in \gerb{E}_{\bar s}$.
$$A \Into \pi_1^{\et}(\gerb{E}, \bar z) \Onto \Gal_S.$$
This extension is also endowed with a family of splittings over finite-index subgroups $\Gal_U \subset \Gal_S$, arising from base points $z \in \gerb{E}[U]$ mapping to $\bar z$.  Having such splittings is useful for topological purposes, e.g., $\pi_1^{\et}(\gerb{E}, \bar z) \To \Gal_S$ is naturally a continuous homomorphism of profinite groups when $A$ is finite.

Consider a second geometric base point $\bar z_0 \in \gerb{E}_{\bar s}$ (over the {\em same} $\bar s \To S$).  There exists an isomorphism $\bar f \From \bar z_0 \To \bar z$ in $\gerb{E}_{\bar s}$.  For a sufficiently large Galois cover $U \To S$, we may assume that $\bar f \From \bar z_0 \To \bar z$ arises from a morphism $f \From z_0 \To z$ in $\gerb{E}[U]$.  

Define $\iota_f \From \Aut_\gamma(z_0) \To \Aut_\gamma(z)$ to be the bijection
$$\iota_f(\eta) = \gamma^\ast f \circ \eta \circ f^{-1}, \text{ for all } \eta \in \Aut_\gamma(z_0),$$
making the following diagram commute (in the groupoid $\gerb{E}[U]$).
$$\begin{tikzcd}
z_0 \arrow{r}{f} \arrow{d}{\eta} & z \arrow{d}{\iota_f(\eta)} \\
\gamma^\ast z_0 \arrow{r}{\gamma^\ast f} & \gamma^\ast z 
\end{tikzcd}$$
As $\gamma$ varies over $\Gal_S$, this provides an isomorphism of extensions, $\iota_f \From  \pi_1^{\et}(\gerb{E}, z_0) \To  \pi_1^{\et}(\gerb{E}, z)$.  Passing to the direct limit, we find an isomorphism of extensions depending only on $\bar f \From \bar z_0 \To \bar z$,
$$\begin{tikzcd}
A \inarrow{r} \arrow{d}{=} & \pi_1^{\et}(\gerb{E}, \bar z_0) \onarrow{r} \arrow{d}{\iota_{\bar f}} & \Gal_S \arrow{d}{=} \\
A \inarrow{r} & \pi_1^{\et}(\gerb{E}, \bar z) \onarrow{r} & \Gal_S
\end{tikzcd}$$

Given {\em another} isomorphisms $\bar g \From \bar z_0 \To \bar z$, there exists a unique element $\alpha \in A$ such that $\bar g = \alpha_{\bar z} \circ \bar f$.  As for $f$, we may assume that $\bar g$ arises from $g \From z_0 \To z$ in $\gerb{E}[U]$.  It follows that, for all $\eta \in \Aut_\gamma(z_0)$,
\begin{align*}
\iota_g(\eta) &= \gamma^\ast g \circ \eta \circ g^{-1} \\
&= \gamma^\ast (\alpha_z \circ f) \circ \eta \circ f^{-1} \circ \alpha_z^{-1} \\
&= \gamma^\ast \alpha_z \circ \iota_f(\eta) \circ \alpha_z^{-1} \\
&= \alpha_{\gamma^\ast z} \circ \iota_f(\eta) \circ \alpha_z^{-1}
\end{align*}
In other words, we have $\iota_{\bar g} = \Int(\alpha) \circ \iota_{\bar f}$.

To summarize the relationship between gerbes banded by $\sheaf{A}$ and extensions of $\Gal_S$ by $A$, we have the following.
\begin{thm}
To each geometric base point $\bar z \in \gerb{E}_{\bar s}$, we obtain an extension
$$A \Into \pi_1^{\et}(\gerb{E}, \bar z) \Onto \Gal_S,$$
known as the fundamental group of the gerbe $\gerb{E}$ at $\bar z$.  For any two geometric base points $\bar z_0, \bar z$, we obtain a {\em family} of isomorphisms of extensions
$$\pi_1^{\et}(\gerb{E}, \bar z_0) \To \pi_1^{\et}(\gerb{E}, \bar z),$$
any two of which are related by $\Int(\alpha)$ for a uniquely determined $\alpha \in A$.
\end{thm}

This theorem may seem more natural using 2-categorical language as follows:  consider the 2-category $\Cat{OpExt}(\Gal_S, \sheaf{A})$ whose objects are extensions of $\Gal_S$ by $A$ in which the $\Gal_S$ action on $A = \sheaf{A}_{\bar s}$ coincides with that which arises from the local system $\sheaf{A}$.  The morphisms in this category are homomorphisms of extensions (giving equality on $\Gal_S$ and $A$).  Given two such morphisms $\iota, \iota'$ sharing the same source and target, a natural transformation $\iota \Rightarrow \iota'$ is an element $\alpha \in A$ such that $\iota' = \Int(\alpha) \circ \iota$.  A restatement of the above theorem is the following.
\begin{thm}  
\label{AppxFundGpDefined}
A gerbe $\gerb{E}$ banded by $\sheaf{A}$, and a geometric base point $\bar s \To S$, yield an object $\pi_1^{\et}(\gerb{E}, \bar s)$ of $\Cat{OpExt}(\Gal_S, \sheaf{A})$, well-defined up to equivalence, the equivalence being uniquely determined up to unique natural isomorphism.
\end{thm}

\begin{remark}
If $\sheaf{A}$ is a constant sheaf, then the extension $A \Into \pi_1^{\et}(\gerb{E}, \bar z) \Onto \Gal_S$ is a central extension.  It follows quickly from the theorem that we may define an extension
$$A \Into \pi_1^{\et}(\gerb{E}, \bar s) \Onto \Gal_S$$
up to unique isomorphism (without choice of geometric base point $\bar z$).  In this case, the 2-category $\Cat{OpExt}(\Gal_S, A)$ is an ordinary category:  the only transformations are the identities.
\end{remark}

Next, consider a functor of gerbes $\phi \From \gerb{E}_1 \To \gerb{E}_2$ lying over a homomorphism $c \From \sheaf{A}_1 \To \sheaf{A}_2$ of local systems of abelian groups.  If $\bar z_1$ is a geometric base point for $\gerb{E}_1$ over $\bar s$, arising from $z_1 \in \gerb{E}_1[U]$, then let $z_2 = \phi[U](z_1)$.  Define $\bar z_2$ to be the resulting geometric base point of $\gerb{E}_2$.  Our strictness assumption for functors of gerbes implies that the geometric base point $\bar z_2$ depends only on the geometric base point $\bar z_1$, and not on the choice of $z_1$.

For any $\gamma \in \Gal_S$, we obtain a map $\phi_\gamma \From \Aut_\gamma(z_1) \To \Aut_\gamma(z_2)$, given by
\begin{align*}
\Aut_\gamma(z_1) = \Hom(z_1, \gamma^\ast z_1) & \xrightarrow{\phi[U]} \Hom(\phi[U](z_1), \phi[U](\gamma^\ast z_1) ) \\
& = \Hom( z_2, \gamma^\ast \phi[U](z_1) ) \\
& = \Hom(z_2, \gamma^\ast z_2) = \Aut_\gamma(z_2).
\end{align*}
Putting these together yields a homomorphism of extensions,
$$\begin{tikzcd}
A_1 \inarrow{r} \arrow{d}{c} & \pi_1^{\et}(\gerb{E}_1, \bar z_1) \onarrow{r} \arrow{d}{\phi} & \Gal_S \arrow{d}{=} \\
A_2 \inarrow{r} & \pi_1^{\et}(\gerb{E}_2, \bar z_2) \onarrow{r} & \Gal_S
\end{tikzcd}$$

If $\phi, \phi' \From \gerb{E}_1 \To \gerb{E}_2$ are two functors of gerbes lying over $c \From \sheaf{A}_1 \To \sheaf{A}_2$, and $\phi(\bar z_1) = \phi'(\bar z_1) = \bar z_2$, then we find two such homomorphisms of extensions,
$$\phi, \phi' \From \pi_1^{\et}(\gerb{E}_1, \bar z_1) \To \pi_1^{\et}(\gerb{E}_2, \bar z_2),$$
lying over $c \From A_1 \To A_2$.  

If $N \From \phi \Rightarrow \phi'$ is a natural isomorphism of functors, then $N$ determines an isomorphism $\phi(\bar z_1) \To \phi'(\bar z_1)$, whence an isomorphism $\bar z_2 \To \bar z_2$.  Such an isomorphism is given by an element $\alpha_2 \in A_2 = \Aut(\bar z_2)$, and one may check that
$$\phi' = \Int(\alpha_2) \circ \phi \From \pi_1^{\et}(\gerb{E}_1, \bar z_1) \To \pi_1^{\et}(\gerb{E}_2, \bar z_2).$$

%\bibliography{CovLang2014.bib}
\printbibliography
\end{document}